\documentclass[12pt,a4paper,reqno]{amsart}

\usepackage{amssymb}
\usepackage{amscd}
\usepackage{amsfonts}
\usepackage{setspace}
\usepackage{version}

\numberwithin{equation}{section}

\theoremstyle{plain}

\newtheorem{theorem}[subsection]{Theorem}
\newtheorem{proposition}[subsection]{Proposition}
\newtheorem{lemma}[subsection]{Lemma}
\newtheorem{corollary}[subsection]{Corollary}
\newtheorem{conjecture}[subsection]{Conjecture}

\newtheorem{question}[subsection]{Question}

\newtheorem*{approx-hom-prop-again}{Proposition \ref{approx-hom-prop}}
\newtheorem*{1step-repeat}{Lemma \ref{1-step-approx}}
\newtheorem*{23step-repeat}{Lemma \ref{23-step-approx}}
\newtheorem*{inv-nec-again}{Proposition \ref{inv-nec}}

\theoremstyle{definition}

\newtheorem{definition}[subsection]{Definition}

\renewcommand{\leq}{\leqslant}
\renewcommand{\geq}{\geqslant}

\newsavebox{\proofbox}
\savebox{\proofbox}{\begin{picture}(7,7)%
  \put(0,0){\framebox(7,7){}}\end{picture}}

\newcommand{\md}[1]{\ensuremath{(\operatorname{mod}\, #1)}}
\newcommand{\mdsub}[1]{\ensuremath{(\mbox{\scriptsize mod}\, #1)}}
\newcommand{\mdlem}[1]{\ensuremath{(\mbox{\textup{mod}}\, #1)}}
\newcommand\E{\mathbb{E}}
\newcommand\Z{\mathbb{Z}}
\newcommand\R{\mathbb{R}}

\newcommand\C{\mathbb{C}}
\newcommand\N{\mathbb{N}}

\newcommand\F{\mathbb{F}}

\newcommand\Lip{\operatorname{Lip}}
\newcommand\GI{\operatorname{GI}}

\newcommand\struct{\operatorname{struct}}
\newcommand\ab{\operatorname{ab}}
\newcommand\unstruct{\operatorname{unstruct}}

\newcommand\2{{\bf 2}}

\newcommand\g{{\bf g}}
\newcommand\h{{\bf h}}

\newcommand\eps{\varepsilon}

\newcommand\id{\operatorname{id}}

\def\proof{\noindent\textit{Proof. }}
\def\endproof{\hfill{\usebox{\proofbox}}\vspace{11pt}}

\begin{document}
\title{An inverse theorem for the Gowers $U^4$-norm}
\author{Ben Green}
\address{Centre for Mathematical Sciences\\
Wilberforce Road\\
Cambridge CB3 0WA\\
England }
\email{b.j.green@dpmms.cam.ac.uk}
\author{Terence Tao}
\address{UCLA Department of Mathematics, Los Angeles, CA 90095-1555}
\email{tao@math.ucla.edu}
\author{Tamar Ziegler}
\address{Department of Mathematics \\
Technion - Israel Institute of Technology\\
Haifa, Israel 32000}
\email{tamarzr@tx.technion.ac.il}

\begin{abstract} We prove the so-called \emph{inverse conjecture for the Gowers $U^{s+1}$-norm} in the case $s = 3$ (the cases $s < 3$ being established in previous literature). That is, we show that if $f : [N] \rightarrow \C$ is a function with $|f(n)| \leq 1$ for all $n$ and $\Vert f \Vert_{U^4} \geq \delta$ then there is a bounded complexity $3$-step nilsequence $F(g(n)\Gamma)$ which correlates with $f$. The approach seems to generalise so as to prove the inverse conjecture for $s \geq 4$ as well, and a longer paper will follow concerning this.

By combining the main result of the present paper with several previous results of the first two authors one obtains the generalised Hardy-Littlewood prime-tuples conjecture for any linear system of complexity at most 3. In particular, we have an asymptotic for the number of 5-term arithmetic progressions $p_1 < p_2 < p_3 < p_4 < p_5 \leq N$ of primes.
\end{abstract}

\maketitle
\tableofcontents

\setcounter{tocdepth}{1}

\noindent\textsc{notation}. By a $1$-bounded function on a set $X$ we mean a function $f : X \rightarrow \C$ with $|f(x)| \leq 1$ for all $x \in X$. If the cardinality $|X|$ of $X$ is finite and non-zero, we write $\E_{x \in X} f(x)$ for $|X|^{-1} \sum_{x \in X} f(x)$. Throughout the paper the letter $M$ will refer to a large positive ``complexity'' quantity, normally introduced in each statement of a lemma, proposition or theorem. The letters $c$ and $C$ are reserved for absolute constants with $0 < c < 1 < C$; different instances of the notation will generally denote different absolute constants. If $x \in \R$ we will write $\lfloor x\rfloor$ for the greatest integer less than or equal to $x$, and $\{x\} := x - \lfloor x\rfloor$. If $N$ is a positive integer then we write $[N] := \{1,\dots,N\}$.

\section{Introduction}

This paper concerns a special case of a family of conjectures named the \emph{Inverse Conjectures for the Gowers norms} by  the first two authors. For each integer $s \geq 1$ the inverse conjecture $\GI(s)$, whose statement we recall shortly, describes the structure of $1$-bounded functions $f : [N] \rightarrow \C$ whose $(s+1)$st Gowers norm $\Vert f \Vert_{U^{s+1}}$ is large. These conjectures together with a good deal of motivation and background to them are discussed in \cite{green-icm,green-tao-u3inverse,green-tao-linearprimes}. The conjectures $\GI(1)$ and $\GI(2)$ are already known, the former being straightforward application of Fourier analysis and the latter being the main result of \cite{green-tao-u3inverse}. The aim of the present paper is to establish the first unknown case, that of $\GI(3)$, using what is in essence a method which seems to generalise to prove $\GI(s)$ in general. 

We have taken advantage of some shortcuts and explicit calculations that are specific to the $s = 3$ case, hoping that this will render the paper somewhat appetising as an \emph{hors d'{\oe}uvres} for the general case. The general case will, furthermore, be phrased in the language of non-standard analysis since this provides a very effective framework in which to manage the complicated hierarchies of parameters that appear here. We offer the present paper to those readers who are not immediately comfortable with the nonstandard language; it also serves as an illustration of the point, to be made in the longer paper to follow, that our arguments may be taken out of the choice-dependent realm of nonstandard analysis and, in particular, can lead to effective bounds (albeit extremely weak ones).

We begin by recalling the definition of the Gowers norms. If $G$ is a finite abelian group and if $f : G \rightarrow \C$ is a function then we define
\[ \Vert f \Vert_{U^{k}(G)} := \big(  \E_{x,h_1,\dots,h_k \in G} \Delta_{h_1} \dots \Delta_{h_k}  f(x) \big)^{1/2^k},\]
where $\Delta_h f$ is the multiplicative derivative
\[ \Delta_h f(x) := f(x+h) \overline{f(x)}.\] 
In this paper we will be concerned with functions on $[N]$, which is not quite a group. To define the Gowers norms of a function $f : [N] \rightarrow \C$, set $G := \Z/\tilde N\Z$ for some integer $\tilde N \geq 2^k N$, define a function $\tilde f : G \rightarrow \C$ by $\tilde f(x) = f(x)$ for $x = 1,\dots,N$ and $\tilde f(x) = 0$ otherwise, and set $\Vert f \Vert_{U^k[N]} := \Vert \tilde f \Vert_{U^k(G)}/\Vert 1_{[N]}\Vert_{U^k(G)}$, where $1_{[N]}$ is the indicator function of $[N]$. It is easy to see that this definition is independent of the choice of $\tilde N$, and so for definiteness one could take $\tilde N := 2^k N$. 
Henceforth we shall write simply $\Vert f \Vert_{U^k}$, rather than $\Vert f \Vert_{U^k[N]}$, since all Gowers norms will be on $[N]$. One can show that $\Vert \cdot \Vert_{U^k}$ is indeed a norm for any $k \geq 2$, though we shall not need this here.

The \emph{Inverse conjecture for the Gowers $U^{s+1}$-norm} posits an answer to the following question.

\begin{question}
Suppose that $f : [N] \rightarrow \C$ is a $1$-bounded function and let $\delta > 0$ be a positive real number. What can be said if $\Vert f \Vert_{U^{s+1}} \geq \delta$?
\end{question}

The conjecture made in \cite{green-tao-linearprimes} is that $f$ must correlate with a certain rather algebraic object called an $s$-step nilsequence. In the light of subsequent work \cite{green-tao-nilratner,green-tao-mobiusnilsequences} it seems natural to work with a somewhat more general object called a degree $s$ \emph{polynomial} nilsequence. We recall now the bald definition; for much more motivation and examples, see  the introduction to \cite{green-tao-nilratner}.

\begin{definition}[Polynomial nilsequence]
Let $G$ be a connected, simply-connected nilpotent Lie group. By a \emph{filtration} $G_{\bullet}$ of degree $s$ we mean a nested sequence $G = G_{(0)} = G_{(1)} \supseteq G_{(2)} \supseteq \dots \supseteq G_{(s+1)} = \{\id\}$ with the property that $[G_{(i)},G_{(j)}] \subseteq G_{(i+j)}$. By a polynomial sequence adapted to $G_{\bullet}$ we mean a map $g : \Z \rightarrow G$ such that $\partial_{h_i} \dots \partial_{h_1}g \in G_{(i)}$ for all $h_1,\dots,h_i \in \Z$, where $\partial_h \psi(n) := \psi(n+h)\psi(n)^{-1}$. Let $\Gamma \leq G$ be a discrete and cocompact subgroup, so that the quotient $G/\Gamma$ is a nilmanifold, and assume that each of the $G_{(i)}$ are \emph{rational} subgroups\footnote{One may define rationality topologically, by stipulating that the $G_{(i)}$ are connected Lie subgroups of $G$ and that $\Gamma \cap G_{(i)}$ is a cocompact subgroup of $G_{(i)}$. Some readers may wish to think more concretely, in terms of the existence of a Mal'cev basis as in \cite[Definition 2.1]{green-tao-nilratner}. }. If $F : G/\Gamma \rightarrow \C$ is a $1$-bounded, Lipschitz function then the sequence $(F(g(n)\Gamma))_{n \in \Z}$ is called a polynomial nilsequence of degree $s$.
\end{definition}

\emph{Remark.} An important example of a filtration of a nilpotent group is the \emph{lower central series} $G_0 \supseteq G_1 \supseteq G_2 \supseteq \dots$, in which $G_0 = G_1 = G$, and $G_{i+1} := [G, G_i]$ for $i \geq 1$. It is classical (see, for example, \cite{bourbaki}) that this is a filtration of degree $s$ whenever $G$ is $s$-step nilpotent. This is the minimal example of a filtration, since for any other filtration $G_{\bullet}$ one has $G_i \subseteq G_{(i)}$.\vspace{11pt}

\emph{Remark.} An important fact about polynomial sequences adapted to a filtration $G_{\bullet}$ is that they form a group under pointwise multiplication: see \cite{leibman-group-1} or \cite[Proposition 6.2]{green-tao-nilratner}. A polynomial sequence $g : \Z \rightarrow G$ can also be uniquely expressed as a Taylor expansion $g(n) = g_0^{\binom{n}{0}} g_1^{\binom{n}{1}} \dots g_s^{\binom{n}{s}}$ for some $g_i \in G_{(i)}$ for $i = 0,1,\dots,s$, where $\binom{n}{i}$ is the usual binomial coefficient; see \cite[Section 6]{green-tao-nilratner}.\vspace{11pt}

\emph{Remark.} If $G$ admits a filtration of degree $s$ then, as we remarked above, $G$ must be $s$-step nilpotent. On the other hand, the degree can exceed the step by an arbitrary amount. For instance, if $P : \Z \rightarrow \R/\Z$ is a polynomial of degree $d \geq 1$, then the function $e(P(n)) := e^{\2\pi i P(n)}$ is a polynomial nilsequence of degree $d$, despite being associated to a nilmanifold $G/\Gamma = \R/\Z$ of step just $1$.\vspace{11pt}

Roughly speaking, the inverse conjecture $\GI(s)$ asserts that a $1$-bounded function $f$ has large $U^{s+1}$-norm if and only if it correlates with a degree $s$ nilsequence. However, every aspect of this statement must be quantified in order to make a precise statement. The key issue here lies in defining the \emph{complexity} of a nilsequence, a matter which was addressed in some detail in \cite[Sec 2]{green-tao-nilratner}.  In this paper (fortunately) we can take a much rougher approach. If $\delta > 0$ is some parameter we shall simply say that the complexity of a polynomial nilsequence $(F(g(n)\Gamma))_{n \in \Z}$ is $O_{\delta}(1)$ if the following list of objects are bounded in a way that depends only on $\delta$:
\begin{itemize}
\item $\dim G$;
\item The rationality of some Mal'cev basis $\mathcal{X}$ for $G/\Gamma$ (see \cite[Definition 2.4]{green-tao-nilratner});
\item The rationality of each subgroup $G_{(i)}$ in the filtration (see \cite[Definition 2.5]{green-tao-nilratner});
\item The Lipschitz norm of $F$, measured using the metric defined in \cite[Definition 2.2]{green-tao-nilratner}.
\end{itemize}
We do not encourage the reader to read those definitions in detail at this stage. The important thing to note is that nothing is said about the polynomial sequence $g$, other than that it is adapted to the filtration $G_{\bullet}$.

We may now state the Inverse Conjecture for the Gowers $U^{s+1}$-norm, $\GI(s)$, properly.

\begin{conjecture}[$\GI(s)$]\label{gis-conj}
Suppose that $f : [N] \rightarrow \C$ is a $1$-bounded function and that $\Vert f \Vert_{U^{s+1}} \geq \delta$. Then there is a degree $s$ polynomial nilsequence $(F(g(n)\Gamma))_{n \in \Z}$ of complexity $O_{\delta}(1)$ such that $|\E_{n \in [N]} f(n) \overline{F(g(n)\Gamma)}| \gg_{\delta} 1$.
\end{conjecture}

As hinted earlier, this is not quite the formulation of $\GI(s)$ originally given in \cite[Section 8]{green-tao-linearprimes}. There, it was posited that $f$ correlates with a \emph{linear}\footnote{We remark that a linear nilsequence is not the same thing as a degree $1$ nilsequence; a typical linear nilsequence on an $s$-step nilmanifold will have degree $s$.} nilsequence $(F(g^n\Gamma))_{n \in \Z}$. One might now relabel this the \emph{strong} inverse conjecture. In the longer paper to come we will show how this in fact follows from Conjecture \ref{gis-conj}.  In the special case of the $U^4$-norm under consideration here, it is possible to verify the strong inverse conjecture quite directly by inspection, and we sketch this in Appendix \ref{strong-appendix}. We would, however, like to impress upon the reader our opinion that Conjecture \ref{gis-conj} is the most natural one, a viewpoint that became apparent to the first two authors in the light of our paper \cite{green-tao-nilratner}. Unfortunately \cite{green-tao-linearprimes} was written before that paper and hence operates under the assumption of the strong inverse conjecture. Relatively simple changes would be required to make all of the arguments there work under the assumption of Conjecture \ref{gis-conj} however, the key issue being \S 11 of that paper.  

The evidence for the inverse conjectures prior to the present work was a ``local version'' due to Gowers \cite{gowers-longaps}, its truth in the cases $s = 1$ and $2$ (see \cite{green-tao-u3inverse}) as well as the truth of analogues of the conjecture in both ergodic theory \cite{host-kra, ziegler} and in the ``finite field model'' in which $[N]$ is replaced by $\F^n$ for some small prime field $\F$ \cite{bergelson-tao-ziegler,tao-ziegler}. 

It is also known that this conjecture is necessary, in the following sense. 

\begin{proposition}[Necessity of inverse conjecture]\label{inv-nec}
Suppose that $f :[N] \rightarrow \C$ is a $1$-bounded function, that $(F(g(n)\Gamma))_{n \in \Z}$ is a polynomial nilsequence of degree $s$ and complexity $O_{\delta}(1)$, and that $|\E_{n \in [N]} f(n) \overline{F(g(n)\Gamma)}| \geq \delta$. Then $\Vert f \Vert_{U^{s+1}} \gg_{\delta} 1$.
\end{proposition}

There is currently no proof of this written in the literature. In the case of \emph{linear nilsequences} $(F(g^n\Gamma))_{n \in \Z}$ there are two different (albeit related) proofs in the literature: one in \cite[Proposition 12.6]{green-tao-u3inverse} and the other in \cite[Section 11]{green-tao-linearprimes}. The second of these proofs would generalise rather easily to the more general setting of degree $s$ polynomial nilsequences $(F(g(n)\Gamma))_{n \in \Z}$, the key issue being to note that \cite[Lemma E.4]{green-tao-linearprimes} is true for the values $(g(n + \omega \cdot h)\Gamma)_{\omega \in \{0,1\}^{s+1}}$, this being essentially \cite[Proposition 6.5]{green-tao-nilratner}. The reader will doubtless be relieved to hear that we recently discovered a very short proof of Proposition \ref{inv-nec}, and we give this in Appendix \ref{nec-appendix}. Note, however, that this proposition is included for motivation and interest only, and is not actually required in this paper.

Here, then is the main result of our paper.

\begin{theorem}[$\GI(3)$]\label{mainthm}
The inverse conjecture for the $U^4$-norm, $\GI(3)$, is true.
\end{theorem}

As already remarked, in Appendix \ref{strong-appendix} we will also establish the strong form of the inverse conjecture for the $U^4$-norm, in the form given in \cite[Section 8]{green-tao-linearprimes}.

By combining this result with the previous results in \cite{green-tao-linearprimes, green-tao-mobiusnilsequences} we obtain a proof of what was referred to in \cite{green-tao-linearprimes} as the generalised Hardy-Littlewood conjecture for linear systems of complexity at most $3$. In particular we have the following.

\begin{theorem}
The number of quintuples of primes $p_1 < p_2 < p_3 < p_4 < p_5 \leq N$ in arithmetic progression is asymptotic to $\gamma N^2/\log^5 N$, where 
\[ \gamma = \frac{27}{16} \prod_{p \geq 5} \frac{p^3(p-4)}{(p-1)^4}.\]
\end{theorem}

We refer the reader to \cite{green-tao-linearprimes} for further discussion. Several further applications of the $\GI(s)$ conjectures will be given in a forthcoming paper of the first two authors \cite{green-tao-arithmetic-regularity}.\vspace{11pt}

\emph{Acknowledgements.} BG was, for some of the period during which this work was carried out, a fellow of the Radcliffe Institute at Harvard. He is very grateful to the Radcliffe Institute for providing excellent working conditions. TT is supported by NSF Research Award DMS-0649473, the NSF Waterman award and a grant from the MacArthur Foundation. TZ is supported by ISF grant 557/08, an Alon fellowship, and a Landau fellowship of the Taub foundation . All three authors are very grateful to the University of Verona for allowing them to use classrooms at Canazei during a week in July 2009 when this work was largely completed.

\section{An outline of the proof}\label{inductive-sec}

In this section we outline the argument we use to establish the inverse conjecture for the $U^4$-norm.

It is easy to show, and well-known, that if $\Vert f \Vert_{U^4} \geq \delta$ then there are $\gg \delta^C N$ values of $h$ for which $\Delta_h f(n) := f(n+h)\overline{f(n)}$ has $U^3$-norm at least $\delta^C$. Applying $\GI(2)$, it follows that for all these $h$ we have
\begin{equation}\label{derivative-assumption} |\E_n \Delta_h f(n)\overline{\chi_h(n)}| \gg 1,\end{equation} where $\chi_h(n)$ is a $2$-step nilsequence (with complexity bounded uniformly in $h$). 

Very roughly speaking, the aim is to show that these $2$-step nilsequences ``line up'' in such a way that they may be interpreted as the derivatives of a single $3$-step object. To make this work and for ease of exposition it is convenient to assume that $\chi_h(n)$ is in fact equal to $e(\psi_h(n))$, where $\psi_h(n)$ is a bracket quadratic phase: a sum of terms of the form $\alpha_1 n\lfloor \alpha_2 n\rfloor$, $\alpha_3 n^2$ and $\alpha_4 n$. The link between these objects and 2-step nilsequences was explored in \cite[Section 10]{green-tao-u3inverse} and will be recalled later in this paper. For the purposes of this discussion let us suppose that $\psi_h(n) = \alpha_h n \lfloor \beta_h n\rfloor $; this is something of a simplification of the true situation.

%\begin{equation}\label{chi-form} \chi_h(n) =  \prod_{1 \leq i < i' \leq k} F_{[i',i]}^{m_{[i',i]}(h)}(g_h(n)\Gamma) \cdot e(\alpha_h n^2 + \beta_h n)\end{equation} for integers $m_{[i',i]}(h) = O_{\delta}(1)$. 

Here is a rough outline of the main steps we shall be taking to control the dependence of $\alpha_h$ and $\beta_h$ on $h$. Suppose that \eqref{derivative-assumption} holds with $\chi_h(n) = e(\alpha_h n \lfloor \beta_h n\rfloor)$.

\begin{enumerate}
\item[\emph{Step 1}] (Reducing the $h$-dependence) We may assume (possibly after refining the set of $h$ and modifying $\alpha_h$ and $\beta_h$ somewhat) that $\beta_h$ does not depend on $h$. 

\item[\emph{Step 2}](Approximate linearity of $h$-dependent frequency) We may assume (possibly after refining the set of $h$ again) that $\alpha_h$ is approximately equal to a bracket linear form $\theta_1\{\theta'_1 h\} + \dots + \theta_d \{\theta'_d h\} + \theta h$.

\item[\emph{Step 3}](Symmetry argument) Following \emph{Step 2}, $\chi_h(n)$ is essentially $e(\psi_h(n))$ with the phase $\psi_h(n)$ being of the form $T(h,n,n)$, where $T(n_1,n_2,n_3)$ is a sum of terms of the form $\{\theta_1 n_1\}\theta_2 n_2\lfloor \theta_3 n_3\rfloor$. Not every such function $T(h,n,n)$ can be obtained as the ``derivative'' of a 3-step object, however, and in order to make this assertion we need some additional symmetry properties of the ``generalised trilinear form'' $T(n_1,n_2,n_3)$. 
\end{enumerate}

It may be of some interest to make a comparison between this strategy and that used in the proof of the $U^3$-inverse theorem \cite{green-tao-u3inverse}. If $\Vert f \Vert_{U^3} \geq \delta$ then for many $h$ we have, once again, 
\[ |\E_n \Delta_h f(n)\overline{\chi_h(n)} | \gg 1,\] but now $\chi_h(n)$ may be assumed to be nothing more complicated than a linear phase $e(\alpha_h n)$. The argument runs roughly as follows:

\begin{enumerate}
\item[\emph{Step 2}'](Approximate linearity of frequencies) At the possible expense of passing to a subset of the $h$, the frequencies $\alpha_h$ are approximately ``bracket-linear'' in $h$, as above;
\item[\emph{Step 3}'](Symmetry argument) Following \emph{Step 2}, $\chi_h(n)$ is essentially $e(\psi_h(n))$ with the phase $\psi_h(n)$ being of the form $T(h,n)$ where $T(n_1,n_2)$ is a sum of terms of the form $\{\theta_1 n_1\}\theta_2 n_2$. Not every such function $T(h,n)$ can be obtained as the ``derivative'' of a $2$-step object, however, and in order to make this assertion we need some additional symmetry properties of the form $T(n_1,n_2)$.
\end{enumerate}

We note that \emph{Step 2'} is essentially due to Gowers \cite[Chapter 7]{gowers-longaps}, although one must apply a little extra geometry of numbers to get the precise conclusion we hint at here. \emph{Step 3'} is due to the first two authors and is the main new result of \cite{green-tao-u3inverse}, specifically Lemma 9.4 of that paper. Note that \emph{Step 1} in the outline above did not feature at all in the proof of the $U^3$-inverse theorem and it is new to this paper.

Let us say a few words about how \emph{Steps 1, 2} and \emph{3} are accomplished. The key to almost all of our analysis is a straightforward adaption of a fundamental idea of Gowers \cite{gowers-4aps}, which proceeds from the assumption that 
\begin{equation}\label{derivative-correlation} |\E_n \Delta_h f(n)\overline{\chi_h(n)}| \gg 1\end{equation} for many $h$ and draws a conclusion involving just the $\chi_h(n)$, and not the function $f$. This argument is valid for any bounded functions $\chi_h(n)$ and we give it in \S \ref{gowers-sec}. 

The conclusion of that argument is that 
\begin{equation}\label{gowers-conclusion} \E_{n \in [N]} \chi_{h_1}(n)\chi_{h_2}(n + h_1 - h_4)\overline{\chi_{h_3}(n) \chi_{h_4}(n + h_1 - h_4)} \gg 1\end{equation} for many additive quadruples $h_1,h_2,h_3,h_4$, that is to say quadruples satisfying $h_1 + h_2 = h_3 + h_4$.

\emph{Steps 1,2} and \emph{3} all involve interpreting this in the case that $\chi_h(n)$ is a 2-step object such as a bracket quadratic phase. One way to do this is to visualise 
\[ \chi_{h_1}(n)\chi_{h_2}(n + h_1 - h_4)\overline{\chi_{h_3}(n) \chi_{h_4}(n + h_1 - h_4)}\] as a certain $h_1,h_2,h_3,h_4$-dependent nilsequence on a product of four nilmanifolds (one for each of the $h_i$), in which case \eqref{gowers-conclusion} states that the underlying polynomial sequence $g_{h_1,h_2,h_3,h_4}(n)$ is far from equidistributed. This situation may then be studied using the distributional results on nilsequences contained in \cite{green-tao-nilratner} in order to draw conclusions concerning the dependence on $h$ of ``leading order'' terms in the $\chi_h(n)$.

\emph{Steps 1} and \emph{2} really only use the ``top-order'' structure of \eqref{gowers-conclusion} -- that is to say the shifts $h_1 - h_4$ are not relevant. To handle \emph{Step 3} these shifts cannot be ignored. In the general case the treatment of \emph{Step 3} will involve another appeal to the distributional results on nilmanifolds in \cite{green-tao-nilratner}, but in the case of the $U^4$-norm a much more hands-on approach involving Bohr sets may be employed, and it is this argument that we give here.

The following deliberately vague discussion may perhaps be helpful. Suppose that $|f(n)| = 1$ for all $n$ and that $\Delta_h f \sim \chi_h$ (where we are not attaching any real meaning to $\sim$). Then we have the ``cocycle identity'' $\Delta_{h+k}f(n) = \Delta_h f(n+k) \Delta_k f(n)$, which translates to $\chi_{h + k}(n) \sim \chi_h(n + k) \chi_k(n)$. Imagining that the shift $n \mapsto n + k$ does not affect the ``top-order structure'' of $\chi_h(n+k)$, we have the approximate linearity condition
\[ \chi_{h+k} \sim \chi_h\chi_k \qquad \mbox{to top order}.\]
Roughly speaking, \emph{Steps 1} and \emph{2} are concerned with exploiting this rigorously.
On the other hand we also have the symmetry relation $\Delta_h \Delta_k f = \Delta_k \Delta_h f$, which suggests that $\Delta_h \chi_k \sim \Delta_k \chi_h$; \emph{Step 3} may be thought of in terms of exploiting this kind of information.

\section{Almost nilsequences}

In this paper we will be dealing with various objects which are ``almost'' nilsequences but not quite. They can invariably be represented as $F(g(n)\Gamma)$ for some function $F$ which is only \emph{piecewise} Lipschitz, the discontinuities being on sets which are somehow ``polynomial''. Rather than formalise these notions, we instead introduce the notion of an approximate nilsequence, give some examples, and point out a number of consequences of the definition.

\begin{definition}[Almost nilsequences]\label{almost-nil-def}
Suppose that $\Psi : [N] \rightarrow \C$ is a $1$-bounded function and that $M > 1$ is a complexity parameter. Then we say that $\Psi$ is a \emph{degree $s$ almost polynomial nilsequence} of complexity $O_M(1)$ if, for any $\eps > 0$, there is a genuine degree $s$ polynomial nilsequence $\Psi_{\eps}$ with complexity $O_{s,\eps,M} (1)$ such that $\E_{n \in [N]}|\Psi(n) - \Psi_{\eps}(n)| \leq \eps$. 
\end{definition}

\emph{Remarks.} That is, $\Psi$ can be approximated arbitrarily well, in $L^1$, by genuine nilsequences. We will not specify the function $O_{s,\eps,M}(1)$ exactly (and indeed it does not make sense to do so, in view of the loose manner in which we have defined complexity). The reader should just imagine that there is \emph{some} fixed function which may be taken in this definition and which makes all statements that we make later on true. Let us also remark that the non-standard analogue of this definition, which will feature in our forthcoming paper on the general case $\GI(s)$, is much cleaner and does not involve any unspecified complexity parameters $O_{s,\eps,M}(1)$.

We make the following easily verified, but rather useful, claim:

\begin{lemma}[Algebra properties]\label{alg-lem}  If $\Phi, \Psi$ are degree $s$ almost polynomial nilsequences, then their sum $\Phi+\Psi$ and product $\Phi \Psi$, and complex conjugate $\overline{\Phi}$ are also degree $s$ almost polynomial nilsequences \textup{(}with a slightly different complexity bound $O_{s,\eps,M}(1)$ on the approximants, of course\textup{)}.
\end{lemma}

The utility of Definition \ref{almost-nil-def} is made clear by the following lemma, which states that correlation with almost nilsequences is essentially the same thing as correlation with genuine nilsequences.

\begin{lemma}\label{almost-to-genuine}
Suppose that $f : [N] \rightarrow \C$ is a 1-bounded function and that \[|\E_{n \in [N]} f(n) \Psi(n)| \geq \delta\] for some degree $s$ almost polynomial nilsequence $\Psi$ of complexity $O_{\delta}(1)$. Then there is a genuine degree $s$ polynomial nilsequence $F(g(n)\Gamma)$ of complexity $O_{s,\delta}(1)$ such that $|\E_{n \in [N]} f(n) F(g(n)\Gamma)| \geq \delta/2$.
\end{lemma}
\proof Simply take $\eps = \delta/2$ in Definition \ref{almost-nil-def} and set $F(g(n)\Gamma) = \Psi_{\eps}(n)$.\endproof

A particular consequence, which we shall make use of later, is that it suffices to establish Conjecture \ref{gis-conj} with almost nilsequences instead of genuine ones. 

For 1-step nilsequences there is a further, very helpful, reduction that can be made.

\begin{lemma}[$1$-step correlation]\label{1-step-to-phase}
Suppose that $f : [N] \rightarrow \C$  is a 1-bounded function and that $|\E_{n \in [N]} f(n)\Psi(n)| \geq \delta$ for some degree $1$ almost nilsequence $\Psi$ of complexity $O_{\delta}(1)$. Then there is a $\theta \in \R/\Z$ such that $|\E_{n \in [N]} f(n) e(\theta n)| \gg_{\delta} 1$.  \textup{(}The implied constants here depend of course on the implied constants in the definition of an almost nilsequence.\textup{)}
\end{lemma}
\proof By the previous lemma we may assume that $\Psi$ is a genuine degree $1$ nilsequence of complexity $O_{\delta}(1)$, that is to say a sequence of the form $(F(n \alpha))_{n \in \Z}$ where $\alpha \in (\R/\Z)^k$ for some $k = O_{\delta}(1)$ and $F : (\R/\Z)^k \rightarrow \C$ is a function with Lipschitz constant $O_{\delta}(1)$. Standard Fourier analysis (see, for example, \cite[Lemma A.9]{green-tao-u3mobius}) implies that we may expand 
\[ 
F(t) = \sum_{m \in \Z^k, |m| \leq O_{\delta}(1)} c_m e(m \cdot t) + O(\delta/10),
\] where the $c_m$ are complex numbers with $|c_m| = O_{\delta}(1)$. The result follows quickly from this.\endproof

The next two lemmas collect together various examples of almost nilsequences. The proofs, which are somewhat technical and tedious, are given in Appendix \ref{approx-nil-app}.

\begin{lemma} \label{1-step-approx}Suppose that $\alpha,\beta \in [0,1]$ and that $M > 1$ is a complexity parameter. The following are all examples of almost nilsequences of degree $1$ and complexity $O_M(1)$:
\begin{enumerate}
\item  the set of $1$-step Lipschitz nilsequences of complexity at most $M$;
\item  the set of characteristic functions $1_P$, where $P \subseteq [N]$ is a progression of length at least $N/M$;
\item  the set of functions of the form $n \mapsto e(\alpha \{\beta n\})$, with $\alpha \in \R$ and $\beta \in \R/\Z$;
\item  the set of functions of the form $n \mapsto e(\{\alpha n\} \{\beta n\})$, with $\alpha,\beta \in \R/\Z$;
\item  the set of functions of the form $n \mapsto e(\alpha n\lfloor\beta n\rfloor)$, where $\Vert \beta \Vert_{\R/\Z} \leq M/N$.
\end{enumerate}
In particular (by Lemma \ref{1-step-to-phase}), if $f: [N] \to \C$ is a $1$-bounded function such that $|\E_{n \in [N]} f(n) \Psi(n)| \geq \delta$, where $\Psi$ is one of the functions on the above list, then there exists $\theta \in \R/\Z$ such that $|\E_{n \in [N]} f(n) e(\theta n)| \gg_{\delta,M} 1$.
\end{lemma}

\begin{lemma}\label{23-step-approx} Suppose that $\alpha,\beta,\gamma \in [0,1]$. Then the following are all examples of almost nilsequences of degree $s \geq 2$ and complexity $O(1)$:
\begin{enumerate}
\item $n \mapsto e(\{\alpha n\}\beta n)$, of degree $2$;
\item $n \mapsto e(\{\alpha n\}\beta n^2)$, of degree $3$;
\item $n \mapsto e(\{\alpha n\} \{\beta n \} \gamma n)$, of degree $3$.
\end{enumerate}
\end{lemma}

Although the proof of this last lemma is little tedious, it is also important in the sense that this is the only place in our paper where a 3-step nilsequence is actually constructed.

\section{Distributional results concerning nilsequences}\label{nil-dist-sec}

We will rely heavily on the quantitative distribution results concerning polynomial nilsequences $(g(n)\Gamma)_{n \in [N]}$ established by the first two authors in \cite{green-tao-nilratner}. There were two main results in that paper, the first of which was used in the proof of the second. We are aware that the paper \cite{green-tao-nilratner} is long and somewhat difficult. However, the reader wishing to understand the present paper need only be \emph{au fait} with the statements of the results there, which means that she need only read Chapters 1 and 2 of the paper. We will assume familiarity with those chapters throughout this paper, and in particular will use notation from them without further comment.  We will also revisit these results in a non-standard setting in the sequel to this paper, in which we will give more detailed proofs.

The first result we refer to gives a criterion for $(g(n)\Gamma)_{n \in [N]}$ being equidistributed. This is \cite[Theorem 2.9]{green-tao-nilratner}. This theorem is a quantitative version of a polynomial equidistribution theorem for nilmanifolds. The qualitative version basically claims
that equidistribution of polynomial sequences is determined on the  abelianization $G/[G,G]\Gamma$.
For linear sequences this is a classical result, and for polynomial sequences the result is due to Leibman \cite{leibman-poly}. 

\begin{theorem}[Quantitative Leibman dichotomy]\label{quant-leib}\cite[Theorem 2.9]{green-tao-nilratner}
Let $m,s \geq 0$, $0 < \delta < 1/2$ and $N \geq 1$. Suppose that $G/\Gamma$ is an $m$-dimensional nilmanifold together with a filtration $G_{\bullet}$ of degree $s$ and that $\mathcal{X}$ is a $\frac{1}{\delta}$-rational Mal'cev basis adapted to $G_{\bullet}$. Suppose that $g : \Z \rightarrow G$ is a polynomial sequence adapted to $G_{\bullet}$. If $(g(n)\Gamma)_{n \in [N]}$ is not $\delta$-equidistributed, then there is some $k \in \Z^{m'}$, where $m' := \dim G - \dim [G,G]$ is the dimension of the horizontal torus of $G/\Gamma$, with $|k| \ll \delta^{-O_{m,s}(1)}$ and 
\begin{equation}\label{kpig}
\Vert k \cdot (\pi \circ g)\Vert_{C^{\infty}[N]} \ll \delta^{-O_{m,s}(1)}, 
\end{equation}
where $\pi : G \rightarrow G/[G,G]\Gamma \equiv (\R/\Z)^{m'}$ is projection onto the horizontal torus of $G/\Gamma$.
\end{theorem}

The second result we allude to, proved in sections 9 and 10 of \cite{green-tao-nilratner} by iterating the preceding theorem, is a certain factorization result. We will need a variant of it in the present paper involving an arbitrary growth function $\omega : \R^+ \rightarrow \R^+$; this may be established\footnote{We feel rather sorry for our readers at this point. One particular advantage of the non-standard analysis approach to be taken in the more general paper to follow is that the need for arbitrary growth functions $\omega$ is eliminated.} by exactly the same iterative argument that is used in the proof of \cite[Theorem 10.2]{green-tao-nilratner}.

\begin{theorem}[Factorization result]\label{nil-distribution}
Let $s,N \geq 0$ be integers, let $M \geq 1$ be a real number, and let $\omega : \R^+ \rightarrow \R^+$ be an arbitrary growth function. Suppose that $G/\Gamma$ is a nilmanifold of complexity at most $M$ together with a filtration $G_{\bullet}$ of degree $s$. Suppose that $\mathcal{X}$ is an $M$-rational Mal'cev basis adapted to $G_{\bullet}$ and that $g : \Z \rightarrow G_{\bullet}$ is a polynomial map adapted to $G_{\bullet}$. Then there is an integer $M_0$ with $M \leq M_0 = O_{M,s,\omega}(1)$, a rational subgroup $G' \subseteq G$, a Mal'cev basis $\mathcal{X}'$ for $G'/\Gamma'$ in which each element is an $M_0$-rational combination of the elements of $\mathcal{X}$, and a decomposition $g = \eps g' \gamma$ into polynomial sequences $\eps,g',\gamma : \Z \rightarrow G$ adapted to $G_{\bullet}$ with the following properties:
\begin{enumerate}
\item $\eps : \Z \rightarrow G$ is $(M_0,N)$-smooth;
\item $g' : \Z \rightarrow G'$ takes values in $G'$, and the finite sequence $(g'(n)\Gamma'')_{n \in [N]}$ is totally $1/\omega(M_0)$-equidistributed in $G'/\Gamma''$, whenever $\Gamma''$ is a sublattice of $\Gamma'$ of index at most $\omega(M_0)$, and using the metric $d_{\mathcal{X}}$ on $G'/\Gamma''$;
\item $\gamma : \Z \rightarrow G$ is $M_0$-rational, and $(\gamma(n)\Gamma)_{n \in \Z}$ is periodic with period at most $M_0$.
\end{enumerate}
\end{theorem}

\emph{Remark.} The terms ``smooth'' and ``totally equidistributed'' in this sort of context will not feature elsewhere in the paper, as we shall rely only on this theorem to prove Theorem \ref{quant-ratner} below.  An extremely similar deduction was utilised (and proved in some detail) in \S 2 of \cite{green-tao-mobiusnilsequences}.\vspace{11pt}

\emph{Sketch proof.}  The main idea is to iterate Theorem \ref{quant-leib} using the following ``dimension reduction argument''.  At any given stage of the argument, one has an initial factorisation $g = \eps g' \gamma$ obeying all the properties claimed in the theorem for some $M_0$, except for the equidistribution conclusions on $g'$.  (Note that one can trivially obtain such an initial factorisation by setting $\eps$ and $\gamma$ to be the identity, and $G' = G$.)  If $g'$ obeys the stated equidistribution properties, then we are done.  Otherwise, by appealing to Theorem \ref{quant-leib} and refining to a finite sublattice of $\Gamma'$ if necessary, the horizontal coefficients of $g'(n)$ will contain an approximate linear dependence in the sense of \eqref{kpig}.  One can then use this, following the arguments used to prove \cite[Theorem 10.2]{green-tao-nilratner}, in order to factorise $g' = \eps' g'' \gamma'$, where $\eps', \gamma'$ satisfy similar properties to $\eps, \gamma$ but with a worse value of $M_0$, and $g''$ takes values in a connected subgroup $G''$ of $G'$ of strictly lower dimension.  We then absorb the $\eps'$ and $\gamma'$ factors to $\eps, \gamma$, replace $G'$ by $G''$, increase $M_0$ to a larger quantity depending on $M_0$ and $\omega$, and continue the argument.  Since one cannot have an infinite descent of connected subgroups of $G$, the argument must eventually terminate with a factorisation with the desired properties.\endproof

The theorem below is  a quantitative version of an equidistribution result of Leibman \cite{leibman-poly}
stating that the orbit closure of a polynomial sequence $(g(n)\Gamma)_{n \in \N}$ is a finite union of subnilmanifolds $Y_j$, each a closed orbit of  a connected closed subgroup $H_j$ of $G$; moreover the polynomial sequence visits each $Y_j$ periodically, and is well distributed there with respect to the normalized Haar measure.  Much the same argument (with more details) is given in Section 2 of \cite{green-tao-mobiusnilsequences}.

\begin{theorem}[``Quantitative Ratner'' result]\label{quant-ratner}
Let $s, N \geq 0$ be integers, let $M \geq 1$ be a real number, and let $\omega : \R^+ \rightarrow \R^+$ be an arbitrary growth function. Suppose that $G/\Gamma$ is a nilmanifold of complexity at most $M$ together with a filtration $G_{\bullet}$ of degree $s$. Suppose that $\mathcal{X}$ is an $M$-rational Mal'cev basis adapted to $G_{\bullet}$ and that $g : \Z \rightarrow G_{\bullet}$ is a polynomial map adapted to $G_{\bullet}$. Then there is an integer $M_0$ with $M \leq M_0 = O_{M,s,\omega}(1)$ and a decomposition of $[N]$ into subprogressions $P_j$, each of length at least $N/M_0$, together with $M_0$-rational connected subgroups $H_j \leq G$ and elements $x_j \in G$ with coordinates at most $M_0$ such that $(g(n)\Gamma)_{n \in P_j}$ is $1/\omega(M_0)$-equidistributed on $x_j H_j\Gamma/\Gamma$ for each $j$.
\end{theorem}

\emph{Sketch proof.} In Theorem \ref{nil-distribution}, take a growth function $\omega' : \R^+ \rightarrow \R^+$ even more rapidly growing than the $\omega$ in the statement here. Let $g = \eps g'\gamma$ be the resulting decomposition. Take the progressions $P_j$ to have common difference $q$, the period of $\gamma(n)\Gamma$, and length sufficiently small that the smooth term $\eps(n)$ is almost constant on each $P_j$. Choose $y_j,\gamma_j$ such that $\eps(n) \approx y_j$ and $\gamma(n)\Gamma = \gamma_j\Gamma$ for $n \in P_j$. Then the theorem holds with $H_j := \gamma_j^{-1} G' \gamma_j$ and $x_j := y_j\gamma_j$. Note that the action of conjugation by $\gamma_j$ moves $\Gamma$ to a slightly different subgroup of $G$, but this new group intersects $\Gamma$ in a subgroup of index $O_{M_0}(1)$, and so one can proceed by using the fact that $g'$ is assumed equidistributed with respect to such subgroups also.  
\endproof

\section{Free nilpotent Lie groups and free nilcharacters}\label{sec3}

In previous papers in additive combinatorics in which nilsequences have been discussed, such as \cite{green-tao-u3inverse,green-tao-nilratner}, the Heisenberg nilmanifold has been the central example and readers have been encouraged to think of upper triangular matrix groups as the archetypal nilpotent Lie groups. A key innovation in this paper and the sequel \cite{uk-inverse}, strongly inspired by the recent work of Leibman on bracket polynomials \cite{leibman}, is a shift away from this viewpoint. Instead, it seems that \emph{free} nilpotent Lie groups and certain functions on them play a crucial r\^ole. 

In this section we give some basic definitions in this regard in the 2-step case. In Appendix \ref{approx-nil-app} we will briefly meet  an example of the 3-step case, but for the most part we will be working with 2-step objects in which case it is not a particularly onerous task to proceed very explicitly. The definitions in the higher step case are similar but necessarily require some more general discussion of bases in free nilpotent Lie algebras. 

\begin{definition}[Free 2-step nilpotent Lie group and nilmanifold]\label{free-2-step}
By the free $2$-step nilpotent Lie group on generators $e_1,\dots,e_k$ we mean 
\[ G :=\{e_1^{t_1} \dots e_k^{t_k} e_{[2,1]}^{t_{[2,1]}} \dots  e_{[k,k-1]}^{t_{[k,k-1]}} : t_1,\dots,t_k, t_{[2,1]},\dots, t_{[k,k-1]} \in \R\},\]
subject to the relations $e_i^{-1}e_j^{-1} e_ie_j = [e_i, e_j] = e_{[i,j]}$ for $1 \leq j < i \leq k$. By the \emph{standard filtration} $G_{\bullet}$ we mean simply the lower central series filtration with $G_{(0)} = G_{(1)} = G$, $G_{(2)} = [G,G]$ and $G_{(3)} = \{\id\}$.
Inside $G$ we take the standard lattice 
\[ \Gamma := \{e_1^{m_1} \dots e_k^{m_k} e_{[2,1]}^{m_{[2,1]}} \dots  e_{[k,k-1]}^{m_{[k,k-1]}} : m_1,\dots,m_k, m_{[2,1]},\dots, m_{[k,k-1]} \in \Z\}.\]
The quotent $G/\Gamma$ is then called the free $2$-step nilmanifold on $k$ generators.
\end{definition}

A Mal'cev basis for $G/\Gamma$ consists of the elements $X_i = \log e_i$ and $X_{[i',i]} = \log e_{[i',i]}$; the Mal'cev coordinates of an element of $G$ are simply the elements \[ (t_1,\dots,t_k,t_{[2,1]},\dots,t_{[k,k-1]}).\] As in \cite{green-tao-nilratner}, such a basis may be used to coordinatise $G/\Gamma$ by identifying $[0,1]^{k + \binom{k}{2}}$ as a fundamental domain for the right action of $\Gamma$ on $G$. 
Let us perform a calculation. In Mal'cev coordinates it is easy to check that the multiplication law on $G$ corresponds to the operation
\[ (t_i, t_{[i',i]})  \ast (u_i, u_{[i',i]}) = (t_i + u_i, t_{[i',i]} + u_{[i',i]} + t_{i'} u_i).\] For a given element $g  \in G$ with coordinates $(t_i,t_{[i',i]})$ we may pick some $\gamma  \in \Gamma$ such that $g\gamma$ has coordinates in the fundamental domain $\mathcal{F} = [0,1]^{k + \binom{k}{2}} \in \R^{k + \binom{k}{2}}$. Possible coordinates for $\gamma$ are 
\[u_i = -[t_i ],u_{[i',i]} = -[  t_{[i',i]}- t_{i'}[t_i]],\]  where $[\; ]$ is the floor function. These coordinates are unique if $g\gamma$ lies in the interior of the fundamental domain $\mathcal{F}$. The coordinates of $g\gamma$ are then
\[  ( \{ t_i \},  \{ t_{[i',i]}- t_{i'}[t_i]\}).\] 
\begin{definition}[Coordinates]
Suppose that $G/\Gamma$ is the free $2$-step nilmanifold on $k$ generators. Suppose that an element $g \in G$ has Mal'cev coordinates $t_i, t_{[i',i]}$. Then the coordinates of $g\Gamma \in G/\Gamma$ are the entries of the vector
\[ \big( \{ t_i \}, \{t_{[i',i]}- t_{i'}[t_i]\}  \big).\]
We write them as $(t_{i},t_{[i',i]})$.
\end{definition}

\begin{definition}[Coordinate functions]\label{def5.3}
By the \emph{basic coordinate functions} $F_{i}, F_{[i',i]} : G/\Gamma \rightarrow \C$ we mean the functions $F_{i}(t) = e(t_{i})$ and $F_{[i',i]}(t) = e(t_{[i',i]})$. The \emph{top order} basic coordinate functions $F_{[i',i]}$ will have a particularly important role to play.
\end{definition}

Consider now a polynomial sequence $g : \Z \rightarrow G$ of the form
\begin{equation}\label{2-step-poly-seq}
g(n) := (\xi_1 n, \dots, \xi_k n, q_{[2,1]}(n), \dots, q_{[k,k-1]}(n)),
\end{equation}
where the $q_{[i',i]}$ are quadratic polynomials. By the theory developed towards the end of \S 6 of \cite{green-tao-nilratner} (or simply by a short direct calculation), these are degree two polynomial sequences adapted to the standard filtration $G_{\bullet}$ based on the lower central series. The objects $F_{i}(g(n)\Gamma),F_{[i',i]}(g(n)\Gamma)$ are then called \emph{free $2$-step nilcharacters}, and they will be basic building blocks in this paper. The top-order nilcharacters involving $F_{[i',i]}$ will play a particularly crucial role. In the light of the above computations these top-order free $2$-step nilcharacters may be computed quite explicitly, and indeed we have
\begin{equation}\label{5.33} F_{[i',i]}(g(n)\Gamma) = e( \xi_i n\lfloor \xi_{i'} n\rfloor)e(\alpha_{[i',i]}n^2 + \beta_{[i',i]}n),\end{equation} for some $\alpha_{[i',i]}, \beta_{[i',i]} \in \R$. By altering the quadratics $q_{[i',i]}(n)$ we may make the coefficients $\alpha_{[i',i]}, \beta_{[i',i]}$ arbitrary. These quadratic phases $e(\alpha n^2 + \beta n)$ should be thought of as essentially 1-step objects, albeit of degree 2, and the most important feature of our 2-step nilcharacters are the bracket monomials $\xi_{i'} n\lfloor \xi_{i} n\rfloor$. We will often use explicit bracket-quadratics in this paper. In the longer paper to come, dealing with the general case, it will not be possible to proceed so explicitly and indeed the main new innovation of that paper (following the work of Leibman) is to develop a kind of ``calculus'' of bracket polynomials.

Let us note that $F_{[i',i]}(g(n)\Gamma)$ is not actually a $2$-step nilsequence, because the function $F_{[i',i]}$ is only piecewise Lipschitz. From the explicit form given above and Lemma \ref{23-step-approx}, however, one sees that it \emph{is} an almost $2$-step nilsequence.

We now give a variant of the $U^3$ inverse theorem involving 2-step free nilcharacters.

\begin{theorem}[Inverse theorem for $U^3$, variant]\label{u3-variant}
Suppose that $f : [N] \rightarrow \C$ is a 1-bounded function with $\Vert f \Vert_{U^3} \geq \delta$. Then we have $|\E_{n \in [N]} f(n)\overline{\chi(n)}| \gg_{\delta} 1$, where 
\[ \chi(n) = e(\alpha n^2 + \beta n)\prod_{1 \leq i < i' \leq k} F_{[i',i]}^{m_{[i',i]}}(g(n)\Gamma)\] is the product of some free 2-step nilcharacters with a quadratic phase. Here, $k = O_{\delta}(1)$ and the $m_{[i',i]}$ are integers bounded by $O_{\delta}(1)$. 
\end{theorem}
\proof In \cite[Theorem 10.9]{green-tao-u3inverse} it is shown that a function $f$ with $\Vert f \Vert_{U^3} \geq \delta$ has inner product $\gg_{\delta} 1$ with a function which is the product of $O_{\delta}(1)$ bracket quadratics $e(\alpha_i n\lfloor \beta_i n\rfloor)$, $i = 1,\dots,m$ and a quadratic phase $e(\alpha n^2 + \beta n)$. But such a function already has the form $\chi$ given in the statement of the theorem, simply by taking $k = 2m$ and horizontal frequencies $\xi_{2i-1} = \beta_i$ and $\xi_{2i} = \alpha_i$, $i = 1,2,\dots,m$. \endproof

\emph{Remark.} The proof of \cite[Theorem 10.9]{green-tao-u3inverse} was actually a stepping stone on the way to the proof of the $U^3$ inverse theorem itself, which requires these 2-step nilcharacters to be assembled into a \emph{Lipschitz} Heisenberg nilsequence.\vspace{11pt}

\emph{Remark.} It is possible to proceed directly from the $U^3$ inverse theorem, that is to say from the formulation given in Conjecture \ref{gis-conj}, although -- as the previous remark suggests -- it would be a little perverse to do so. To do this requires one to do a slightly odd kind of Fourier decomposition in the coordinate space $(t_i,t_{[i',i]}) = [0,1]^{k + \binom{k}{2}}$,mapped onto the torus $(\R/\Z)^{k + \binom{k}{2}}$, but there is an issue because a function $F$ which is Lipschitz on $G/\Gamma$ need not even be continuous on this torus. We have a way around this difficulty involving the introduction of a random shift to the fundamental domain $\mathcal{F}$. However we do not believe this argument will be necessary even in the more general paper to come, since our plan is to first prove a variant form of Conjecture \ref{gis-conj}, akin to Theorem \ref{u3-variant}, by induction and only then to deduce Conjecture \ref{gis-conj} itself.\vspace{11pt}

To conclude this section we give some crucial identities involving bracket quadratics. It is the proper understanding and generalisation of these that we referred to above when we talked about the development of a ``calculus'' of bracket polynomials in the forthcoming longer paper.

The key identity we shall rely on is
\begin{equation}\label{key-bracket} 
    X[Y] = XY - \{X\} \{Y\} - [X]Y + [X][Y],
\end{equation} 
valid for all $X,Y \in \R$. This implies that the map $\phi : (X,Y) \mapsto X[Y] \md{1}$ is ``antisymmetric and bilinear modulo lower order terms''. Specifically, $\phi(X_1 + X_2 ,Y) - \phi(X_1,Y) - \phi(X_2,Y) = 0$, whilst $\phi(X,Y_1 + Y_2) - \phi(X,Y_1) - \phi(X,Y_2) = \{X\}\{Y_1\} + \{X\} \{Y_2\} - \{X\}\{Y_1 + Y_2\}$, and $\phi(X,Y) = \phi(Y,X) + XY-\{X\}\{Y\}$.

Let us say a clarify to some extent what we mean by ``lower order''. We shall be applying these identities when $X_i = \alpha_i n$ and $Y_j = \beta_j n$, and we shall also be considering $e(\phi)$ rather than $\phi$ itself. Then these obstructions to antisymmetric bilinearity take the form $e(\{\alpha n\}\{\beta n\})$, an almost 1-step nilsequence (cf. Lemma \ref{1-step-approx} (iv)) and $e(\theta n^2)$, another 1-step object (but of degree $2$).

Let us record these observations in the form of a lemma.

\begin{lemma}[Bracket quadratic identities]\label{brack-identities}
Suppose that $\alpha,\alpha_1,\alpha_2,\beta,\beta_1,\beta_2,\gamma \in \R$. Then 
\begin{enumerate}
\item $e((\alpha_1 + \alpha_2)n \lfloor \beta n\rfloor) = e(\alpha_1 n\lfloor \beta n \rfloor) e(\alpha_2 n \lfloor \beta n \rfloor)$;
\item $e(\alpha n\lfloor (\beta_1 + \beta_2)n\rfloor) = e(\alpha n\lfloor \beta_1 n\rfloor) e(\alpha n\lfloor \beta_2 n\rfloor)$ up to a product of terms of the form $e(\{\theta n\}\{\theta' n\})$; 
\item $e(\alpha n\lfloor \beta n \rfloor) = e(-\beta n\lfloor \alpha n \rfloor)$ up to a product of terms of the form $e(\theta n^2)$ and $e(\{\theta' n\}\{\theta'' n\})$. 
\item $e(\gamma n\lfloor \gamma n\rfloor)$ is a product of terms of the form $e(\theta n^2)$ and $e(\{\theta' n\}\{\theta'' n\})$.
\end{enumerate}
\end{lemma}
\proof The first three of these follow immediately from \eqref{key-bracket} and the subsequent discussion. Part (iv) perhaps requires some comment: to prove it, first choose $\gamma' $ so that $2\gamma' = \gamma \md{1}$. Then take $\alpha = \beta = \gamma'$ in (iii) to obtain the fact that $e(\gamma n\lfloor \gamma' n\rfloor)$ is a product of terms of the required type. Now apply (ii) to conclude the same thing for $e(\gamma n\lfloor \gamma n\rfloor)$. \endproof

\section{Some arguments of Gowers}\label{gowers-sec}

In this section we give the observation of Gowers\cite{gowers-4aps} described in \S \ref{inductive-sec}, whereby one proceeds from the assumption that $|\E_{n} \Delta_h f(n)\overline{ \chi_h(n)} | \geq \delta$ for many $f$ to get \eqref{gowers-conclusion}, a kind of weak linearity statement concerning the map $h \mapsto \chi_h$. Here is a more precise statement. 

\begin{proposition}[Gowers]\label{gowers-prop}Suppose that $f : [N] \rightarrow \C$ is a 1-bounded function, that $H \subseteq [N]$ is a set with cardinality $\eta N$ and that for each $h \in H$ we have a function $\chi_h : [N] \rightarrow \C$ with $|\chi_h(n)| \leq 1$ for all $n$, and that 
\begin{equation}\label{derivative} |\E_{n \in [N]} \Delta_h f(n) \overline{\chi_h(n)}| \geq \delta\end{equation} for all $h \in H$.
Then for at least $\eta^8\delta^4 N^3/2$ of the quadruples $(h_1,h_2,h_3,h_4)$ such that $h_1 + h_2 = h_3 + h_4$ we have
\[ |\E_{n} \chi_{h_1}(n) \chi_{h_2} (n + h_1 - h_4) \overline{\chi_{h_3}(n)} \overline{\chi_{h_4}(n + h_1 - h_4)}| \geq c\eta^4 \delta^2.\]
\end{proposition} 

\emph{Remark.}  In the original paper \cite{gowers-4aps}, attention is restricted to the linear case $\chi_h(n) = e(\xi_h n)$, but the argument extends without difficulty to the general case, as we shall see in the proof.\vspace{11pt}

\proof As in many arguments of analytic number theory and additive combinatorics in which a function that one does not wish to understand is to be eliminated, our main tool is the Cauchy-Schwarz inequality. Two applications of that inequality give that
\begin{equation}\label{to-cbs} |\E_{n,m} a_n b_m \Phi(n,m)|^4 \leq \E_{n,n',m,m'} \Phi(n,m)\overline{\Phi(n',m)\Phi(n,m')}\Phi(n',m')\end{equation} whenever $(a_n)_{n \in X}, (b_m)_{m \in Y}, (\Phi(n,m))_{n \in X, m \in Y}$ are 1-bounded sequences of complex numbers.

Returning to the proposition itself, the assumptions imply that
\[ \E_h |\E_{n} \Delta_h f (n) \overline{\chi_h(n)}|^2 \gg \eta \delta^2,\] where we have taken the expectation over some group $\Z/N'\Z$ with $N' \sim 2N$ (say) and define all functions to be zero outside of $[N]$ and $\chi_h$ to be identically zero if $h \notin H$.
Expanding out and making some obvious substitutions this yields
\[ \E_k \E_{n,m} f(m) \overline{f(m+k)} \overline{f(n)} f(n+k) \Delta_k \chi_{m-n}(n) \gg \eta \delta^2.\] Applying H\"older's inequality this means that
\[ \E_k |\E_{n,m} f(m) \overline{f(m+k)} \overline{f(n)} f(n+k) \Delta_k \chi_{m-n}(n)|^4 \gg \eta^4 \delta^8.\] 
Applying \eqref{to-cbs} for each $k$, we obtain 
\[ \E_k \E_{n,n',m,m'} \Delta_k \chi_{m-n}(n)  \overline{\Delta_k \chi_{m'-n}(n) \Delta_k \chi_{m-n'}(n')} \Delta_k \chi_{m'-n'}(n') \gg\eta^4 \delta^8.\] 

This is more suggestively written as
\[ \E_{h_1 + h_2 = h_3 + h_4} \E_{n,k} \Delta_k \chi_{h_1}(n) \Delta_k \chi_{h_2}(n + h_1 - h_4) \overline{\Delta_k \chi_{h_3}(n)} \overline{\Delta_k \chi_{h_4} (n + h_1 - h_4)} \gg \eta^4 \delta^8,\]
which is the same as
\[ \E_{h_1 + h_2 = h_3 + h_4} |\E_n \chi_{h_1}(n) \chi_{h_2}(n + h_1 - h_4) \overline{\chi_{h_3}(n) \chi_{h_4}(n + h_1 - h_4)}|^2 \gg \eta^4 \delta^8.\]
This immediately implies the stated result by a trivial averaging argument.\endproof

We now give a corollary of this in the specific case that the $\chi_h(n)$ have the form appearing in the statement of Theorem \ref{u3-variant}, that is to say
\begin{equation}\label{chi-form} \chi_h(n) :=  e(\alpha_h n^2 + \beta_h n)\prod_{1 \leq i < i' \leq k} F_{[i',i]}^{m_{[i',i](h)}}(g_h(n)\Gamma).\end{equation}
\begin{corollary}\label{gowers-cor}
Suppose that $f : [N] \rightarrow \C$ is a 1-bounded function and that 
\[ |\E_{n \in [N]} \Delta_h f(n) \overline{\chi_h(n)}| \geq \delta\] for all $h$ in some set $H$, $|H| \geq \delta N$. Suppose now that the functions $\chi_h(n)$ have the specific form \eqref{chi-form}. Then for at least $\delta^C N^3$ additive quadruples $h_1 + h_2 = h_3 + h_4 \in H^4$ there are frequencies $\alpha_{h_1,h_2,h_3,h_4},\beta_{h_1,h_2,h_3,h_4} \in \R/\Z$ such that 
\[ |\E_{n \in [N]} \chi_{h_1}(n) \chi_{h_2}(n)\overline{ \chi_{h_3}(n) \chi_{h_4}(n)} e(\alpha_{h_1,h_2,h_3,h_4} n^2 + \beta_{h_1,h_2,h_3,h_4} n)| \gg_{\delta} 1.\]
\end{corollary}
\proof Apply Proposition \ref{gowers-prop} and then use the bracket identities in Lemma \ref{brack-identities} to expand out terms such as $\chi_{h_4}(n + h_1 - h_4)$. This exhibits 
\[  \chi_{h_1}(n) \chi_{h_2}(n + h_1 - h_4) \overline{\chi_{h_3}(n) \chi_{h_4}(n + h_1 - h_4)}\]
as a product of 
\[ \chi_{h_1}(n) \chi_{h_2}(n)\overline{\chi_{h_3}(n) \chi_{h_4}(n)} \] times various (possibly $h_i$-dependent) terms of the form $e(\alpha n^2)$ or  $e(\{\alpha n\} \{\beta n\})$. By Lemma \ref{1-step-approx} (iv) the latter are almost 1-step nilsequences; the conclusion then follows from Lemma \ref{1-step-to-phase}. \endproof

\emph{Remark.} That this computation worked was no accident. In fact from the general theory in \cite{green-tao-nilratner} one knows that if $\chi(n)$ is a Lipschitz  $s$-step nilsequence with a vertical character then $\chi(n+k)\overline{\chi(n)}$ is an $(s-1)$-step nilsequence. We did not apply this general theory here, since we are being forced to deal with the coordinate functions $F_{[i,i']}$ which are not Lipschitz.

\section{Step 1: Reducing the $h$-dependence}\label{step1-sec}

The aim of this rather long and technical section is to handle \emph{Step 1} of the outline in \S \ref{inductive-sec}. Our first task is to formulate properly exactly what it is we intend to do. Recall that if $\Vert f \Vert_{U^4} \geq \delta$ then, from the fact that $\Vert \Delta_h f \Vert_{U^3} \gg \delta^C$ for $\gg \delta^C N$ values of $h$ and Theorem \ref{u3-variant} we have
\begin{equation}\label{der-nil-cor}
\E_{n \in [N]} \Delta_h f (n)\overline{\chi_h(n)} \gg_{\delta} 1
\end{equation} 
where $\chi_h$ is an object having the form \eqref{chi-form}, that is to say
\begin{equation}\label{chi-form-repeat}
 \chi_h(n) = e(\alpha_h n^2 + \beta_h n)\prod_{1 \leq i < i' \leq k} F_{[i',i]}^{m_{[i',i]}(h)}(g_h(n)\Gamma).\end{equation} 
 Each term involving an $F_{[i,i']}$ is, by the calculations in \S \ref{sec3} and in particular those around \eqref{5.33}, essentially a bracket quadratic $e(\xi n\lfloor \xi' n\rfloor)$ involving $\xi,\xi'$, two of the frequencies  in the ``horizontal'' part of the polynomial sequence $g_h(n)$. 
 
Let us be a little more precise and write $g_h(n) = (\xi_{h,1}n,\dots,\xi_{h,k}n,\dots)$; thus the numbers $\xi_{h,i}$ are the horizontal frequencies just alluded to. Write $\Xi_h := \{\xi_{h,1},\dots,\xi_{h,k}\}$ for this set. When we outlined \emph{Step 1} earlier on, we did little more than suggest that our aim was to show that no bracket quadratic $e(\xi_{h,i} n\lfloor \xi_{h,i'} n\rfloor)$ involving two genuinely $h$-dependent frequencies $\xi_{h,i}$ actually occurs in the formula for $\chi_h(n)$.

To attach meaning to this, we will split $\Xi_h$ as a union $\Xi_* \cup \Xi'_h$ of a ``core'' set $\Xi_* = \{\xi_{h,1},\dots,\xi_{h,k_*}\}$ and a ``petal'' set $\Xi'_h = \{\xi_{h,k_* + 1},\dots,\xi_{h,k}\}$ in such a way that the frequencies $\xi_{h,i}$, $i = 1,\dots, k_*$, do not actually depend on $h$. Our task, then, is to show that \eqref{der-nil-cor} and \eqref{chi-form-repeat} may be achieved in such a way that $m_{[i,i']}(h) = 0$ when $i,i' > k_*$. In other words, no bracket quadratic $e(\xi_{h,i} n\lfloor \xi_{h,i'} n\rfloor)$ with $i, i > k_*$ actually occurs in the expression for $\chi_h(n)$.

We will not prove that any situation such as \eqref{der-nil-cor} and \eqref{chi-form-repeat} has this form automatically. Rather, we will perform an inductive procedure in which the underlying frequency sets are slowly modified so that they take on more and more characteristics of the above ``sunflower'' decomposition into core and petals. At the same time, the set of $h$ for which \eqref{der-nil-cor} holds will be gradually reduced, although it will always have cardinality $\gg_{\delta} N$.

Here is a precise statement.

\begin{proposition}[\emph{Step 1}]\label{step-1-statement}
Suppose that $\Vert f \Vert_{U^4} \geq \delta$. Then for $\gg_{\delta} N$ values of $h$ we have $\E_{n \in [N]} \Delta_h f (n)\overline{\chi_h(n)} \gg_{\delta} 1$, where
where \[ \chi_h(n) = e(\alpha_h n^2 + \beta_h n)\prod_{1 \leq i < i' \leq k} F_{[i',i]}^{m_{[i',i]}(h)}(g_h(n)\Gamma)\]
with $k, |m_{[i,i']}(h)| = O_{\delta}(1)$. Furthermore there is a ``sunflower'' decomposition of the frequency sets $\Xi_h = \{\xi_{h,1},\dots,\xi_{h,k}\}$ of $g_h(n)$ into a ``core'' $\Xi_* = \{\xi_{h,1},\dots,\xi_{h,k_*}\}$ which does not depend on $h$ together with ``petals'' $\Xi'_h = \{ \xi_{h,k_* + 1},\dots,\xi_{h,k}\}$, in such a way that $m_{[i,i']}(h) = 0$ if $i, i' > k_*$. 
\end{proposition}

The last statement -- that is to say the assertion that there are no bracket quadratics with two petal frequencies -- is of course the beef here.

Here is a plan of the rest of this section. Proposition \ref{step-1-statement} is proved by a kind of induction (on the ``complexity'' of the core-petal decomposition). The inductive step is stated as Proposition \ref{prop-h-dep} below, and we give the full derivation of Proposition \ref{step-1-statement} shortly after the proof of that. Proposition \ref{prop-h-dep} is itself deduced from Corollary \ref{cor-h-dep}, which is in turn an easy deduction from Lemma \ref{h-dep-key}. This latter result is the main business of this section, and indeed is probably the hardest part of the entire argument. For that reason we will, between stating it and proving it, give a kind of model variant of the argument to illustrate the underlying algebraic structure. 

Before we can begin we require a definition which will also feature later in the paper.

\begin{definition}[Approximate relations and dissociativity]\label{approx-rel-def}
Suppose that \[\Xi = \{\xi_1,\dots,\xi_k\} \subseteq \R/\Z\] is a finite set of frequencies. We say that this set \emph{satisfies an $M$-linear relation} up to $\eps$ if there are integers $m_1,\dots,m_k$, $|m_i| \leq M$, not all zero, such that $\Vert m_1 \xi_1 + \dots + m_k \xi_k\Vert_{\R/\Z} \leq \eps$.  If a set $\Xi$ satisfies no such linear relation then we say that it is \emph{$(M,\eps)$-dissociated}. We say that a further frequency $\xi$ lies in the $M$-linear span of $\Xi$ up to $\eps$ if there are integers $m_1,\dots,m_k$, $|m_i| \leq M$, such that $\Vert \xi - m_1 \xi_1 - \dots - m_k \xi_k\Vert_{\R/\Z} \leq \eps$. 
\end{definition}

Let us now state the main lemma of this section. We remark that the hypothesis of this lemma comes from applying Proposition \ref{gowers-prop} to the assumption \eqref{der-nil-cor}. However we shall revisit this point later on when we actually perform the inductive application of the lemma.

\begin{lemma}\label{h-dep-key} 
Fix $h_1,h_2,h_3,h_4 \in H$  and suppose that for $j = 1,\dots,4$ we have a decomposition $\Xi_{h_j} = \Xi_* \cup \Xi'_{h_j}$ of the frequency set $\Xi_{h_j}$ into a core $\Xi_* = \{\xi_{h_j,1},\dots,\xi_{h_j,k_*}\}$ not depending on $j$ and a petal set $\Xi'_{h_j} = \{\xi_{h_j,k_* + 1},\dots,\xi_{h_j,k}\}$. Suppose that the functions $\chi_{h}(n)$ have the form \eqref{chi-form} above, where both $k$ and the indices $m_{[i',i]}(h)$ are bounded by $M$, and suppose that we have 
\begin{equation}\label{eq-h-dep} 
|\E_{n \in [N]} \chi_{h_1}(n)\chi_{h_2}(n)\chi_{h_3}(n)\chi_{h_4}(n)e(\alpha_{h_1,h_2,h_3,h_4} n^2 + \beta_{h_1,h_2,h_3,h_4} n) |  \geq 1/M.
\end{equation} 
 Suppose that  $m_{[i',i]}(h_1) \neq 0$ for some pair $i,i' > k_*$. Then either there is an $O_{M}(1)$-linear relation, up to $O_{M}(1/N)$, between the elements in $\Xi_* \cup \Xi'_{h_1} \cup \Xi'_{h_2} \cup \Xi'_{h_3}$, or else there is such a relation between the elements of $\Xi_* \cup \Xi'_{h_1} \cup \Xi'_{h_2} \cup \Xi'_{h_4}$. 
\end{lemma}

\proof The main idea is to apply the distributional results on nilsequences, and in particular the ``Quantitative Ratner'' result, Theorem \ref{quant-ratner}, to the assumption \eqref{eq-h-dep}. There is a very natural way to do this, which is to 
write \eqref{eq-h-dep} as
\begin{equation}\label{eq-h-dep-nil} |\E_{n \in [N]} \tilde F(\tilde g(n)\tilde \Gamma)| \geq 1/M,\end{equation}
where of course
\[ 
\tilde F(\tilde g(n)\tilde \Gamma) := e(\alpha_{h_1,h_2,h_3,h_4} n^2 + \beta_{h_1,h_2,h_3,h_4} n)\prod_{1 \leq i < i' \leq k} \prod_{j=1,2,3,4}F_{[i',i]}^{m_{[i',i]}(h_j)}(g_{h_j}(n)\Gamma) .
\]  
We may interpret the left-hand side as one big polynomial nilsequence on the $2$-step nilmanifold $\tilde G/\tilde \Gamma$, where $\tilde G = G \times G \times G \times G \times\R$ and $\tilde \Gamma = \Gamma \times \Gamma \times \Gamma \times \Gamma \times \Z$, and the polynomial sequence $\tilde g = \tilde g_{h_1,h_2,h_3,h_4}(n)$ is given by 
\[ 
\tilde g(n) = \tilde g_{h_1,h_2,h_3,h_4}(n) = g_{h_1}(n) \times g_{h_2}(n) \times g_{h_3}(n) \times g_{h_4}(n) \times (\alpha_{h_1,h_2,h_3,h_4}n^2 + \beta_{h_1,h_2,h_3,h_4}n).
\] 

The Quantitative Ratner results are a little complicated, and so before continuing with the proof we sketch how it goes in what might be termed the \emph{asymptotic limit case}, in which we work not with any given scale $N$, but rather with the limiting behaviour as $N \to \infty$.  More precisely, instead of \eqref{eq-h-dep-nil} we assume merely that\footnote{The convergence of all limits involving polynomial nilsequences was established in \cite{leibman-poly}, at least in the case when $\tilde F$ is continuous.} 
\[ \lim_{N \rightarrow \infty} \E_{n \in [N]} \tilde F(\tilde g(n)\tilde \Gamma) \neq 0,\] and instead of finding quantitative relations amongst the frequency sets we merely conclude that the frequencies in 
either $\Xi_* \cup \Xi_{h_1}\cup \Xi_{h_2} \cup \Xi_{h_3}$ or $\Xi_* \cup \Xi_{h_1}\cup \Xi_{h_2} \cup \Xi_{h_4}$ are rationally dependent. The main difference between the model case and the actual one is that the corresponding nilmanifold distribution results, due to Leibman \cite{leibman-poly}, are much cleaner in this setting. For simplicity of notation (in this sketch) let us suppose that $\Xi_* = \emptyset$. 

Suppose, then that $\Xi_{h_1}\cup \Xi_{h_2} \cup \Xi_{h_3}$ and $\Xi_{h_1}\cup \Xi_{h_2} \cup \Xi_{h_4}$ are both rationally independent. 
Consider the orbit $(\tilde g(n)\tilde \Gamma)_{n \in \N}$ in $\tilde G$. Roughly speaking\footnote{In actual fact this is only true after subdividing $\N$ into finitely many subprogressions, and furthermore we need to work with a translate $x_0 H\tilde \Gamma/\tilde \Gamma$. Both of these points are merely technical. The finitary analogue of this result is, of course, the Quantitative Ratner Theorem, Theorem \ref{quant-ratner}.}, the results of \cite{leibman-poly} assert that this orbit is equidistributed on a subnilmanifold of the form $H\tilde \Gamma/\tilde \Gamma$, where $H$ is a closed connected rational subgroup of $\tilde G$. 

If $\tilde F$ were continuous then this would imply that 
\begin{equation}\label{dyn-equi}
\lim_{N \to  \infty}\E_{n \in [N]}\tilde F(\tilde g(n)\tilde \Gamma) =\int \tilde F(x) dm_H(x) \neq 0
\end{equation}
where $m_H$ is the Haar measure on $H\tilde \Gamma/\tilde \Gamma$. Unfortunately $\tilde F$ is not quite continuous, a further technicality we will have to handle when discussing the proof of Lemma \ref{h-dep-key} proper. For the purposes of this sketch, however, let us assume that \eqref{dyn-equi} holds.

Let $\pi_1,\pi_2,\pi_3,\pi_4 :\tilde G \to G$ be the projections of $\tilde G$ onto each of the four factors of $G$ comprising $\tilde G$, and by abuse of notation use the same notation for the projection maps from $\tilde G/\tilde \Gamma$ to the factors $G/\Gamma$. Now the projection $(\pi_1 \times \pi_2 \times \pi_3)(\tilde g(n))$ has, as its set of horizontal frequencies, $\Xi_{h_1} \cup \Xi_{h_2} \cup \Xi_{h_3}$, a set which is rationally independent. But these frequencies are precisely those occurring in the projection of $(\pi_1 \times \pi_2 \times \pi_3)(\tilde g(n)\tilde \Gamma)$ onto the horizontal torus (abelianisation) of $G/\Gamma \times G/\Gamma \times G/\Gamma$, and hence the orbit of this (abelian) nilsequence is dense. However Leibman's criterion\footnote{The finitary analogue of this is the Quantitative Leibman dichotomy, Theorem \ref{quant-leib}.} asserts that a polynomial nilsequence is dense if and only if its abelianisation is, and so $((\pi_1 \times \pi_2 \times \pi_3)(\tilde g(n)\tilde \Gamma))_{n \in \N}$ is dense in $G/\Gamma \times G/\Gamma \times G/\Gamma$. 

Since $(\tilde g(n)\tilde \Gamma)_{n \in \N}$ equidistributes in $H\tilde \Gamma/\tilde \Gamma$, we must have
\[ (\pi_1 \times \pi_2 \times \pi_3)(H\tilde \Gamma/\tilde \Gamma) = G/\Gamma \times G/\Gamma \times G/\Gamma.\]
Topological arguments\footnote{In the finitary world these are somewhat painful and involve, for example, some quantitative linear algebra; see Appendix \ref{lift-app}.} using the fact that $H$ is closed and connected let us lift this statement to $G$ to conclude that 
\begin{equation}\label{123} (\pi_1 \times \pi_2 \times \pi_3)(H) = G \times G \times G.\end{equation}
By exactly the same argument we have
\begin{equation} \label{124}(\pi_1 \times \pi_2 \times \pi_4)(H) = G \times G \times G.\end{equation}
We claim that as a consequence of these observations we have \[ [G,G] \times \id \times \id \times \id \times \id \subseteq H.\]
To see this, let $g,g' \in G$ be arbitrary. Then \eqref{123} implies that $H$ contains an element of the form $(g,\id,\id,x,z)$, for some $x \in G$ and some $z \in \R$, whilst \eqref{124} implies that $H$ contains an element of the form $(\id, g',x', \id,z')$ for some $x' \in G$ and some $z' \in \R$. The commutator of these two elements is $([g,g'],\id,\id,\id,\id)$, thereby establishing the claim.

\emph{Remark.} This idea has appeared in related contexts before, for example in the work of Furstenberg and Weiss \cite{furst-weiss}, as well as in less related contexts such as a paper of Hrushovski \cite[Lemma 4.11]{hrushovski}.

As a special case of the above claim, we see that for each pair $i,i'$ with $1 \leq i < i' \leq k$ and for each $t \in \R$ the element $z := (e_{[i,i']}^t,\id,\id,\id,\id)$ lies in $H$. It follows that
\[ \int \tilde F(x) dm_H(x) = \int\tilde F(zx) d m_H(x).\]
However a direct calculation using the definition of $\tilde F$ confirms that 
\[ \tilde F(zx) = e(t m_{[i,i']}(h_1)) \tilde F(x).\]
Since $t$ is arbitrary, the only way to reconcile this with \eqref{dyn-equi} is to conclude that $m_{[i,i']}(h_1) = 0$. Thus in this case (in which there is no core $\Xi_*$) we see that either the functions $\chi_h(n)$ are somewhat trivial in the sense that all of the $m_{[i,i']}(h)$ vanish, or else we were wrong to assume that the frequencies in both $\Xi_{h_1} \cup \Xi_{h_2} \cup \Xi_{h_3}$ and $\Xi_{h_1} \cup \Xi_{h_2} \cup \Xi_{h_4}$ are rationally independent.

This concludes our sketch of the asymptotic limit case, and we now return to our original task of proving Lemma \ref{h-dep-key}. The underlying idea is the same as in the above sketch except that everything must be made quantitative, without any recourse to limits. Furthermore there was one point in the above sketch where we treated a special case (the core set $\Xi_*$ is empty) and others where we waived our hands somewhat (the function $\tilde F$ is not Lipschitz, the orbit only equidistributes on a \emph{coset} of a nilmanifold, and then only after passing to a subprogression). These issues must, of course, be dealt with properly.

Consider the orbit $(\tilde g(n)\tilde \Gamma)_{n \in [N]}$.  
Let $\omega : \R^+ \rightarrow \R^+$ be a growth function to be specified later. By Theorem \ref{quant-ratner} there is some $M_0 = O_{\omega,M}(1)$ (which we may clearly assume to be at least $\max(M,\# \tilde G/\tilde \Gamma)$, since both of these quantities are $O_M(1)$) with the following property. We may partition $[N]$ into subprogressions $P_j$ with lengths at least $N/M_0$, such that corresponding to each progression $P_{j}$ the uniform measure
\[ 
\mu_{j} := \frac{1}{|P_{j}|} \sum_{n \in P_{j}} 1_{\tilde g(n)\tilde\Gamma}
\] 
is $1/\omega(M_0)$-close to the Haar measure $m_{H_j}$ on $b_{j}H_{j}\tilde\Gamma/\tilde\Gamma$, where $H_{j} \leq \tilde G$ is some closed, connected, $M_0$-rational subgroup. 
Namely for any Lipschitz function $F$ on $\tilde G/\tilde \Gamma$ we have
\begin{equation}\label{equiHj}
|\E_{n \in P}  F(\tilde g(n)\tilde \Gamma) - \int F dm_{H_j}| \leq \frac{1}{\omega(M_0)}\|F\|_{\Lip}
\end{equation}

By a trivial averaging argument, condition (\ref{eq-h-dep-nil}) implies that there is some $P = P_j$ such that  
\begin{equation}\label{eq-h-dep-nil*}
|\E_{n \in P} \tilde F(\tilde g(n)\tilde \Gamma)| \geq 1/M;
\end{equation}  
let $H = H_{j}$ be the corresponding group, and $m_H$ the Haar measure on  $bH \tilde\Gamma/\tilde\Gamma$.

Let  $z \in [H,H]$ be an element, all of whose coordinates are bounded by $O_{M_0}(1)$, and let
$F:\tilde G/\tilde \Gamma \to \C$ be a Lipschitz function. Then 
$F_z(x\tilde \Gamma)=F(zx\tilde \Gamma)$  is also Lipschitz and $\|F_z\|_{\Lip}=O_{M_0}(1)\|F\|_{\Lip}$.
Furthermore since $m_H$ is invariant under translation by $z$ (which lies in the centre of $G$) we have
\[
\int F_z dm_{H}=\int F dm_{H},
\]
and thus from (\ref{equiHj}) we get 
\[
|\E_{n \in P}  F_z(\tilde g(n)\tilde \Gamma) - \int F dm_{H}| =O_{M_0}(1/\omega(M_0))\|F\|_{\Lip}.
\]
And by the triangle inequality
\begin{equation}\label{center-inv}
|\E_{n \in P}  F_z(\tilde g(n)\tilde \Gamma) -\E_{n \in P}  F(\tilde g(n)\tilde \Gamma) | =O_{M_0}(1/\omega(M_0))\|F\|_{\Lip},
\end{equation}
thus if $\omega$ is sufficiently rapidly-growing then the error term here is negligible and thus 
\begin{equation}\label{approx} 
\E_{n \in P}  F(z\tilde g(n) \tilde \Gamma) \approx \E_{n \in P}F(\tilde g(n)\tilde \Gamma).
\end{equation}

Let $\epsilon_{M_0}$ be the quantity from the lifting Proposition \ref{lift-lemma}. Namely any element
of $x$ of $H[\tilde G,\tilde G]\tilde \Gamma/ [\tilde G,\tilde G]\tilde \Gamma=\R^m/\Z^m$ whose coordinates are bounded by $\epsilon_{M_0}$  has a lift under the natural projection $\tilde G \to \tilde G/ [\tilde G,\tilde G]\tilde \Gamma$  to an element in $H$ with coordinates $O_{M_0}(1)$, whose first $m$ coordinates are the reduced coordinates of $x$.  

We now deal with the issue of $\tilde F$ not being Lipschitz. Fix $\delta_0=\frac{1}{10}\epsilon_{M_0}$.  We first need to modify the function $\tilde F$.
We will choose two parameters $\delta_1,\delta_2$, such that $\delta_1$ is much smaller that $\delta_0$, and  $\delta_2$ still smaller depending on  $\delta_1$. However, both these quantities will
be $\gg_{M,w} 1$.  Consider the distribution of some fixed coordinate $t_{h_j,[i',i]}$ of $\tilde g(n) \tilde \Gamma$ as $n$ varies over $P$.  We may clearly suppose that there is no $O_{M,\omega}(1)$-linear relation, up to $O_{M,\omega}(1/N)$, amongst the frequencies $\Xi_* \cup \Xi_{h_j}$ since otherwise the conclusion of the lemma is trivially satisfied. If there is no such relation, and if the implicit constants in the $O_{M,\omega}(1)$ notation above are chosen sufficiently large, then by the quantitative Leibman dichotomy, Theorem \ref{quant-leib}, the sequence $(g_{h_j}(n)\Gamma)_{n \in P}$ is $\delta_2$-equidistributed in $G/\Gamma$. Fix a $j \in \{1,2,3,4\}$ and a pair $i,i'$ with $1 \leq i < i' \leq k$. Let $\psi = \psi_{j,i,i'} : G/ \Gamma \rightarrow [0,1]$ be supported where $t_{h_j,[i',i]} \leq 2\delta_1$ or $t_{h_j,[i',i]} \geq 1- 2\delta_1$ and be equal to $1$ whenever $t_{h_j,[i',i]} \leq \delta_1$ or $t_{h_j,[i',i]} \geq 1- \delta_1$ and have $\Vert \psi \Vert_{\Lip} = O_{M,\delta_1}(1)$. Let $\tilde \psi = \tilde \psi_{j,i,i'} : \tilde G/\tilde \Gamma \rightarrow [0,1]$ be the pullback of $\psi_{j,i,i'}$ under the natural projection from $\tilde G$ to the $j$th copy of $G$.

Our preceding observation about the distribution of $(g_{h_j}(n)\Gamma)_{n \in P}$ implies that 
\begin{align*} 
|\E_{n \in P} \tilde \psi(\tilde g(n)\tilde\Gamma)|  & =  |\E_{n \in P} \psi (g_{h_j}(n)\Gamma)|   \leq \int_{G/\Gamma} \psi dm_{G/\Gamma} + \delta_2 \Vert \psi \Vert_{\Lip} \\ & = o_{M;\delta_1\rightarrow 0}(1) + O_{M,\delta_1}(\delta_2),
\end{align*}
where $o_{M; \delta_1 \rightarrow 0}(1)$ denotes a quantity that is bounded in magnitude by $c_M(\delta_1)$ for some $c_M(\delta_1)$ that goes to zero as $\delta_1 \to 0$ for any fixed $M$.

Let $z$ be an element in $[H,H]$ with $O_{M_0}(1)$-bounded coordinates, and suppose  $z = (z_{h_1},z_{h_2},z_{h_3},z_{h_4},w)$ under the decomposition of $\tilde G$ as $G \times G \times G \times G \times \R$. Then $\psi_{z_{h_j}} (g_{h_j}(n)\Gamma)=
\psi(z_{h_j} g_{h_j}(n)\Gamma)$ is Lipschitz with $\|\psi_{z_{h_j}}\|_{\Lip}=O_{M_0}(1)$,
and by the invariance of $m_{G/\Gamma}$ under multiplication by $z_{h_j}$ we get
\[ 
|\E_{n \in P} \tilde \psi(z \tilde g(n)\tilde\Gamma)| = o_{M_0;\delta_1 \rightarrow 0}(1) + O_{M_0,\delta_1}(\delta_2).
\]
By adding, the same type of bounds hold for the function $\Psi := \sum_{j=1}^4 \sum_{1 \leq i < i' \leq k} \tilde \psi_{j,i,i'}$, that is to say
\[  
 |\E_{n \in P} \Psi(\tilde g(n)\tilde\Gamma)| , |\E_{n \in P}  \Psi(z \tilde g(n)\tilde\Gamma)| = o_{M_0;\delta_1\rightarrow 0}(1) + O_{M_0,\delta_1}(\delta_2).
 \]
Let us, at this point, fix $\delta_1 \gg_{M_0} 1$ in such a way that the $o_{M_0;\delta_1 \rightarrow 0}(1)$ term here is bounded by $\delta_0^{10}$ (say), and let us then choose $\delta_2 \gg_{M_0} 1$ in such a way that the $O_{M_0,\delta_1}(\delta_2)$ term is also bounded by $\delta_0^{10}$.
Then the last displayed equation becomes
\begin{equation}\label{lip-fix}   |\E_{n \in P} \Psi(\tilde g(n)\tilde\Gamma)| , |\E_{n \in P}  \Psi(z \tilde g(n)\tilde\Gamma)| = O(\delta_0^{10}).\end{equation}
Note that by construction $\Psi$ is equal to $1$ in a $\delta_1$-neighbourhood of all of the discontinuities of our function $F$. As a result of this it is clear that we may find a function $\tilde F_0$ with the property that \[ \Vert \tilde F_0\Vert_{\Lip} = O_{M_0}(1)\] whilst \[ |\tilde F - \tilde F_0| \leq \Psi\] pointwise. By \eqref{center-inv} we have
\[
 |\E_{n \in P} \tilde F_0(z\tilde g(n)\tilde \Gamma) - \E_{n \in P} \tilde F_0(\tilde g(n)\Gamma)|  = O_{M_0}(1/\omega(M_0))
 \] 
and by \eqref{lip-fix} we have
\begin{align*}
| \E_{n \in P} (\tilde F - \tilde F_0)(z\tilde g(n)\tilde \Gamma) |, | \E_{n \in P} (\tilde F - \tilde F_0)(\tilde g(n)\tilde \Gamma) | = O(\delta_0^{10}).
\end{align*}
Adding, we obtain
\[
|\E_{n \in P} \tilde F(z \tilde g(n)\tilde \Gamma) - \E_{n \in P} \tilde F( \tilde g(n)\tilde \Gamma)|  = O_{M_0}(1/\omega(M_0)) +O(\delta_0^{10}) .
\]
Recall that $\delta_0 = \frac{1}{10}\eps_{M_0}$ depends only on $M_0$. By choosing $\omega : \R^+ \rightarrow \R^+$ to be sufficiently rapidly-growing, the whole of the right-hand side can therefore be made $O(\delta_0^{10})$, that is to say
\begin{equation}\label{approx*}
\E_{n \in P} \tilde F(z \tilde g(n)\tilde \Gamma) = \E_{n \in P} \tilde F( \tilde g(n)\tilde \Gamma)+ O(\delta_0^{10}).
\end{equation}

Now that $\omega$ has been fixed, we have $M_0 = O_M(1)$ and $\delta_0 \gg_M 1$. As stated before, our aim now is to assume that $\Xi_* \cup \Xi_{h_1} \cup \Xi_{h_2} \cup \Xi_{h_3}$ and $\Xi_* \cup \Xi_{h_1} \cup \Xi_{h_2} \cup \Xi_{h_4}$ are highly dissociated and use this to produce an element $z \in [H,H]$ which, in conjunction with \eqref{eq-h-dep-nil}, contradicts \eqref{approx*}. We shall require a further parameter $\delta_3 \gg_M 1$, much smaller than $\delta_0$. We will specify it later on.
  
Let $\pi_1,\pi_2,\pi_3,\pi_4 :\tilde G \to G$ be the projections from $\tilde G$ onto the four copies of $G$ (recall, of course, that $\tilde G = G \times G \times G \times G \times \R$). Once again we abuse notation and use the same notation for the corresponding projections from $\tilde G/\tilde \Gamma$ to $G/\Gamma$. 
Suppose that $\Xi_* \cup \Xi_{h_1}  \cup \Xi_{h_2} \cup \Xi_{h_3}$ is $O_M(1)$-dissociated up to $O_M(1/N)$. Let us examine the abelian part of $(\pi_1 \times \pi_2 \times \pi_3)(\tilde{g}(n)\Gamma)_{n \in P}$, that is to say the image of $(\tilde g(n)\tilde\Gamma)_{n \in P}$ under the projection \[ \pi^{\ab}_{123} : \tilde G/\tilde\Gamma \rightarrow G/[G,G]\Gamma \times G/[G,G]\Gamma \times G/[G,G]\Gamma \cong (\R/\Z)^k \times (\R/\Z)^k \times (\R/\Z)^k.\]  This image takes the form
\[ 
\big( (\xi_{h_1,i}n)_{i=1}^{k_*},(\xi_{h_1,i}n)_{i=k_*+1}^k, (\xi_{h_2,i}n)_{i=1}^{k_*},(\xi_{h_2,i}n)_{i=k_*+1}^k, (\xi_{h_3,i}n)_{i=1}^{k_*},(\xi_{h_3,i}n)_{i=k_*+1}^k \big) \md{1}
\]
Recalling that $\Xi_* = \{\xi_{h,1},\dots,\xi_{h,k_*}\}$ and that $\Xi_h = \{ \xi_{h,k_* + 1},\dots,\xi_{h,k}\}$, it follows from the asserted dissociativity (assuming the implicit $O_M(1)$ terms are large enough) and  Kronecker's theorem in quantitative form (cf. Lemma \ref{kronecker}), that this image is $\delta_3$-equidistributed in the subtorus
\[ 
\{ (t,u_1,t,u_2,t,u_3) : t \in (\R/\Z)^{k_*}, u_1,u_2,u_3 \in (\R/\Z)^{k - k_*} \} \subseteq (\R/\Z)^{k} \times (\R/\Z)^{k} \times (\R/\Z)^k.
\]

In particular there is an element in $(\pi^{\ab}_{123}(\tilde g(n)\tilde \Gamma))_{n \in P}$ within $O(\delta_3)$ of 
\[
 (0, (0,\dots,\delta_0,\dots,0), 0,0,0,0) + \pi_{123}^{\ab}(b),
 \]
where the $\delta_0$ lies in the $i$th position (note that $i > k_*$ by assumption).  Now since the uniform probability measure on $(\tilde g(n)\tilde \Gamma)_{n \in P}$ is $O_{M_0}(1/\omega(M_0)) = \delta_0^{10}$-close to the Haar measure on $b H\tilde \Gamma/\tilde \Gamma$, the projection $(\pi^{\ab}_{123}(\tilde g(n)\tilde \Gamma))_{n \in P}$ is $\delta_0^{10}$-equidistributed in  $\pi^{\ab}_{123}(b H\tilde \Gamma/\tilde \Gamma)$. This means that there is an element $x$ of $\pi^{\ab}_{123}(H\tilde \Gamma/\tilde \Gamma)$ within $O(\delta_3)$ of 
\[ 
(0, (0,\dots,\delta_0,\dots,0), 0,0,0,0) .
\]

Recall that we chose $\delta_0$ so that the lifting property \ref{lift-lemma} holds.
Recalling the relationship between distance in coordinates and distance in $\tilde G$ (cf. \cite[Lemma A.4]{green-tao-nilratner}) we can thus find an element in $(\pi_1 \times \pi_2 \times \pi_3)(H)$ at distance $O_M(\delta_3)$ from 
$e_i^{\delta_0}z_1 \times z_2 \times z_3$ where $z_1,z_2,z_3 \in [G,G]$ are arbitrary (with coordinates  bounded by $O_{M_0}(1)$). It follows that we can find an $g \in H$ with 
\[
d_{\tilde G}(g, e_{i}^{\delta_0}z_1 \times  z_2 \times  z_3 \times  w_4 \times  u ) = O_M(\delta_3)
\]
where  $w_4 \in G$ and $u \in \R$ are arbitrary.

Similarly, if  $\Xi_* \cup \Xi_{h_1}  \cup \Xi_{h_2} \cup \Xi_{h_4}$ is $O_M(1)$-dissociated up to $O_M(1/N)$ then we may locate inside $H$ an element $g'$ with
\[
d_{\tilde G}(g', e_{i'}^{\delta_0}z_1' \times  z'_2 \times w'_3 \times  z'_4 \times u' ) = O_M(\delta_3),
\]
 where  $z_1',z'_2,z'_4 \in [G,G]$, and $w'_3 \in G$,  $u' \in \R$ are arbitrary. 
 
We then take for our element $z \in [H,H]$ the commutator $[g,g']$. Noting that\footnote{Here we have used the fact, specific to the $2$-step case, that $[x^t,y^{t'}] = [x,y]^{tt'}$. One way to check this would be to verify it for $t,t' \in \Z$ and use the fact that both sides are polynomials in a suitable coordinate system. In the higher step case, the more general \emph{Baker-Campbell-Hausdorff formula} could be used instead.}
\[ 
[e_{i'}^{\delta_0}z_1 \times z_2 \times  z_3 \times w_4 \times u,e_{i}^{\delta_0}z'_1 \times  z'_2 \times w'_3 \times  z'_4 \times  u' ] = [e_{i'}^{\delta_0},e_{i}^{\delta_0}] = [e_{i'},e_i]^{\delta_0^2} \times \id \times \id \times \id \times \id
\] 
and that the maps $g \mapsto [g,g_0]$ are uniformly Lipschitz for $g_0$ in any bounded set, we have
\begin{equation}\label{center-elem} 
d_{\tilde G}(z, [e_{i'},e_{i}]^{\delta_0^2} \times \id \times \id \times \id \times \id) = O_{M_0}(\delta_3).\end{equation}
Now, as we have remarked, the coordinate functions $F_{[l',l]} : G/\Gamma \rightarrow \C$ are not Lipschitz. However, they \emph{are} $O_M(1)$-Lipschitz when restricted to $[G,G]\Gamma/\Gamma$, as an easy computation confirms. It follows from this observation, \eqref{center-elem} and the definition of the functions $F_{[l',l]}$ that for $j = 1,2,3,4$ and for any $x \in G/\Gamma$ we have
\[ F_{[l',l]}(zx)^{m_{[l',l]}(h_j)} = F_{[l',l]}(x)^{m_{[l',l]}(h_i) }+ O_{M_0}(\delta_3)\] unless $l = i$, $l' = i'$ and $j = 1$ in which case
\[ F_{[i',i]}(zx)^{m_{[i',i]}(h_1)} = e(\delta_0^2 m_{[i',i]}(h_1)) F_{[i',i]}(x)^{m_{[i',i]}(h_1)} + O_{M_0}(\delta_3).\]
Taking products over all choices of $i,i'$, it follows that 
\[ \tilde F(zx) = e(\delta_0^2 m_{[i',i]}(h_1)) \tilde F(x) + O_{M_0}(\delta_3),\] from which it of course follows that 
\[ \E_{n \in P} \tilde F(z \tilde g(n)\tilde\Gamma) = e(\delta_0^2m_{[i',i]}(h_1)) \E_{n \in P} \tilde F(\tilde g(n)\Gamma) + O_{M_0}(\delta_3).\]
Choosing $\delta_3$ so small that the error term here is $O(\delta_0^{10})$, we obtain upon comparison with \eqref{approx*} that 
\[ |1 - e(\delta_0^{2}m_{[i',i]}(h_1))| \E_{n \in P} \tilde F(\tilde g(n)\Gamma)| = O(\delta_0^{10}).\]
Recalling that $m_{[i',i]}(h_1)$ is an integer bounded in magnitude by $M$, that \[ |\E_{n \in P} \tilde F(\tilde g(n)\tilde \Gamma)| \geq 1/M,\] and that $\delta_0$ may certainly be assumed to be much smaller than $1/M$, we are forced to conclude (at last!) that $m_{[i',i]}(h_1)= 0$.\endproof

We may put Lemma \ref{h-dep-key} together with Corollary \ref{gowers-cor} in a straightforward manner.

\begin{corollary}\label{cor-h-dep}
Suppose that $|\E_{n \in [N]} \Delta_h f(n) \overline{\chi_h(n)}| \geq 1/M$ for all $h \in H$, where $H \subseteq [N]$, $|H|\geq N/M$ and $\chi_h(n)$ has the form \eqref{chi-form} with complexity at most $M$, and with decompositions of the frequency sets $\Xi_h = \{\xi_{h,1},\dots,\xi_{h,k}\}$ into cores $\Xi_* = \{\xi_{h,1},\dots,\xi_{h,k_*}\}$ which do not depend on $h$ and petal sets $\Xi'_h = \{\xi_{h,k_* + 1},\dots,\xi_{h,k}\}$. Then one of the following two alternatives holds true:
\begin{enumerate}
\item There is a set $H' \subseteq H$, $|H'| \gg_M |H|$, such that $m_{[i,i']}(h) = 0$ whenever $i,i' > k_*$ and $h \in H'$;
\item For $\gg_M N^3$ triples $h,h',h'' \in H^3$ the set $\Xi_* \cup \Xi'_{h} \cup \Xi'_{h'} \cup \Xi'_{h''}$ fails to be $O_M(1)$-dissociated up to $O_M(1/N)$.
\end{enumerate}
\end{corollary}
\proof By Corollary \ref{gowers-cor} there are $\gg_M N^3$ additive quadruples $h_1 + h_2 = h_3 + h_4$ such that there are $\alpha_{h_1,h_2,h_3,h_4},\beta_{h_1,h_2,h_3,h_4}$ for which
\[ |\E_{n \in [N]} \chi_{h_1}(n)\chi_{h_2}(n)\chi_{h_3}(n)\chi_{h_4}(n) e(\alpha_{h_1,h_2,h_3,h_4}n^2 + \beta_{h_1,h_2,h_3,h_4} n)| \gg_M 1.\]
By pigeonhole there must either be $\gg_M N^3$ of these quadruples such that $m_{[i',i]}(h_1) = 0$ for all $i,i' > k_*$, in which case we are clearly in alternative (i), or else there must be some choice of $i,i' > k_*$ such that there are $\gg_M N^3$ quadruples with $m_{[i',i]}(h_1) \neq 0$. By Lemma \ref{h-dep-key} it follows that for each of these quadruples at least one of the sets $\Xi_* \cup \Xi'_{h_1} \cup \Xi'_{h_2} \cup \Xi'_{h_3}$ or $\Xi_* \cup \Xi'_{h_1} \cup \Xi'_{h_2} \cup \Xi'_{h_4}$ fails to be $O_M(1)$-dissociated up to $O_M(1/N)$. It follows immediately that we are in case (ii).\endproof

Now if alternative (i) holds in this last corollary then \emph{Step 1} is complete (that is, Proposition \ref{step-1-statement} is proven). If alternative (ii) holds, then it is possible to replace the core-petal decomposition $\Xi_h = \Xi_* \cup \Xi'_h$ by one in which some of the petal behaviour is absorbed into the core. The precise statement of this, which follows now, is slightly long:

\begin{proposition}\label{prop-h-dep} Let $H \subseteq [N]$ be a set with $|H| \geq N/M$. 
Suppose that \[ |\E_{n \in [N]} \Delta_h f(n) \overline{\chi_h(n)}| \geq 1/M\] for all $h \in H$, where $|H| \geq N/M$ and the nilcharacter $\chi_h(n)$ has the form \eqref{chi-form} with complexity at most $M$ and there is a decomposition of the underlying frequency set $\Xi_h = \{\xi_{h,1},\dots,\xi_{h,k}\}$ into \begin{itemize}\item a core component $\Xi_* = \{\xi_{h,1},\dots,\xi_{h,k_*}\}$ which does not depend on $h$ and \item a petal component $\Xi'_h = \{\xi_{h,k_* +1},\dots,\xi_{h,k}\}$.\end{itemize} Then either \begin{itemize}
\item there is a set $H' \subseteq H$, $|H'| \gg_M |H|$, such that $m_{[i',i]}(h) = 0$ for all $i, i' > k_*$ and for all $h \in H'$, or 
\item there is a set $\tilde H \subseteq H$, $|\tilde H| \gg_M |H|$ and nilcharacters $\tilde\chi_h(n)$ of complexity $O_M(1)$, $h\in H''$, such that \begin{equation}\label{new-chi-cor}\E_{n \in [N]} \Delta_h f(n)\overline{\tilde\chi_h(n)} \gg_M 1\end{equation} for all $h \in \tilde H$.  Here the nilcharacters $\tilde\chi_h(n)$ have the form 
\[ 
\tilde\chi_h(n) = e(\alpha_h n^2 + \beta_h n) \prod_{1 \leq i < i' \leq \tilde k} F_{[i',i]}^{\tilde m_{[i',i]}(h)}(\tilde g_h(n)\tilde\Gamma) ,
\] 
where $\tilde g_h(n) = (\tilde\xi_{h,1},\dots,\tilde\xi_{h,\tilde k},0,\dots,0)$. Furthermore writing \[ \tilde \Xi_h := \{ \tilde \xi_{h,1},\dots,\tilde \xi_{h,\tilde k}\}\] we have a decomposition $\tilde \Xi_h = \tilde \Xi_* \cup \tilde \Xi'_h$, where either
\begin{enumerate}
\item \textup{(core decreases)} $|\tilde \Xi_*| < |\Xi_*|$ and $|\tilde \Xi'_h| = |\Xi'_h|$ or
\item \textup{(petals decrease)} $|\tilde \Xi_*| \leq |\Xi_*| + 1$ and $|\tilde \Xi'_h| < |\Xi'_h|$.
\end{enumerate}
\end{itemize}
\end{proposition}

\emph{Proof of Proposition \ref{step-1-statement}, a.k.a. Step 1.} Before embarking on the proof of this last proposition, we remark how a simple iteration of it leads to Proposition \ref{step-1-statement}. One starts with the trivial decomposition $\Xi_h = \Xi_* \cup \Xi'_h$ where $\Xi_* = \emptyset$ and $\Xi'_h = \Xi_h$, and with the initial value of $M$ being $O_{\delta}(1)$. It is not hard to see that there cannot be more than $O_M(1)$ iterations of alternatives (i) (core decreases) or (ii) (petals decrease).\endproof

\emph{Proof of Proposition \ref{prop-h-dep}.} By Corollary \ref{cor-h-dep} we may assume that there are $\gg_M N^3$ triples $h,h',h'' \in H$ such that $\Xi_* \cup \Xi_{h} \cup \Xi_{h'} \cup \Xi_{h''}$ fails to be $O_M(1)$-dissociated up to $O_M(1/N)$. To each such triple is associated a $k_* + 3(k - k_*)$ tuple 
\[ q_{*,1},\dots,q_{*,k_*}, q_{h,k_* + 1},\dots,q_{h,k},q_{h',k_* + 1},\dots,q_{h',k},q_{h'',k_* + 1},\dots,q_{h'',k}\] of integers, all at most $O_M(1)$ in magnitude, such that 
\begin{align*} \Vert q_{*,1} \xi_{*,1} + \dots & + q_{*,k_*} \xi_{*,k_*} + q_{h,k_* + 1} \xi_{h,k_* + 1} + \dots + q_{h,k} \xi_{h,k} + q_{h',k_*+1} \xi_{h,k_* + 1} + \dots \\ & + q_{h',k}\xi_{h',k}  + q_{h'',k_* + 1} \xi_{h'',k_* + 1} + \dots + q_{h'',k} \xi_{h'',k} \Vert_{\R/\Z} = O_M(1/N).\end{align*}
By pigeonholing we may pass to a further subcollection of triples $h,h',h''$ for which these integers $q_{h,j},q_{h',j'},q_{h'',j''}$ have no $h,h',h''$-dependence. If at least one of these latter quantities (with $j > k_*$) is nonzero then by relabeling we may assume it is $q_{h,k}$. 
All this having been done, let us fix $h'$ and $h''$ appearing in $\gg_M N$ of these triples. We then have integers $q_{1},\dots q_k = O_M(1)$, not all zero, and some frequency $\xi_0$ such that 
\[ \Vert \xi_0 + q_1 \xi_{*,1} + \dots + q_{k_*} \xi_{*,k_*} + q_{k_* + 1} \xi_{h,k_* + 1} + \dots + q_k \xi_{h,k}\Vert_{\R/\Z} = O_M(1/N)\] for at least $\gg_M N$ values of $h$. Furthermore (case 1) we have $\xi_0 = 0$ if $q_{k_* + 1} = \dots = q_k = 0$; otherwise (case 2) we have $q_k \neq 0$.

Suppose we are in case 1 and that, without loss of generality, we have $q_{k_*} \neq 0$. Then $\xi_{*,k_*}$ is in the $O_M(1)$-linear span, up to $O_M(1/N)$, of the set $\tilde \Xi_* := \{\frac{1}{Q}\xi_{*,1},\dots,\frac{1}{Q}\xi_{*,k_* - 1}\}$, where $Q$ is the lowest common multiple of the integers up to $O_M(1)$. Taking $\tilde \Xi'_h = \Xi'_h$, we see that (i) is satisfied and also that $\Xi_h$ is in the $O_M(1)$-linear span, up to $O_M(1/N)$, of $\tilde \Xi_* \cup \tilde \Xi'_h$. Suppose now that we are in case 2; then take $\tilde \Xi_* = \frac{1}{Q}\Xi_* \cup \{\frac{1}{Q} \xi_0\}$ and $\tilde\Xi'_h = \frac{1}{Q}\Xi'_h \setminus \{\frac{1}{Q} \xi_{h,k}\}$. Now condition (ii) is satisfied, and once again $\Xi_h$ is in the $O_M(1)$-linear span, up to $O_M(1/N)$, of $\tilde \Xi_* \cup \tilde \Xi'_h$.

The treatment of the two cases is, henceforth, the same and at this point we revert to the bracket quadratic expressions 
\[ F_{[i',i]}^{m_{[i',i]}(h)}(g_h(n)\Gamma) = e(m_{[i',i]}(h)\xi_{h,i'} n\lfloor \xi_{h,i} n\rfloor).\]
For each $\xi_{h,i}$ we substitute in the expression for this frequency as an $O_M(1)$-linear combination of the frequencies in $\tilde \Xi'_* \cup \tilde \Xi'_h$, plus an error which is $O_M(1/N)$. To simplify this we use the bracket identities of Lemma \ref{brack-identities} repeatedly to express the whole product $\chi_h(n)$ as a product of terms $e(\tilde m_{[i',i]}(h)\tilde\xi_{h,i'} n\lfloor \tilde \xi_{h,i} n\rfloor)$ with $i < i'$, where the exponents $\tilde m_{[i',i]}(h)$ are still $O_M(1)$, together with various terms of the form $e(\theta n^2)$, $e(\{\alpha n\} \{\beta n\})$ and $e(\alpha n\lfloor \beta n \rfloor)$ with $\beta = O_M(1/N)$.

Now we may use Lemma \ref{1-step-approx} (ii), (iii) and (iv) repeatedly, bearing in mind the assumption $|\E_{n \in [N]} \Delta_h f(n)\overline{\chi_h(n)}| \geq 1/M$, to remove all terms of these last two types and replace them by a single linear term $e(\theta' n)$. Doing this and then taking the new bracket quadratics $e(\tilde m_{[i',i]}(h)\tilde\xi_{h,i'} n\lfloor \tilde \xi_{h,i} n\rfloor)$ and writing them as nilcoordinate functions $F_{[i',i]}^{\tilde m_{[i',i]}(h)}(\tilde g_h(n)\tilde \Gamma)$, we obtain precisely the desired conclusion \eqref{new-chi-cor}.\endproof

\section{Step 2: Approximate linearity}

In this section we address \emph{Step 2} of the outline in \S \ref{inductive-sec}. In the last section we decomposed the underlying frequency sets $\Xi_h = \{\xi_{h,1},\dots,\xi_{h,k}\}$ into a core set $\Xi_*$ and a petal set $\Xi'_h$, in such a way that no nilcharacter $F_{[i,i']}(g_h(n))$ corresponding to two petal frequencies $\xi_{h,i}, \xi_{h,i'}$ appears in the expression for $\chi_h(n)$.  Our task now is to proceed from here to show that, at least for many $h$, the petal set $\Xi'_h$ has a weak linear structure. There follows a precise statement of what we shall prove. By a \emph{bracket-linear form} of complexity $M$ we mean a function $\psi : \Z \rightarrow \R/\Z$ of the form
\[ \psi(h) = \beta_0 + \alpha_1\{\beta_1 h\} + \dots + \alpha_m\{\beta_m h \} + \theta h,\] where the $\alpha_j,\beta_j,\beta$ lie in $\R$ and $m \leq M$.

\begin{proposition}\label{approx-lin}
Suppose that $f : [N] \rightarrow \C$ is a $1$-bounded function with $\Vert f \Vert_{U^4} \geq \delta$. Then there is a set $H \subseteq [N]$, $|H| \gg_{\delta} N$, such that for all $h \in H$ we have
\[ |\E_{n \in [N]} \Delta_h f(n) \overline{\chi_h(n)}| \gg_{\delta} 1.\]
Here we have
\begin{equation}\label{chi-form-2} \chi_h(n) =  e(\alpha_h n^2 + \beta_h n)\prod_{1 \leq i < i' \leq k} F_{[i',i]}^{m_{[i',i]}(h)}( g_h(n))  \end{equation} with $k, |m_{[i,i']}(h)| = O_{\delta}(1)$,
where $g_h(n) = (\xi_{h,1} n,\dots,\xi_{h,k}n,0,\dots,0)$, $m_{[i',i]}(h) = 0$ if $i,i' > k_*$, the frequency set $\Xi_h$ decomposes as $\Xi_* \cup \Xi'_h$ with $\Xi_* = \{\xi_{h,1},\dots,\xi_{h,k_*}\}$ independent of $h$, and every frequency $\xi_{h,i}$, $i > k_*$, in the petal set $\Xi'_h$ is a bracket linear form in $h$ of complexity $O_{\delta}(1)$.
\end{proposition}

We shall establish this proposition inductively in a manner not too dissimilar to that in the last section. The inductive step which drives Proposition \ref{approx-lin} is the following; it might be compared to Proposition \ref{prop-h-dep} in the last section.

\begin{proposition}\label{approx-lin*}
Suppose that $H \subseteq [N]$ is a set with $|H| \geq N/M$. Suppose that for all $h \in H$ we have
\[ |\E_{n \in [N]} \Delta_h f(n) \overline{\chi_h(n)}| \geq 1/M,\] where $\chi_h(n)$ has the form \eqref{chi-form-2} and the frequency set $\Xi_h$ is decomposed as $\Xi_* \cup \Xi_h^{\struct} \cup \Xi_h^{\unstruct}$, where the frequencies in $\Xi_*$ do not depend on $h$ and those in $\Xi_h^{\struct}$ are bracket-linear in $h$ with complexity at most $M$.  Then there is a set $\tilde H \subseteq H$, $|\tilde H| \gg_M 1$, such that 
\[ |\E_{n \in [N]}\Delta_h f(n) \overline{\tilde \chi_h(n)}| \gg_M 1,\] 
where $\tilde \chi_h(n)$ has the form 
\[ \tilde\chi_h(n) =  e(\tilde\alpha_h n^2 + \tilde\beta_h n)\prod_{1 \leq i < i' \leq \tilde k} F_{[i',i]}^{\tilde m_{[i',i]}(h)}(\tilde g_h(n)\tilde\Gamma) ,\] 
a nilcharacter with complexity $O_M(1)$ in which the frequency set $\tilde \Xi_h$ decomposes as $\tilde \Xi_* \cup \tilde \Xi_h^{\struct} \cup \tilde \Xi_h^{\unstruct}$ where either
\begin{enumerate}
\item \textup{(core decreases)} $|\tilde \Xi_*| < |\Xi_*|$, $|\tilde \Xi_h^{\struct}| \leq |\Xi_h^{\struct}|$, $|\tilde \Xi_h^{\unstruct}| \leq |\tilde \Xi_h^{\unstruct}|$;
\item \textup{(unstructured part decreases)} $|\tilde \Xi_*| = O_{M, |\Xi_*|}(1)$, $|\tilde \Xi_h^{\struct}| = |\Xi_h^{\struct}| + 1$, \\$|\tilde \Xi_h^{\unstruct}| = | \Xi_h^{\unstruct}| - 1$.
\end{enumerate}
\end{proposition}

\noindent\emph{Proof of Proposition \ref{approx-lin} given Proposition \ref{approx-lin*}.} To prove Proposition \ref{approx-lin} one first, of course, applies \emph{Step 1}. With that in hand one may pick $M = O_{\delta}(1)$ and initialise the inductive use of Proposition \ref{approx-lin*} by taking $\Xi_h^{\unstruct}$ to equal to the entire petal frequency $\Xi'_h$ and $\Xi_h^{\struct} = \emptyset$. It is not hard to see that this repeated application of Proposition \ref{approx-lin*} terminates in time $O_M(1)$, at which point the unstructured component $\Xi_h^{\unstruct}$ is empty.\endproof

It remains, of course, to prove Proposition \ref{approx-lin*}, and this will be the main business of this section. Once again the key tool is Proposition \ref{gowers-prop}, of which we require the following variant.

\begin{lemma}\label{gowers-unstruct}
Suppose that $H \subseteq [N]$ is a set with $|H|\geq N/M$ and that \[ |\E_{n \in [N]} \Delta_h f(n) \overline{\chi_h(n)}| \geq 1/M\] for all $h \in H$, where $\chi_h(n)$ has the form \eqref{chi-form-2} with $m_{[i,i']}(h) = 0$ if $i,i' > k_*$ and the underlying frequency set $\Xi_h$ has been decomposed as $\Xi_* \cup \Xi_h^{\struct} \cup \Xi_h^{\unstruct}$, where the core $\Xi_*$ does not depend on $h$ and $\Xi_h^{\struct}$ consists of bracket linear forms of complexity at most $M$. Write $\chi_h(n) = \chi_h^{\struct}(n)\chi_h^{\unstruct}(n)$, where the two parts here correspond to the structured and unstructured frequencies in $\Xi_h$.  Then there is a set $\tilde H \subseteq H$, $|\tilde H| \gg_M |H|$, and frequencies $\alpha_{h_1,h_2,h_3,h_4},\beta_{h_1,h_2,h_3,h_4} \in \R/\Z$ such that 
\[ \E_{n \in [N]}\chi^{\unstruct}_{h_1}(n)\chi^{\unstruct}_{h_2}(n)\chi^{\unstruct}_{h_3}(n)\chi^{\unstruct}_{h_4}(n)e(\alpha_{h_1,h_2,h_3,h_4} n^2 + \beta_{h_1,h_2,h_3,h_4} n) \gg_M 1\]
for $\gg_M N^3$ additive quadruples $h_1 + h_2 = h_3 + h_4 \in \tilde H$.
\end{lemma}
\proof The idea is to apply Proposition \ref{gowers-prop} and then simply observe that the contribution from the structured parts $\chi^{\struct}_h(n)$ can be made to cancel out.  Bracket linear forms are not quite genuinely linear, but if $\psi(h) = \alpha_1\{\beta_1 h \} + \dots + \alpha_m\{\beta_mh\} + \theta h$ then we have $\psi(h_1) + \psi(h_2) = \psi(h_3) + \psi(h_4)$ whenever the tuple $(\beta_1 h,\dots,\beta_m h) \md{1}$ lies in some cube $\prod_{j=1}^m [ i_j/10, \gamma_j + (i_j+1)/10]$ (say), where the $i_j$ are integers between $0$ and $9$. By pigeonholing we may pass to a set $\tilde H \subseteq H$ such that for each bracket-linear form $\psi(h)$ in $\Xi_h^{\struct}$, and for all $h \in \tilde H$, the corresponding tuple always lies in a cube of this form depending only on $\psi$, and not on $h$. 

By Proposition \ref{gowers-prop} there are $\gg_M N^3$ additive quadruples $h_1 + h_2 = h_3 + h_4 \in \tilde H$ and frequencies $\alpha_{h_1,h_2,h_3,h_4},\beta_{h_1,h_2,h_3,h_4} \in \R/\Z$ such that 
\[ |\E_n \chi_{h_1}(n)\chi_{h_2}(n)\overline{\chi_{h_3}(n)\chi_{h_4}(n)}e(\alpha_{h_1,h_2,h_3,h_4} n^2 + \beta_{h_1,h_2,h_3,h_4} n)| \gg_M 1.\]
Now the contribution to this from the structured parts, 
\[ \chi_{h_1}^{\struct}(n)\chi_{h_2}^{\struct}(n)\overline{\chi_{h_3}^{\struct}(n)\chi_{h_4}^{\struct}(n)},\] is a product of bracket quadratic terms of the form 
\[ e(\psi(h_1) n \lfloor \theta n\rfloor + \psi(h_2) n \lfloor \theta n\rfloor - \psi(h_3)n \lfloor \theta n\rfloor - \psi(h_4)n \lfloor \theta n\rfloor)\] 
or
\[ e(\theta n \lfloor \psi(h_1) n\rfloor + \theta n \lfloor \psi(h_2) n\rfloor  - \theta n \lfloor \psi(h_3) n\rfloor - \theta n \lfloor \psi(h_4) n\rfloor).\]
For the quadruples $h_1,h_2,h_3,h_4$ under consideration we have $\psi(h_1) + \psi(h_2) = \psi(h_3) + \psi(h_4)$, and so the first of these expressions is identically 1. The second is not, but by applying Lemma \ref{brack-identities} we see that it is merely a combination of terms of the form $e(\theta n^2)$, $e(\{\alpha n\}\{\beta n\})$ and $e(\alpha n\lfloor \beta n \rfloor)$ with $\Vert \beta \Vert_{\R/\Z} = O_M(1/N)$, where $\alpha,\beta$ and $\theta$ depend on $h_1,h_2,h_3,h_4$. Applying Lemma \ref{1-step-approx}, it follows that we may completely ignore the contribution from these structured parts, although we may need to modify the frequencies $\alpha_{h_1,h_2,h_3,h_4},\beta_{h_1,h_2,h_3,h_4}$.\endproof

The next task is to use a similar (but much simpler) argument to that used for Lemma \ref{h-dep-key} to study the conclusion of Lemma \ref{gowers-unstruct} for a particular quadruple $h_1 + h_2 = h_3 + h_4$.

\begin{lemma}\label{simpler-nildistribution-lem}
Let $h_1,h_2,h_3,h_4$ be fixed and suppose that nilcharacters $\chi_{h_j}(n)$ have the form \eqref{chi-form-2}. Suppose that for each $j = 1,2,3,4$ the underlying frequency set $\Xi_{h_j}$ is decomposed as $\Xi_* \cup \Xi_{h_j}^{\struct} \cup \Xi_{h_j}^{\unstruct}$, where the core set $\Xi_*$ does not depend on $h_j$ and each element of $\Xi_{h_j}^{\struct}$ is a bracket linear form $\psi(h_j)$, again not depending on $h_j$. Suppose that 
\[ |\E_{n \in [N]} \chi^{\unstruct}_{h_1}(n) \chi^{\unstruct}_{h_2}(n)\overline{\chi_{h_3}^{\unstruct}(n)\chi_{h_4}^{\unstruct}(n)}e(\alpha_{h_1,h_2,h_3,h_4} n^2 + \beta_{h_1,h_2,h_3,h_4} n) |\geq 1/M.\]
 Suppose that not all of the integers $m_{[i,i']}(h_j)$ corresponding to frequencies $\xi_{h_j,i},\xi_{h_j,i'}$, one of which is in $\Xi_{h_j}^{\unstruct}$, vanish. Then some there is some $O_{M}(1)$-rational relation, up to $O_{M}(1/N)$, amongst the elements of $\Xi_* \cup \Xi^{\unstruct}_{h_1} \cup \Xi^{\unstruct}_{h_2} \cup \Xi^{\unstruct}_{h_3} \cup \Xi^{\unstruct}_{h_4}$.
\end{lemma}
\proof Once again we interpret the assumption as an assertion about a 2-step nilsequence. Perhaps the ''correct'' way to do this (and the manner more amenable to generalisation) would be to mimic the construction of the last section and apply the Quantitative Ratner theorem once again. However in the special case of the $U^4$-norm being addressed by this paper a shortcut in which only the (simpler) quantitative Leibman dichotomy, Theorem \ref{quant-leib}, is needed and we give this here. Let us take $G$ to be the free 2-step nilpotent Lie group on the ordered generating set $\{e_{\xi} : \xi  \in \Xi_* \cup \Xi^{\unstruct}_{h_1} \cup \Xi^{\unstruct}_{h_2} \cup \Xi^{\unstruct}_{h_3} \cup \Xi^{\unstruct}_{h_4}\}$. As in \S \ref{sec3} we identify the ``coordinate'' functions $F_{\xi,\xi'} : G/\Gamma \rightarrow \C$, and we take a polynomial sequence $g : \Z \rightarrow G$ whose coordinate at $e_{\xi}$ is $\xi n$, for all $\xi$ in the above indexing set, and all of whose other coordinates are zero except for that at $[e_{\xi},e_{\xi'}]$ for some arbitrary pair of frequencies $\xi,\xi' $ in the above set, where the coordinate of $g$ is some quadratic $q = q_{h_1,h_2,h_3,h_4}(n)$ to be specified shortly. Inside $G$ take $\Gamma$ to be the lattice of integer points in the free 2-step nilpotent Lie group. Finally, take 
\[ F := \prod_{i=1}^4 F_{\xi_{*,j},\xi_{h_i,j'}}^{m_{j,j'}(h_i)}.\]  By an appropriate choice of the quadratic term $q$ we may ensure that
\[ F(g(n)\Gamma) = \chi^{\unstruct}_{h_1}(n) \dots \chi^{\unstruct}_{h_4}(n) e(\alpha_{h_1,h_2,h_3,h_4} n^2 + \beta_{h_1,h_2,h_3,h_4} n).\]
Note that we have $\int_{G/\Gamma} F = 0$. Although $F$ is only piecewise Lipschitz, it is nonetheless the case that if $(g(n)\Gamma)_{n \in P}$ is $\delta$-equidistributed for an appropriate $\delta \gg_M$ then $|\E_{n \in [N]} F(g(n)\Gamma)| \leq 1/10M$, contrary to assumption. This is because, as in the last section, we may decompose $F$ as a sum $F_0 + F_1$ where $\Vert F_0 \Vert_{\Lip} = O_{M,\eps}(1)$ and $|F_1|$ is bounded above pointwise by a function $\Psi$ with $\int_{G/\Gamma} \Psi = O(\eps)$ and $\Vert \Psi\Vert_{\Lip} = O_{M,\eps}(1)$. 

Thus we are forced to conclude that $(g(n)\Gamma)_{n \in P}$ is not $\delta$-equidistributed on $G/\Gamma$, for some $\delta \gg_M 1$. By the quantitative Leibman dichotomy, Theorem \ref{quant-leib}, this implies that there is some $k \in \Z^{\dim (G : [G,G])}$, $0 < |k| = O_{M}(1)$, such that $\Vert k \cdot (\pi \circ g)\Vert_{C^{\infty}[N]} = O_{M}(1/N)$. In view of the way that $\pi \circ g$ was constructed, namely the fact that the horizontal part $\pi \circ g$ contains only the terms $\xi n$ with $\xi \in \Xi_* \cup \Xi^{\unstruct}_{h_1} \cup \Xi^{\unstruct}_{h_2} \cup \Xi^{\unstruct}_{h_3} \cup \Xi^{\unstruct}_{h_4}$, this is precisely the result claimed.\endproof

The conclusion of Lemma \ref{simpler-nildistribution-lem} looks rather weak, but using the tools of additive combinatorics pioneered in this context by Gowers (particularly in \cite[Ch. 7]{gowers-longaps}) it turns out to be enough for us to be able to impose some bracket linear behaviour on some of the unstructured sets $\Xi_h^{\unstruct}$. The following result concerning approximate homomorphisms is our key tool. We know of no source for this precise result in the literature, though we feel it should be somehow be regarded as ``known''. It is appropriate to associate the names of Fre\u{\i}man, Ruzsa and Gowers with results of this kind.

\begin{proposition}[Approximate homomorphisms]\label{approx-hom-prop}
Let $\delta,\eps \in (0,1)$ be parameters and suppose that $f_1,f_2,f_3,f_4 : S \rightarrow \R/\Z$ are functions defined on some subset $S \subseteq [N]$ such that there are at least $\delta N^3$ quadruples $(x_1,x_2,x_3,x_4) \in S^4$ with $x_1 + x_2 = x_3 + x_4$ and $\Vert f_1(x_1) + f_2(x_2) - f_3(x_3) + f_4(x_4)\Vert_{\R/\Z} \leq \eps$. Then there is a bracket linear phase $\psi : \Z \rightarrow \R/\Z$ of complexity $O_{\delta}(1)$ and a set $S' \subseteq S$, $|S'| \gg_{\delta} N$, such that $f_1(x) = \psi(x) + O(\eps)$ for all $x \in S'$.
\end{proposition}
\proof See Appendix \ref{approx-hom-app}.\endproof

\begin{lemma}\label{bsg-app}
Let $H \subseteq [N]$ be a set of size at least $N/M$, and suppose that we have a core set $\Xi_*$ and, for each $h \in H$, sets $\Xi^{\unstruct}_h$. Suppose that  $|\Xi_*|, |\Xi^{\unstruct}_h| \leq M$.  Suppose that for at least $N^3/M$ additive quadruples $h_1 + h_2 = h_3 + h_4$ in $H$ there is an $M$-linear relation, up to $O(M/N)$, in $\Xi_* \cup \Xi^{\unstruct}_{h_1} \cup \Xi^{\unstruct}_{h_2} \cup \Xi^{\unstruct}_{h_3} \cup \Xi^{\unstruct}_{h_4}$. Then either
\begin{enumerate}
\item There is some element of the core $\Xi_*$ which lies in the $O_M(1)$-span of the others, up to $O_M(1/N)$, or
\item There is a bracket linear form $\psi$ of degree $O_M(1)$ and a set $H' \subseteq H$, $|H'| \gg_M |H|$, such that $\psi(h)$ lies in the $O_M(1)$-linear span up to $O_M(1/N)$ of $\Xi^{\unstruct}_h$ for all $h \in H'$.
\end{enumerate}
\end{lemma}
\proof Let the elements of the core set $\Xi_*$ be $\{\xi_{*,1},\dots,\xi_{*,M}\}$ and those of the petal set $\Xi_h$ be $\{\xi_{h,1},\dots,\xi_{h,M}\}$. Suppose, for a given quadruple $h_1 + h_2 = h_3 + h_4$, that the approximate linear relation between the elements of $\Xi_* \cup \Xi_{h_1} \cup \Xi_{h_2} \cup \Xi_{h_3} \cup \Xi_{h_4}$ is
\begin{align*} 
& \Vert q_{*,1}(h_1,h_2,h_3,h_4) \xi_{*,1} + \dots + q_{*,M}(h_1, h_2, h_3, h_4) \xi_{*,M} \\ 
& + q_{1,1}(h_1,h_2,h_3,h_4) \xi_{h_1,1} + \dots + q_{1,M}(h_1,h_2,h_3,h_4) \xi_{h_1,M} \\ & + \dots + q_{4,1}(h_1,h_2,h_3,h_4) \xi_{h_4,1} + \dots + q_{4,M} (h_1,h_2,h_3,h_4)\xi_{h_4, M} \Vert_{\R/\Z}  = O(M/N),
\end{align*}
where each integer $q$ has magnitude at most $M$. There are only $(2M + 1)^{5M}$ choices for these integers and so we may pass to a subcollection of $\gg_M N^3$ quadruples for which there is such a relation and for which none of the $q$'s depend on $h_1,h_2,h_3,h_4$. Since $\Xi_*$ is $M$-dissociated, at least one of the $q_{i,j}$ must be nonzero, $i = 1,2,3,4$; without loss of generality, suppose that $q_{1,1} \neq 0$. 

Writing $f_1(h_1) := q_{1,1} \xi_{h_1,1} + \dots + q_{1,M} \xi_{h_1, M}$,  we see that we have found functions $f_2,f_3,f_4 : H \rightarrow \R/\Z$ such that 
\[\Vert f_1(h_1) + f_2(h_2) - f_3(h_3) - f_4(h_4) \Vert_{\R/\Z} \leq 1/N'\] for $\gg_M N^3$ additive quadruples $h_1 + h_2 = h_3 + h_4 \in H$, for some $N' \gg N/M$.  Now we apply Proposition \ref{approx-hom-prop} to conclude that there is a bracket linear phase $\psi$ of complexity $O_M(1)$ such that $f_1(h) = \psi(h) + O_{M}(1/N)$ for all $h$ in some set $H' \subseteq H$, $|H'| \gg_M N$. This concludes the proof of the lemma.
\endproof

We are now in a position to prove Proposition \ref{approx-lin*} which, recall, was the inductive step driving the main result of this section, namely Proposition \ref{approx-lin}. The argument is very similar to that employed in the proof of Proposition \ref{prop-h-dep}, hingeing on repeated use of the bracket identities of Lemma \ref{brack-identities} to expand out linear combinations of frequencies.\vspace{11pt}

\emph{Proof of Proposition \ref{approx-lin*}.} The assumption that $|\E_{n \in P} \Delta_hf(n)\overline{\chi_h(n)}| \geq 1/M$ may be fed into Lemma \ref{gowers-unstruct} to conclude the existence of a set $H' \subseteq H$ with $|H'| \gg_M |H|$ such that 
\[ \E_n\chi^{\unstruct}_{h_1}(n)\chi^{\unstruct}_{h_2}(n)\chi^{\unstruct}_{h_3}(n)\chi^{\unstruct}_{h_4}(n)e(\alpha_{h_1,h_2,h_3,h_4} n^2 + \beta_{h_1,h_2,h_3,h_4} n)) \gg_M 1\] for $\gg_M N^3$ additive quadruples $h_1 + h_2 = h_3 + h_4$ in $H'$. This in turn may be fed into Lemma \ref{simpler-nildistribution-lem}, which allows us to conclude that for each of these additive quadruples there is an $O_M(1)$ linear relation, up to $O_M(1/N)$, between the elements of $\Xi_* \cup \Xi^{\unstruct}_{h_1} \cup \Xi_{h_2}^{\unstruct} \cup \Xi_{h_3}^{\unstruct} \cup \Xi_{h_4}^{\unstruct}$. There is one other possibility here, namely that in the attempt to apply Lemma \ref{simpler-nildistribution-lem} we find that, for many quadruples $h_1 + h_2 = h_3 + h_4$,  all of the integers $m_{[i,i']}(h_j)$ corresponding to frequencies $\xi_{h_j,i}, \xi_{h_j,i'}$, one of which is in $\Xi_{h_j}^{\unstruct}$, are zero. This is a rather trivial case, however, for we may then pass to the set $\tilde H$ of $h_1$ (say) appearing here, and simply delete the unstructured frequencies $\Xi_h^{\unstruct}$, which play no actual role in the expression for $\chi_h(n)$. The conclusion of Proposition \ref{approx-lin*} is then immediate in this case.

Returning to the main line of the argument, we may then apply Lemma \ref{bsg-app} to conclude that either 
\begin{enumerate}
\item There is some element $\xi \in \Xi_*$ which lies in the $O_M(1)$-linear span of the others, up to $O_M(1/N)$, or
\item There is a bracket linear form $\psi$ of degree $O_M(1)$ and a set $\tilde H \subseteq H'$, $|\tilde H| \gg_M |H|$, so that $\psi(h)$ lies in the $O_M(1)$-linear span of $\Xi^{\unstruct}_h$ for all $h \in \tilde H$.
\end{enumerate}
These two possibilities will correspond to alternatives (i) and (ii) respectively in Proposition \ref{approx-lin*}. To see this we proceed rather as in the proof of Proposition \ref{prop-h-dep}, making use once again of Lemma \ref{brack-identities} as well as extensive use of Lemma \ref{1-step-approx} to handle the somewhat annoyingly non-Lipschitz 1-step objects which arise. The treatment of (i) is exactly analogous to the aforementioned argument, so we only describe (ii) in any detail. 

Assume that the sets $\Xi^{\unstruct}_h$ are ordered as $\xi_{h,k_0+1},\dots,\xi_{h,k}$. We are assuming that there is a bracket-linear form $\psi(h)$ having the form $q_{h,k_0 + 1} \xi_{h,k_0 + 1} + \dots + q_{h,k} \xi_{h,k} + O_M(1/N)$, for all $h \in H'$. Here the integers $q_{h,j}$ are all bounded in magnitude by $O_M(1)$ and so we may, by passing to a further subset $H'' \subseteq H'$, assume that they do not depend on $h$. Without loss of generality let us suppose that $q_{h,k} \neq 0$. Then we may write $\xi_{h,k}$ as an $O_M(1)$-linear combination of $\frac{1}{Q}\psi(h)$ and the frequencies $\frac{1}{Q} \xi_{h,k_0 + 1},\dots,\frac{1}{Q}\xi_{h,k-1}$, plus an error of $O_M(1/N)$, where $Q$ is the lcm of the numbers up to $O_M(1)$. Now we replace $\Xi_h^{\struct}$ by $\Xi_h^{\struct} \cup \{\frac{1}{Q}\psi(h)\}$ and $\Xi_h^{\unstruct}$ by $\{\frac{1}{Q} \xi_{h,k_0 + 1},\dots,\frac{1}{Q}\xi_{h,k-1}\}$, and then proceed to rewrite the bracket quadratics $e(\xi n\lfloor \xi' n\rfloor)$ using these new sets of frequencies by means of Lemma \ref{brack-identities} and Lemma \ref{1-step-approx} exactly as we did at the end of \S \ref{step1-sec}.\endproof

Before moving onto the next section we apply one additional piece of analysis to Proposition \ref{approx-lin}. This allows us to conclude that the quadratic frequency $\alpha_h$ varies bracket-linearly in $h$ as well. Thus, once this is done, only the linear term $e(\beta_h n)$ does not have a rigid structure imposed upon it. 

\begin{proposition}\label{approx-lin-new}
In the statement of Proposition \ref{approx-lin}, we may assume that the quadratic frequency $\alpha_h$ varies bracket-linearly in $h$.
\end{proposition}
\proof We may, of course, take for granted the conclusion of Proposition \ref{approx-lin}.  We apply Proposition \ref{gowers-prop} once again, using the same argument we employed at the start of the proof of Lemma \ref{gowers-unstruct} to first pass to a subset $H' \subseteq H$, $|H'| \gg_{M} N$, on which all the bracket linear forms $\psi$ in the petals $\Xi'_h$ are linear in the sense that $\psi(h_1) + \psi(h_2) = \psi(h_3) + \psi(h_4)$ whenever $h_1 + h_2 = h_3 + h_4$ with $h_1,h_2,h_3,h_4 \in H'$. This gives
\[ \E_{n \in [N]} \chi_{h_1}(n)\chi_{h_2}(n + h_1 - h_4)\chi_{h_3}(n)\chi_{h_4}(n + h_1 - h_4)\gg_M 1.
\]

As in Corollary \ref{gowers-cor}, this implies that $\chi_{h_1}(n)\chi_{h_2}(n)\chi_{h_3}(n)\chi_{h_4}(n)$ correlates with a quadratic phase $e(\alpha_{h_1,h_2,h_3,h_4} n^2 + \beta_{h_1,h_2,h_3,h_4}n)$. Moreover a careful analysis of the proof of that corollary, looking at the decomposition $\chi_h(n) = \chi'_h(n) e(\alpha_h n^2 + \beta_h n)$, where 
\[ \chi'_h(n) = \prod_{1 \leq i < i' \leq k} F_{[i',i]}^{m_{[i',i](h)}}(g_h(n)\Gamma),\] reveals that we can take $\alpha_{h_1,h_2,h_3,h_4} = \alpha_{h_1} + \alpha_{h_2} - \alpha_{h_3} - \alpha_{h_4}$. That is, the genuinely bracket-quadratic objects comprising $\chi'_h(n)$ only give rise to \emph{linear} terms. 

The term $\chi'_{h_1}(n)\chi'_{h_2}(n)\chi'_{h_3}(n)\chi'_{h_4}(n)$ arising from the genuinely bracket quadratic parts is a product of terms of the form $e(\alpha n\lfloor \psi(h_1) n\rfloor + \alpha n\lfloor \psi(h_2) n\rfloor - \alpha n \lfloor \psi(h_3) n\rfloor - \alpha n \lfloor \psi(h_4) n\rfloor)$ where, recall, $\psi(h_1) + \psi(h_2)  = \psi(h_3) + \psi(h_4)$. Using Lemma \ref{brack-identities} (iii) to move the $\psi$ terms to the outside of the brackets and applying Lemma \ref{1-step-approx} repeatedly, we conclude that 
\[ \E_{n \in [N]} e((\alpha'_{h_1} + \alpha'_{h_2} - \alpha'_{h_3} - \alpha'_{h_4})n^2 + \theta_{h_1,h_2,h_3,h_4} n) \gg_M 1\] for all these quadruples $h_1 + h_2 = h_3 + h_4$, where $\alpha'_h - \alpha_h$ is a bracket-linear form of complexity $O_M(1)$.
By Lemma \ref{weyl-application} it follows that there is some $q = O_{M}(1)$ such that 
\[ \Vert q(\alpha'_{h_1} + \alpha'_{h_2} - \alpha'_{h_3} - \alpha'_{h_4})\Vert_{\R/\Z} = O_{M}(1/N^2).\]
By Proposition \ref{approx-hom-prop} there is a further subset $H'' \subseteq H$, $|H''| \gg_{M} |H|$, together with a bracket linear form $\psi'(h)$ of complexity $O_{M}(1)$, such that 
\[ q \alpha'_h = \psi'(h) + O_{M}(1/N^2)\] for all $h \in H''$. This means that 
\[ \alpha_h = \psi''(h) + \frac{r_h}{q} + O_M(1/N^2),\] where $\psi''(h)$ is another bracket linear form and $r_h$ takes integer values. 
Refining $[N]$ into progressions of common difference $q$ and length $\gg_M N$ small enough to make the  $O_{M}(1/N^2)$ error negligible, and then applying Lemma \ref{1-step-approx} (ii), we obtain the claim.\endproof

\label{2-step-sec}

\section{Step 3: The symmetry argument}\label{sym-section}

Finally we turn to \emph{Step 3} of the programme outlined in \S \ref{inductive-sec}, the so-called \emph{symmetry argument}. Here we shall take an approach somewhat different to the one we shall employ in the general case of the $U^{s+1}$-norm, $s \geq 4$, where further use is made of the nilmanifold distribution results of \S \ref{nil-dist-sec} and there are slightly complicated issues concerning the keeping-track of the complexity of various bracket expressions.

In the special case of the $U^4$-norm that this paper is concerned with, a rather direct argument using Bohr sets is possible.  Let $S = \{\theta_1,\dots,\theta_d\} \subseteq \R/\Z$ be a set of frequencies and suppose that $\rho \in (0,1)$. Then we set \[ B(S,\rho,N) := \{n \in [\rho N] : \Vert n \theta_j \Vert_{\R/\Z} \leq \rho \quad \mbox{for all $j = 1,\dots,d$}\}.\]
We shall need a small amount of the theory of such sets, particularly pertaining to the notion of \emph{regularity} -- the idea that there is a plentiful supply of $\rho$ for which the size of $B(S,\rho',N)$ is nicely controlled for $\rho' \approx \rho$. The need to introduce this idea in additive combinatorics was first appreciated in \cite{bourgain-3aps} and it has now appeared in several places, for example \cite{green-tao-u3inverse} where the notion is defined in Definition 2.6 and discussed in more detail in Chapter 8.

For our purposes here we say that a value $\rho$ is \emph{regular} if we have
\[ |B(S,(1 +\kappa)\rho,N)| = |B(S,\rho,N)| (1 + O(d|\kappa|))\] uniformly for $|\kappa| \leq 1/d$. We shall need the following facts about regular Bohr sets. It would be possible to obtain much more precise statements but we shall not need to do so here.

\begin{lemma}[Regular Bohr sets -- Basic Facts]\label{basic-bohr-fact}
Fix a set $S = \{\theta_1,\dots,\theta_d\}$ of frequencies, and write $B := B(S,\rho,N)$. We have the following facts.
\begin{enumerate}
\item \textup{(Ubiquity of regular values)} For any $\rho_0 \in (0,1/2)$ there is a regular value of $\rho$ in the interval $[\rho_0,2\rho_0]$.
\item \textup{(Fourier expansion of Bohr cutoffs)} Suppose that $\rho$ is regular, and that $\eps > 0$ is a parameter. Then we may decompose the cutoff $1_B(n)$ as $\psi_1(n) + \psi_2(n)$, where $\psi_1(n) = \int^1_0 \widehat{\psi}_1(\theta) e(\theta n)\, d\theta$ with $\Vert \widehat{\psi}_1 \Vert_1 := \int^1_0 |\widehat{\psi}_1(\theta)|\, d\theta \leq C_{\epsilon,\rho}$ and $\sum_n |\psi_2(n)| \leq \eps N$.
\item \textup{(Large generalised Fourier coefficients)} Suppose that $\rho$ is regular and that $\phi : B(S,2\rho,N) \rightarrow \R/\Z$ is locally linear on $B$ in the sense that $\phi(x+y) = \phi(x) + \phi(y)$ whenever $x,y \in B$. Suppose that $|\E_{x \in B} e(\phi(x))| \geq \eta$. Then there is a  regular value of $\rho'$, $\rho' \gg_{\eps,\eta,\rho} 1$, such that $\Vert \phi(x)\Vert_{\R/\Z} \leq \eps$ for all $x \in B(S,\rho')$.
\end{enumerate}
\end{lemma}
\emph{Sketch Proof.} The definition of Bohr set we are using here is very slightly different to that used in \cite{green-tao-u3inverse}, in that our Bohr sets are contained in $[N]$ and not in $\Z/N\Z$. Nonetheless, the proofs of the above statements are so close to those in $\Z/N\Z$ that we simply refer to the relevant sections of the aforementioned paper. Statement (i) is \cite[Lemma 8.2]{green-tao-u3inverse}. Statement (ii) is not explicitly mentioned in \cite{green-tao-u3inverse}. To prove it, take $\psi_1(n) = \frac{1}{|B'|}1_B \ast 1_{B'}(n)$, where $B' := B(S,\rho',N)$ for a suitably small $\rho' \gg_{\eps,\rho} 1$. The bound on $\Vert \widehat{\psi}_1\Vert_1$ follows from Plancherel, whilst the bound on $\Vert \psi_2\Vert_1$ is a consequence of the regularity of $B$ and the observation that $\frac{1}{|B'|}1_B \ast 1_{B'}(n) = 1_B(n)$ provided that $n \notin B(S,\rho + \rho',N) \setminus B(S,\rho - \rho', N)$. Finally, (iii) is \cite[Lemma 8.4]{green-tao-u3inverse}.
\endproof

Let us return to the main business of this section, which is to conclude the proof of Theorem \ref{mainthm}. The main result of the last section, Proposition \ref{approx-lin-new}, took us from the assumption that $\Vert f \Vert_{U^4} \geq \delta$ to the conclusion that 
\begin{equation}\label{der-cor}
| \E_{n \in [N]} f(n) \overline{f(n+h)}\chi_h(n)| \gg_{\delta} 1
\end{equation}
for a set $H$ of size $\gg_{\delta} N$, where $\chi_h(n)$ is a product of terms of the form $e(\{\alpha h\} \beta n\lfloor \gamma n\rfloor)$, $e(\alpha\{\beta h\} n^2)$ and $e(\theta_h n)$.  Using the fact that $\lfloor \gamma n\rfloor = \gamma n - \{\gamma n\}$, we may assume that 
\[ 
\chi_h(n) = e(\sum_{j=1}^k \{\alpha_j h\} \beta_j n\{\gamma_j n\}) e(\sum_{j=1}^k \alpha'_j \{\beta'_j h\} n^2) ) e(\theta_h n).
\]
Later on it will be convenient to assume that 
\begin{equation}\label{star-5} \mbox{For all $h \in H$ we have $\Vert \theta h\Vert_{\R/\Z} \geq \rho_1$ for all $\theta \in \{\alpha_1,\dots,\alpha_k,\beta'_1,\dots,\beta'_k\}$},\end{equation} for some small parameter $\rho_1 > 0$ to be specified later. This can be achieved at the expense of thinning out $H$ somewhat to a set of size merely $\gg_{\rho_1,\delta} N$, as we now show.

To demonstrate the last claim we distinguish two types of such $\theta$. We say that $\theta$ is \emph{good} if the number of $h \in H$ such that $\Vert \theta h \Vert_{\R/\Z} < \rho_1$ is at most $10\rho_1 N$. By refining $H$ to a set $H' \subseteq H$ with $|H'| \geq |H| - 20\rho_1 k N$, we may assume that $\Vert \theta h\Vert_{\R/\Z} \geq \rho_1$ for all $h \in H'$ and for all good $\theta$. Note that $|H'| \geq |H|/2$ if $\rho_1$ is chosen small enough as a function of $\delta$, as it will be later on. If $\theta$ is not good then the sequence $\{n\theta(\text{mod } \Z)\}_{n \in [N]}$ is not $\rho_1$-equidistributed, and by well-known results of diophantine approximation (see, for example, \cite[Proposition 3.1]{green-tao-nilratner}) there is some $q \ll \rho_1^{-C}$ such that $\Vert q \theta \Vert_{\R/\Z} \ll \rho_1^{-C}/N$.  This means that the bracket $\{\theta h\}$ takes on only $\rho_1^{-2C}$ values as $h$ ranges over $[N]$, and so there is a subset $H'' \subseteq H'$, $|H''| \gg \rho_1^{Ck} |H|$, on which all these brackets are constant. This means that the corresponding terms in $\chi_h(n)$ may be ignored, for the purpose of \eqref{der-cor}, since they depend just on $n$ and not on $h$. Replacing $H$ by $H''$ gives the claim, and henceforth we assume that \eqref{star-5} holds, remembering that we now only have the weaker bound $|H| \gg_{\rho_1,\delta} N$.

Write 
\[ 
T(x,y,z) := \sum_{j=1}^k \{\alpha_j x\} \frac{\beta_j}{3} y \{\gamma_j z\} + \sum_{j=1}^k \frac{\alpha'_j}{3} \{\beta'_jx\} yz ,
\]
so that $3T(h,n,n)$ is the form appearing in the definition of $\chi_h(n)$. Here, there are three possible choices for each $\beta_j/3,\alpha'_j/3$ and it does not matter which we take; the reason for introducing these 3's will become apparent later. Then $T(x,y,z)$ is trilinear on the Bohr set $B := B(S,\rho_0,N)$, where $S = \{\alpha_1,\dots,\alpha_k,\gamma_1,\dots,\gamma_k,\beta'_1,\dots,\beta'_k\}$ and the parameter $\rho_0 \in [\frac{1}{20},\frac{1}{10}]$ is chosen so that $B$ is regular.  By stating that $T$ is trilinear we mean that, for example, $T(x_1+x_2,y,z) = T(x_1,y,z) + T(x_2,y,z)$ when all of $x_1,x_2,x_1+x_2,y,z$ lie in $B$.  We begin by symmetrising $T$ in the last two variables, a straightforward task. For each $j$ pick some $\tilde\beta_j$ such that $2\tilde\beta_j = \beta_j/3$ (there are two choices) and set 
\[ 
\tilde T(x,y,z) := \sum_{j=1}^k \{\alpha_j x\} \tilde\beta_j y \{\gamma_j z\} + \sum_{j=1}^k\{\alpha_j x\} \tilde\beta_j z \{\gamma_j y\} + \sum_{j=1}^k \frac{\alpha'_j}{3} \{\beta'_jx\} yz.
\]
Then of course $\tilde T(h,n,n) = T(h,n,n)$, but now $\tilde T(x,y,z)$ is symmetric in the last two variables. Dropping the tildes, we assume henceforth that $T$ itself is symmetric in the last two variables.

Our assumption, then, is that 
\[ 
|\E_{n \in [N]} f(n) \overline{f(n+h)} e(3T(h,n,n)) e(\theta_h n)| \gg_{\delta} 1
\] for all $h$ lying in some set $H$ of size at least $\gg_{\rho_1,\delta} N$, where $H$ additionally satisfies \eqref{star-5}. Our immediate goal is to localize the variables $h$ and $n$ to small Bohr sets so that we may properly exploit the trilinearity of $T$.

Let us briefly reprise the heuristic mentioned in the \S \ref{inductive-sec}  to recall why it is that we expect $T$ to be symmetric in the first two coordinates as well (on a``nice set''). 
Suppose we knew that  $f(n) \overline{f(n+h)} =  \chi_h(n)=e(3T(h,n,n))e(\theta_h n)$ for all $n,h$. Then we get 
\[
\chi_h(n+k)\chi_k(n)= f(n) \overline{f(n+h+k)}= \chi_k(n+h)\chi_h(n).
\]
Using the trilinearity of $T$ and symmetry in the last two coordinates we get \[ 6T(h,k,n)=  6T(k,h,n).\]
Now of course we do not have proper equations but only correlations, we don't have correlation for all $h$
but only for "many", and we have trilinearity only when the variables are  restricted to Bohr sets,
so we must work much harder.

We start with the $h$ variable. Set $B_1 := B(S,\rho_1,N)$, where $\rho_1$ is the (as yet unspecified) quantity appearing in \eqref{star-5}. Modifying $\rho_1$ by at most a factor of two, we may assume that $B_1$ is regular. We claim that it is possible to find an $h_0 \in H$ such that the intersection $H \cap (h_0 + B_1)$ has size $\gg_{\rho_1,\delta} N$. A slight trick is necessary to establish this: consider
\[ \sum_{n \in [N]} 1_H \ast 1_{B'} \ast 1_{B'}(n) 1_H(n),\] where $B' := B(S,\rho_1/2,N)$.  On the one hand this equals $\sum_n 1_H \ast 1_{B'}(n)^2$ which, by the Cauchy-Schwarz inequality, is $\gg_{\rho_1,\delta} N^3$. On the other hand we have $1_{B'} \ast 1_{B'}(n) \leq |B_1| 1_{B_1}(n)$ for $n \in [N]$, and from these two inequalities the claim follows immediately.

Our assumption now implies that
\[ |\E_{n \in [N]} f(n) f(n + h_0 + h') e(3T(h_0+ h',n,n)) e(\theta_{h_0 + h'}n)| \gg_{\delta} 1\] for all $h' $ lying in some set $H' \subseteq B_1 = B(S,\rho_1,N)$, $|H'| \gg_{\rho_1,\delta} N$. 
By the careful construction of $H$ (cf. \eqref{star-5}) and the fact that $h_0 \in H$ we have $\{\alpha_j (h_0 + h')\} = \{\alpha_j h_0\} + \{\alpha_j h'\}$, and similarly for the $\beta'_j$, and hence we obtain the linearity property $T(h_0 + h',n,n) = T(h_0,n,n) +  T(h',n,n)$. After relabelling we hence have
\[ |\E_{n \in [N]} f_1(n) f_2(n+h) e(3T(h,n,n)) e(\theta_h n)| \gg_{\delta} 1\] for all $h \in H$, where $H \subseteq B_1$, $|H| \gg_{\rho_1,\delta} N$, $f_1(n) := f(n)e(T(h_0,n,n))$ and $f_2(n) := f(n+h_0)$.

We must now localise the $n$ variable, and for this we use a somewhat different trick. By averaging there is some $n_0$ such that 
\[ \E_{n \in [N]} f_1(n_0 + n) f_2(n_0 + n + h) e(3T(h, n_0 + n, n_0 + n)) e(\theta_h n)1_{B_1}(n) \gg_{\delta} 1.\]
Now we have
\begin{align*} \beta & (n+n_0)\{\gamma (n+n_0)\} \\ &= \beta n_0 \{ \gamma(n + n_0)\} + \beta n \{\gamma n\} + \beta n \{\gamma n_0\} + \beta n (\{\gamma (n + n_0)\} - \{\gamma n\} - \{\gamma n_0\}).\end{align*}
Substituting into the expression for $e(3T(h, n+n_0,n+n_0))$ and expanding, we see that the contribution from the term $e(\beta n_0\{\gamma(n + n_0)\})$ may be absorbed into the linear term $e(\theta_h n)$ (by Lemma \ref{1-step-approx}), as may the term $e(\beta n\{\gamma n_0\})$ (trivially). The term $\{\gamma (n +n_0)\} - \{\gamma n\} - \{\gamma n_0\}$ takes values in $\{-1,0,1\}$ according to whether $\gamma n \md 1$ lies in certain intervals $I_{\gamma}^{-1},I_{\gamma}^0,I_{\gamma}^{+1}$, and so we obtain
\begin{align*}\E_{n \in [N]} & f_1'(n)f'_2(n + h)   e(3T(h,n,n)) e(\theta_h n) \times \\ & \times \prod_{j=1}^k \big(1_{\gamma_j n \in I^{-1}_{\gamma_j}} e(  -\{\alpha_j h\}\beta_j n  ) + 1_{\gamma_j n \in I^0_{\gamma_j}} + 1_{\gamma_j n \in I^{+1}_{\gamma_j} }e( \{\alpha_j h\}\beta_j n  )\big) 1_{B_1}(n)\gg 1,\end{align*}
where $f'_1(n) = f_1(n + n_0)$ and $f'_2(n) = f_2(n  + n_0)$.
It follows that there is a choice of $\eps_j \in \{-1,0,1\}$ and a $\theta'_h$ such that 
\[ \E_{n \in [N]} f_1'(n) f'_2(n+h) e(3T(h,n,n)) e(\theta'_h n) \prod_{j = 1}^k 1_{\gamma_j n \in I^{\eps_j}_{\gamma_j}}1_{B_1}(n) \gg 1.\]
By Lemma \ref{1-step-approx} we may remove the last term at the expense of changing $\theta'_h$ again. Removing the dashes for notational convenience we now obtain\[ \E_{n \in [N]} f_1(n) f_2(n+h) e(3T(h,n,n)) e(\theta_h n) 1_{B_1}(n)\gg 1.\]
Here, $f_1(n) = f(n + n_0) e(T(h_0, n_0 + n, n_0 + n))$ and $f_2(n) = f(n + h_0 + n_0)$, and we recall once more that this is known to hold for $\gg_{\rho_1,\delta} N$ values of $h \in B_1$.

Set $\chi_h(n) := e(3T(h,n,n))e(\theta_h n) 1_{B_1}(n)$. Applying Proposition \ref{gowers-prop}, we obtain
\begin{equation}\label{gowers-cor1} 
\E_{n \in [N]} \chi_{h_1}(n) \chi_{h_2}(n  + h_1 - h_4) \overline{\chi_{h_3}(n) \chi_{h_4}(n + h_1 - h_4)} \gg_{\rho_1,\delta} 1
\end{equation} 
for at least $c_{\rho_1,\delta}N^3$ additive quadruples $(h_1,h_2,h_3,h_4) \in B_1$ with $h_1 + h_2 = h_3 + h_4$.  We have already, in previous sections, extracted ``top order'' information from statements like this and our task here is to exploit the additional structure inherent in \eqref{gowers-cor1}, particularly that present in the terms $h_1 - h_4$.

Parametrising these by $h_1 = h$, $h_2 = h + a+b$, $h_3 = h+ a$, $h_4 = h + b$ we obtain
\[ 
\E_{n \in [N]} \chi_h(n) \chi_{h + a + b}(n+b) \overline{\chi_{h+a}(n)\chi_{h + b}(n+b)} \gg_{\rho_1,\delta} 1
\] 
for at least $c_{\delta}N^3$ triples $h,a,b$ with $h \in t + B'$ and  $a,b \in B(S,3\rho_1,N)$. 
Substituting in the definition of $\chi_h(n)$, and using the trilinearity of $T$ we obtain
\begin{equation}\label{to-deal-with} 
\E_n e(6T(a,b,n)  + (\theta_h  + \theta_{h + a + b} - \theta_{h+a} - \theta_{h+b})n)  1_{B_1}(n)1_{B_1}(n+b)\gg_{\rho_1,\delta} 1
\end{equation} 
for at least $c_{\rho_1,\delta}N^3$ triples $h,a,b$ with $h \in  B_1$ and $a,b \in B(S,3\rho_1,N)$.  Pigeonholing in $h$, one sees that there is some fixed $h$ such that this holds for at least $c_{\rho_1,\delta}N^2$ pairs $a,b \in B(S,3\rho_1,N)$.
Let $\eps = \eps(\rho_1,\delta)$ be a small positive quantity to be specified very shortly. By Lemma \ref{basic-bohr-fact} (ii) and the regularity of $B_1$ we may expand
\[ 
1_{B_1}(n+b) = \int^1_0 \widehat{\psi}_1(\theta) e(\theta (n+b)) \, d\theta + \psi_2(n),
\] 
where $\Vert \widehat{\psi}_1\Vert_1 \leq C_{\eps,\rho_1,\delta}$ and $\sum_n |\psi_2(n)| \leq \eps N$. Choosing $\eps$ so that the contribution to \eqref{to-deal-with} from $\psi_2(n)$ is negligible, we see using the triangle inequality that there is some $\theta \in [0,1]$ such that 
\begin{equation}\label{to-deal-with-2}
\E_{n \in B_1} e(6T(a,b,n)  + (\theta_h  + \theta_{h + a + b} - \theta_{h+a} - \theta_{h+b} + \theta)n) \gg_{\rho_1,\delta} 1 
\end{equation} 
for the same fixed $h$ and many pairs $a,b$ as before.
For each $a,b$ write $\phi_{a,b}(n)$ for the phase appearing here, thus
\[ 
\phi_{a,b}(n) := 6T(a,b,n)  + (\eta_a + \eta'_b + \eta''_{a+b})n
\] 
where
$\eta_a := \theta_h - \theta_{h+a} + \theta$, $\eta'_b := \theta_{h+b}$ and $\eta''_{a + b} := \theta_{h + a + b}$.  Equation \eqref{to-deal-with-2} implies that
\[ 
|\E_{n \in B_1} e(\phi_{a,b}(n))| \gg_{\rho_1,\delta} 1.
\]

Let $\eps = \eps(\delta,\rho_1)$ be a small positive parameter to be specified later.  By Lemma \ref{basic-bohr-fact} (iii) there is some $\rho_2= \rho_2(\eps,\rho_1,\delta)$ such that we have\[ \Vert\phi_{a,b}(n)\Vert_{\R/\Z} \leq \eps\] for all $n \in B_2 := B(S,\rho_2,N)$ and for these same pairs $a,b$, that is to say for at least $c_{\rho_1,\delta} N^2$ pairs $a,b \in B(S,3\rho_1)$. Thus
\[ 
6T(a,b,n) = (\eta_a + \eta'_b + \eta''_{a+b}) n + O(\eps)
\] 
for at least $c_{\rho_1,\delta} N^2$ choices of $a,b \in B(S,3\rho_1,N)$ and for all $n \in B_2$. For at least $c_{\delta} N^3$ triples $a,b,b'$ we thus have
\[ 
6T(a,b-b',n)  = (\eta'_b - \eta'_{b'} + \eta''_{a + b} - \eta''_{a + b'})n + O(\eps)
\] 
for all $n \in B_2$. Writing $c := a + b + b'$ it follows that
\[ 
6T(c - b - b', b-b',n) = (\eta'_b - \eta'_{b'} + \eta''_{c - b'} - \eta''_{c - b})n + O(\eps)
\] 
for at least $c_{\rho_1,\delta}N^3$ triples $c,b,b' \in B(S,9\rho_1,N)$ and for all $n \in B_2$. Fix some $c$ for which this holds for at least $c_{\rho_1,\delta}N^2$ pairs $b,b'$; then by trilinearity of $T$ we have
\[ 
6(T(b,b',n) - T(b',b,n)) =  \kappa_b(n) + \kappa'_{b'}(n) + O(\eps)
\] 
for all these pairs $b,b'$ and for all $n \in B_2$, where $\kappa_b := \eta'_b n - T(c,b,n) + T(b,b,n) - \eta''_{c-b} n$ and $\kappa'_{b'} = -\eta'_{b'}n + \eta''_{c - b'} n - \psi(c,b',n) + \psi(b',b',n)$. The exact form of these expressions is not relevant, as we shall very shortly see.

Indeed for at least $c_{\rho_1,\delta} N^3$ triples $b_1,b_2,b' \in B(S, 3\rho_1)$ we have
\[ 
6T(b_1 - b_2,b',n) - T(b',b_1 - b_2,n) = \kappa_{b_1}(n) - \kappa_{b_2}(n) + O(\eps)
\] 
for all $n \in B_2$, and hence for at least $c_{\rho_1,\delta}N^4$ quadruples $b_1,b_2,b'_1,b'_2 \in B(S,3\rho_1)$ we have
\[ 
6T(b_1 - b_2, b'_1 - b'_2,n) - T(b'_1 - b'_2, b_1 - b_2,n)  = O(\eps)
\] 
for all $n \in B_2$.
There are at least $c_{\rho_1,\delta}N^2$ different pairs $x,y \in B(S,6\rho_1)$ represented as $x = b_1 - b_2, y = b'_1 - b'_2$, and for each of them
\[ 
6(T(x,y,n) - T(y,x,n)) = O(\eps)
\] 
for all $n \in B_2$. 
Write $A \subseteq [N]^2$ for the set of these pairs, thus $|A| \geq c_{\delta} N^2$.  Let us write $A \oplus A$ for the set of all pairs $(x,y_1 \pm y_2)$ where both $(x,y_1)$ and $(x,y_2)$ lie in $A$, together with all pairs $(x_1\pm x_2, y)$ where both $(x_1,y)$ and $(x_2,y)$ lie in $A$. By bilinearity we see that 
\[ 
6(T(x,y,n) - T(y,x,n)) = O(k \eps)
\] 
for all pairs $(x,y)$ in the $k$-fold bilinear sumset $A \oplus A \oplus \dots \oplus A$ and for all $n \in \tilde B$. 

Now by Lemma \ref{bilinear-sarkozy} this $k$-fold bilinear  sumset $A' := A \oplus A \dots \oplus A$ contains a product $P \times P$ provided that $k \geq C_{\delta}$, where $P$ is an arithmetic progression which contains $0$ and has length $N$ and common differences $d = O_{\delta}(1)$. 
Thus for all triples $x,y,z \in P \cap \tilde B$ we have
\begin{equation}\label{approx-sym} 
T(x,y,z) - T(y,x,z) = O(k\eps) + \sigma_{x,y,z},
\end{equation} 
where $\sigma_{x,y,z}$ takes values in $\Z/6\Z$.

Recall that we have
\[ 
|\E_{n\in [N]} f_1(n) f_2(n+h) e(3T(h,n,n)) e(\theta_h n)1_{B_1}(n)| \gg 1
\] 
for many $h \in B_1$. 
By the pigeonhole principle, there are $h_1, n_1$ such that 
\[ 
|\E_{n \in [N]} f_1(n + n_1) f_2(n+h_1 + h) e(3T(h_1 + h,n_1 +n,n_1 + n))1_{P \cap B_1}(n)1_{B_1}(n_1 + n)| \gg 1
\] 
for many $h \in P \cap B_1$. Obviously $n_0 \in B(S, 2\rho_1,N)$, and so we may expand $T(h_1 + h,n_1 + n, n_1 + n)$ using trilinearity. Doing this, absorbing the linear terms into $e(\theta_h n)$ using Lemma \ref{1-step-approx} and expanding the cutoff $1_{B_1}(n_1 + n)$ as a Fourier series using Lemma \ref{basic-bohr-fact} (ii), we obtain
\[ 
|\E_{n \in [N]} f'_1(n) f'_2(n + h) e(3T(h,n,n)) 1_{P \cap B_1}(n)e(\theta'_h n)| \gg 1
\] 
for may $h \in P \cap B_1$. 
Here, $f'_1(n) = f_1(n+n_1) e(T(h_1,n,n)) = f(n + n_0 + n_1) e(T(h_0, n_0 + n_1 + n, n_0 + n_1 + n) + T(h_1,n,n))$ whilst $f'_2(n) = f(n + h_0 + h_1 + n_0 + n_1)$. Once again we drop the dashes in what follows for notational convenience.

By the trilinearity of $T$ and the approximate symmetry \eqref{approx-sym} of $T$ in the first two variables, the genuine symmetry in the last two and another application of Lemma \ref{1-step-approx} to handle the terms which are linear in $n$, it follows that
\[ 
|\E_{n \in [N]} f_1(n) e(-T(n,n,n)) f_2(n+h) e(T(n + h, n+h, n+h)) 1_{P \cap B_1}(n) e(\tilde \theta_h n)  e(2\sigma_{n,h,n})| \gg 1, 
\] 
provided that $\eps$ was chosen sufficiently small in terms of $\delta$. This, recall, is for many $h \in P \cap B_1$.

Now from \eqref{approx-sym} and the smallness of $\eps$ we see that $\sigma : (P \cap B_1)^3 \rightarrow \Z/6\Z$ is trilinear. Thus $\sigma_{n,h,n}$ is constant as $n,h$ vary over any translate of $Q := 6 \cdot (P \cap B_1) := \{6x : x \in P \cap B_1\}$. Since $P \cap B_1$ may be covered by $O_{\delta}(1)$ such translates, we may pigeonhole yet again to conclude the existence of $h_2,n_2$ such that 
\[ 
|\E_{n \in [N]} f_1(n + n_2) e(-T(n_2 + n, n_2 + n, n_2 + n)) F_2(n+h)1_{Q}(n + n_2)e(\theta'_h n)| \gg 1.
\]
for many $h$,
where 
\begin{align*} 
F_1(n) & := f_1(n+n_2) e(-T(n_2 +n,n_2 + n, n_2 + n)) \\ & = f(n + n_0 + n_1 + n_2) e(-T(n_2 + n, n_2 + n, n_2 + n) \\ & \qquad + T(h_0, n_0 + n_1 +n_2 + n, n_0 + n_1 + n_2+  n) + T(h_1,n_2 + n,n_2 +n))\end{align*}
and $F_2$ is a $1$-bounded function whose precise nature is unimportant.
It follows from this and an expansion of $1_Q(n+n_2)$ as a Fourier series that
\[ 
\E_h \Vert F_1(n) F_2(n+h) \Vert_{U^2}^4 \gg 1.
\]
Expanding out implies that the Gowers inner product $\langle F_1,F_1,F_1,F_1,F_2,F_2,F_2,F_2\rangle_{U^3}$ is $\gg 1$. By the Gowers-Cauchy-Schwarz inequality we see that $\Vert F_1 \Vert_{U^3} \gg 1$ which, by the inverse theorem for the $U^3$ norm, implies that 
\[ \E_{n \in [N]} F_1(n)\Psi(n) \gg 1\] for some $2$-step nilsequence $\Psi(n)$.

Now $F_1(n)$ is equal to $f(n + n_0 + n_1 + n_2)$ times a variety of bracket terms. By Lemma \ref{23-step-approx}, each of those bracket terms is a product of almost nilsequences of degree at most $3$.  Thus $f$ itself has inner product $\gg_{\delta} 1$ with a degree $3$ almost nilsequence on $[N]$. As we observed in Lemma \ref{almost-to-genuine}, this is enough to establish (at last!) the inverse conjecture for the $U^4$-norm, that is to say Theorem \ref{mainthm}.\endproof

\appendix

\section{Lifting results for nilmanifolds}\label{lift-app}

In this section we establish some slightly technical results concerning the relationship between points on a connected, simply-connected nilpotent Lie group $G$ and points in the nilmanifold $G/\Gamma$. These results were necessary in \S \ref{step1-sec}.

We begin with a folklore result of quantitative linear algebra type.

\begin{lemma}[Bounded equations have bounded solutions]\label{quant-lin-algebra}
Suppose that $A$ is an $m \times n$ matrix and that $b \in \C^m$. Suppose that all of the entries of $A$ are rational numbers of complexity at most $M$, and that the entries of $b$ are bounded by $M$. Then if the equation $Ax = b$ has a solution over $\C^n$, it has a solution in which each coordinate is bounded by $O_{M,m,n}(1)$.
\end{lemma}
\emph{Sketch proof.} By removing rows of $A$ if necessary we may assume that the rows of $A$ are linearly independent. One may then augment $A$ to a nonsingular $n \times n$ matrix $\tilde A$ by adding appropriate basis vectors $e_i$. Augment $b$ to a vector $\tilde b \in \C^n$ by simply adding $n - m$ zeros to $b$. Then the equation $\tilde A \tilde x = \tilde b$ has a solution given by $\tilde x = \tilde A^{-1}\tilde b$. All entries of $\tilde x$ are bounded by $O_{M,m,n}(1)$ by the construction of $\tilde A^{-1}$, the key point here being to note that $|\det \tilde A|$ is bounded below by $\Omega_{M,m,n}(1)$ since it is a nonzero rational number of complexity $O_{M,m,n}(1)$.\endproof

We record the following special case.

\begin{corollary}[Linear lifting]
Suppose that $V \leq \R^n$ is a vector subspace given by the vanishing of linear forms over $\Z$ with coefficients of magnitude at most $M$. Let $\pi : \R^n \rightarrow \R^m$ be projection onto the first $m$ coordinates. Suppose that the entries of $x \in \R^m$ are bounded by $M$, and that $\pi^{-1}(x) \cap V$ is nonempty. Then $\pi^{-1}(x) \cap V$ contains a vector whose entries are bounded by $O_{M,m,n}(1)$. 
\end{corollary}
\proof The condition that a vector $y$ lies in $\pi^{-1}(x) \cap V$ may be encoded as $Ay = b$, where this linear system includes the equations $y_1 = x_1,\dots,y_m = x_m$ and the equations that $y$ must satisfy in order to lie in $V$. By construction the entries of $A$ are rational numbers of complexity at most $M$ and the entries of $b$ are bounded. The corollary therefore follows from the preceding lemma.\endproof

Using a little Lie theory, this last result has the following further corollary.

\begin{corollary}\label{cor-a3}
Suppose that $G$ is a connected, simply-connected nilpotent Lie group and let $\pi : G \rightarrow G/[G,G]$ be the natural projection. Suppose that the Lie algebra $\g = \log G$ has a basis $\mathcal{X} = \{X_1,\dots,X_m,X_{m+1},\dots,X_n\}$, where $\pi (\mathcal{X}) := \{\pi(X_1),\dots,\pi(X_m)\}$ is a basis for $\g/[\g,\g] = \log(G/[G,G])$ as a vector space over $\R$. Suppose that $H$ is an $M$-rational connected subgroup relative to $\mathcal{X}$, and that $\pi(H)$ contains an element $x \in \R^m$ whose entries, written in the basis $\pi(\mathcal{X})$, are bounded by $M$. Then there is an element $\tilde x \in H$ with $\pi(\tilde x) = x$ whose entries are bounded by $O_{M,n}(1)$. 
\end{corollary}
\proof Let $\h = \log H$ be the Lie algebra of $H$. By standard Lie theory (see, for example, \cite{bourbaki}) the exponential/logarithm maps from $\g$ to $G$ and from $\h$ to $H$ are diffeomorphisms. The result now follows from the preceding corollary upon taking $V = \h$. 

This last corollary took place at the level of Lie groups. The actual result we required in \S \ref{step1-sec} concerned lifting from nilmanifolds. We state it now. 

\begin{proposition}[Lifting from nilmanifolds]\label{lift-lemma}
Let $G/\Gamma$ be a nilmanifold with Mal'cev basis $\mathcal{X} = \{X_1,\dots,X_m,X_{m+1},\dots,X_n\}$ and of complexity at most $M$, and let $H \leq G$ be a closed connected $M$-rational subgroup giving rise to a subnilmanifold $H\Gamma/\Gamma$. Then there is a quantity $\eps_M > 0$ with the following property. Suppose that $H[G,G]\Gamma/[G,G]\Gamma$, identified with the torus $\R^m/\Z^m$ using the Mal'cev basis $\mathcal{X}$, contains an element $x$ whose reduced coordinates \textup{(}those nearest 0\textup{)} are all at most $\eps_M$. Let $\psi : G \rightarrow G/[G,G]\Gamma$ be the natural projection onto the horizontal torus of $G/\Gamma$. Then there is a lift $\tilde x \in H$ with coordinates $O_M(1)$ whose first $m$ coordinates are precisely the reduced coordinates of $x$.
\end{proposition}
\proof The Mal'cev coordinates give a commutative diagram
\begin{equation}
\begin{CD} G/[G,G] @>>> \R^m \\ @VVV @VVV \\ G/\Gamma[G,G]
@>>> \R^m/\Z^m.\end{CD}
\end{equation}
The inclusion of $H[G,G]/[G,G]$ into $G/[G,G]$ identifies the former with a vector subspace $V \leq \R^m$ given by the vanishing of linear forms over $\Z$ with coefficients of magnitude $O_M(1)$, and then $H[G,G]\Gamma/[G,G]\Gamma$ becomes identified with $V\Z^m/\Z^m$. Note that this last object is not in general connected, being a union of a finite number of cosets of a subtorus of $\R^m/\Z^m$. We claim that there is an intermediate lift $x'$ of $x$ to $H[G,G]/[G,G]$ whose coordinates in $\R^m$ are the same as the reduced coordinates of $x$ in $\R^m/\Z^m$. Once this claim is proved we may use the last corollary to lift $x'$ again, under the map $\pi : G \rightarrow G/[G,G]$, thereby confirming the proposition.

The claim is a completely abelian statement concerning tori. To prove it, suppose that the linear relations over $\Z$ which define $V$ as a subspace of $\R^m$ are given by $\sum_{j=1}^m k_{ij} x_j = 0$, $i = 1,\dots,m'$. Suppose that $\eps_M < |k_{ij}|/10m$ (say) and that $x$, written as $(x_1,\dots,x_m)$ in reduced coordinates, lies in $V\Z^m/\Z^m$. By assumption we have $|x_j| \leq \eps_M$ for all $j$. Then $\sum_{j=1}^m k_{ij}x_j$ is an integer, yet it also has magnitude at most $1/10$. It must therefore vanish, which means that element $x' \in G/[G,G]$ whose coordinates in $\R^m$ are precisely those of $x$ must lie in $H[G,G]/[G,G]$, as claimed.\endproof

\section{S\'ark\"ozy-type results}

% {\bf to be removed?}

In this section we prove a lemma that was used in the course of the so-called symmetry argument in \S \ref{sym-section}. It is a familiar principle in additive combinatorics that if one takes some fairly ``dense'' set $A$ in an abelian group then the sumsets $2A = A + A$, $3A = A + A + A$ become progressively more structured, containing longer and longer progressions and ever larger Bohr sets. See, for example, \cite{bogolyubov,bourgain-sumsets, green-sumset}. S\'ark\"ozy \cite{sarkozy} was the first to observe that in \emph{very} high-order sumsets $kA$, one may locate very large amounts of structure indeed. The following rather neat version of his result follows directly from a theorem of Lev (\cite[Theorem $2'$]{lev}):

\begin{theorem}[Lev]
Suppose that $A \subseteq [N]$ is a set of size $\alpha N$. Then for any $k \geq 2/\alpha$ the set $kA - kA:= A + \dots + A - A - \dots - A$ contains an arithmetic progression $\{0,d,2d,\dots,(N-1)d\}$ where $d \leq 1/\alpha$. 
\end{theorem}

In \S \ref{sym-section} we required a kind of ``bilinear'' version of this. Suppose that $A \subseteq [N]^2$ is a set. Let us write $A \oplus A$ for the set of all pairs $(x,y_1 \pm y_2)$ where both $(x,y_1)$ and $(x,y_2)$ lie in $A$, together with all pairs $(x_1\pm x_2, y)$ where both $(x_1,y)$ and $(x_2,y)$ lie in $A$. The importance of this definition for us lies in the fact that if a bilinear form is approximately annihilated by $A$ then it is also also approximately annihiliated by $A \oplus A$.

\begin{proposition}[Bilinear S\'ark\"ozy result]\label{bilinear-sarkozy}
Suppose that $A \subseteq [N]^2$ is a set of size $\alpha N^2$. Then for $k \geq 128/\alpha^3$ the $k$-fold iterated bilinear sumset $A \oplus A \dots \oplus A$ contains a product $P \times P'$, where $P = \{0,d,2d,\dots,(N-1)d\}$ and $P' = \{0,d',2d',\dots,(N-1)d'\}$ with $0 < d,d' \leq 4/\alpha^2$.
\end{proposition}
\proof For each $x\in [N]$ write $A_x := \{y \in [N] : (x,y) \in A\}$ for the vertical fibre of $A$ above $x$. By a simple averaging argument there are at least $\alpha N/2$ values of $x$ for which $|A_x| \geq \alpha N/2$. For each such $x$ the vertical sumset $kA_x - kA_x$, where $k \geq 4/\alpha$, contains a progression $P = \{0,d_x,2d_x,\dots,(N-1)d_x\}$ with $0 < d_x \leq 2/\alpha$. By the pigeonhole principle we may pass to a further set $\{A_x : x \in X\}$ of vertical fibres , $|X| \geq \alpha^2 N/4$, which all have the same value of $d_x$, say $d$.  By a further application of Lev's theorem the set $lX - lX$, $l \geq 8/\alpha^2$, contains a progression $P' = \{0,d',2d',\dots, (N-1)d'\}$ with $0 < d' \leq 4/\alpha^2$.  \endproof

\emph{Remark.} We believe that it ought to be possible to prove a structural result in which only some bounded sum $A \oplus A \oplus \dots \oplus A$ is involved, where the number of summands does not depend on $\alpha$ (and might, for example, be 16). Such a result would deserve to be called a ``bilinear Bogolyubov theorem'' by analogy with Bogolyubov's lemma \cite{bogolyubov}. One would not expect to find a structure as simple and rich as the product $P \times P'$; we expect the relevant structure to be, rather, a ``transverse set'', the intersection of sets of the form $\{(x,y) \in [N]^2 : \Vert \theta x y\Vert_{\R/\Z} \leq \eps\}$.

\section{Structure of approximate homomorphisms}
\label{approx-hom-app}

The aim of this appendix is to indicate a proof of Proposition \ref{approx-hom-prop}, whose statement we recall now. As we said before, this result is somehow ``known'' without being explicitly given anywhere in the literature. The forthcoming \emph{Barbados lectures} of the first author will give a self-contained treatment of results of this type. 

\begin{approx-hom-prop-again}[Approximate homomorphisms]%\label{approx-hom-prop}
Let $\delta,\eps \in (0,1)$ be parameters and suppose that $f_1,f_2,f_3,f_4 : S \rightarrow \R/\Z$ are functions defined on some subset $S \subseteq [N]$ such that there are at least $\delta N^3$ quadruples $(x_1,x_2,x_3,x_4) \in S^4$ with $x_1 + x_2 = x_3 + x_4$ and $\Vert f_1(x_1) + f_2(x_2) - f_3(x_3) - f_4(x_4)\Vert_{\R/\Z} \leq \eps$. Then there is a bracket linear phase $\psi : \Z \rightarrow \R/\Z$ of complexity $O_{\delta}(1)$ and a set $S' \subseteq S$, $|S'| \gg_{\delta} N$, such that $f_1(x) = \psi(x) + O(\eps)$ for all $x \in S'$.
\end{approx-hom-prop-again}
\proof We begin with a ``rounding'' trick to dispose of the error of $\eps$ in the range. Take $N := [1/\eps]$ and for $i = 1,2,3,4$ define $\tilde f_i : S \rightarrow \R/\Z$ by taking $\tilde f_i(x) = r/N$, where $r$, $0 \leq r < N$, is the integer such that $r/N$ is nearest to $f_i(x)$ in $\R/\Z$ (ties being broken arbitrarily). Then of course $\tilde f_i(x) = f_i(x) + O(\eps)$ for all $x \in S$ and so
\[  \tilde f_1(x_1) + \tilde f_2(x_2) - \tilde  f_3(x_3) - \tilde f_4(x_4) = O(\eps) \] for the set of additive quadruples $(x_1,x_2,x_3,x_4) \in S^4$ in the hypothesis of the proposition.
The quantity $\Vert \tilde f_1(x_1) + \tilde f_2(x_2) - \tilde f_3(x_3) - \tilde f_4(x_4)\Vert_{\R/\Z}$ is quantised and restricted to integer multiples of $1/N$, and there are only $O(1)$ such numbers with magnitude $O(\eps)$. It follows that there is some $\theta_0$ such that $\Vert \tilde f_1(x_1) + \tilde f_2(x_2) - \tilde f_3(x_3) - \tilde f'_4(x_4)\Vert_{\R/\Z} = 0$ for $c\delta N^3$ additive quadruples $(x_1,x_2,x_3,x_4) \in S^4$, where $\tilde f'_4(x) = \tilde f_4(x) + \theta_0$. 

Writing $\Gamma_i := \{(x,\tilde f_i(x)) : x \in S\} \subseteq \Z \times \R/\Z$, $i = 1,2,3$, and $\Gamma'_4 := \{(x,\tilde f'_4(x)) : x \in S\} \subseteq \R/\Z$ for the ``graphs'' of $\tilde f_1,\tilde f_2,\tilde f_3$ and $\tilde f'_4$, this means that the additive energy (cf. \cite[Chapter 2]{tao-vu}) $E(\Gamma_1,\Gamma_2,\Gamma_3,\Gamma'_4)$ is at least $c\delta N^3$. By \cite[Corollary 2.10]{tao-vu} (or the Cauchy-Schwarz-Gowers inequality) it follows that the additive energy $E(\Gamma_1,\Gamma_1,\Gamma_1,\Gamma_1)$ is at least $c\delta^C N^3$, or in other words that there are $\geq c\delta^C N^3$ additive quadruples $(x_1,x_2,x_3,x_4) \in S^4$ for which $\Vert \tilde f_1(x_1) + \tilde f_1(x_2) - \tilde f_1(x_3) - \tilde f_1(x_4)\Vert_{\R/\Z} = 0$.

From this point on we give references to the paper \cite{green-tao-u3inverse} of the first two authors, which is reasonably well-adapted to our purposes. Most of the ideas here go back to \cite[Chapter 7]{gowers-longaps} and to earlier work of Ruzsa. Starting from the assumption that the graph $\Gamma$ has large additive energy, the key steps are the following\footnote{Strictly speaking, the tools we are applying here only apply to groups rather than to intervals such as $[N]$.  However, this can be easily addressed by temporarily embedding $[N]$ in, say, $\Z/10N\Z$; we omit the details.}. 

\begin{enumerate}
\item \cite[Proposition 5.4]{green-tao-u3inverse} Apply the Balog-Szemer\'edi-Gowers theorem followed by the Pl\"unnecke-Ruzsa inequalities to conclude that there is a set $S_0 \subseteq S$, $|S_0| \geq c\delta^C N$, such that the graph $\Gamma := \{(x,\tilde f_1(x)) : x\in S_0\}$ satisfies an iterative sumset estimate $|k \Gamma - l\Gamma| \ll_{k,l} N$ for all integers $k,l \geq 1$. 
\item \cite[Proposition 9.1]{green-tao-u3inverse} The function $\tilde f_1$ correlates with a function which is locally linear on a Bohr set. This means that there are is a Bohr set $B = B(\Theta,\rho,N)$ with $\Theta = \{\theta_1,\dots,\theta_d\} \subseteq \R/\Z$, $d = O_{\delta}(1)$ and $\rho \gg_{\delta} 1$ together with a function $\phi : B \rightarrow \R/\Z$ satisfying $\phi(x + y) = \phi(x) + \phi(y)$ whenever $x,y,x+y \in B(\Theta,\rho,N)$, as well as some $x_0\in [N]$ and some $\theta_0 \in \R/\Z$ such that $\tilde f(x + x_0) = \theta_0 + \phi(x)$ for $\gg_{\delta} N$ values of $x \in (S_0 - x_0) \cap B$. The appropriate definitions here are given in full in \cite{green-tao-u3inverse} and are also recalled in \S \ref{sym-section} of the present paper.
\item Apply some geometry of numbers to conclude that any such linear function $\phi$ has the form $\phi(x) = \alpha_1 \{\theta_1 x\} + \dots + \alpha_d \{\theta_d x\} + \theta x$ on some multidimensional progression $P \subseteq B$ with $|P| \gg_{\delta} N$. The proof of this is very similar to, but easier than, that of \cite[Proposition 10.8]{green-tao-u3inverse}, where an analogous statement is established for locally \emph{quadratic} phase functions on Bohr sets. 
\end{enumerate}

It follows from all of this that we have

\[ \tilde f_1(x) = \theta_0 + \alpha_1 \{\theta_1 (x - x_0)\} + \dots + \alpha_d \{\theta_d (x - x_0)\} + \theta(x - x_0)\] for all $x$ in some set $S_1 \subseteq S_0$, $|S_1| \gg_{\delta} N$. 

Now we have $\{\theta_j (x - x_0)\} = \{\theta_j x\} - \{\theta_j x_0\} + \tau_{j,x}$, where $\tau_{j,x}$ takes values in $\{-1,0,1\}$. By the pigeonhole principle we may pass to a further subset $S_2 \subseteq S_1$ with $|S_2| \gg_{\delta} N$ such that, for all $x \in S_2$, each of the $\tau_{j,x}$ is independent of $x$.

Take $S' := S_2$. Then for $x \in S'$ we have
\[ \tilde f_1(x) = \theta'_0 + \alpha_1 \{\theta_1 x\} + \dots + \alpha_d \{\theta_d x\} + \theta x,\] a bracket linear form of complexity $d = O_{\delta}(1)$. Recalling that $\tilde f_1(x) = f_1(x) + O(\eps)$, the result follows.\endproof

\emph{Remark.} The rounding trick we used to remove the $\eps$ errors was a slightly dirty one but makes the argument quite short given known results. It would probably be possible, and more natural in some moral sense, to run through the Balog-Szemer\'edi-Gowers and Freiman arguments carrying an $O(\eps)$ error throughout.\vspace{11pt}

%{\bf Perhaps we should flesh this out a little, as it is the key place in the argument where bracket behaviour actually gets imposed. I didn't quite have the energy.}
 
\section{Some diophantine results}

This section recalls some well-known results from Diophantine approximation which, in the context of this paper, may be naturally viewed as distributional results for abelian (1-step) nilsequences. We will use them repeatedly in the next section. Furthermore Lemma \ref{kronecker} below was crucial in \S \ref{step1-sec}, and Lemma \ref{weyl-application} was required at the end of \S \ref{2-step-sec}.

\begin{lemma}\label{weyl-application}
Let $d \geq 1$ be an integer, let $\eps \in (0,1/2)$ be a parameter, and suppose that $\psi(n) = \alpha_d n^d + \dots + \alpha_0$ is a polynomial of degree $d$ such that $(\psi(n)\mdlem{1})_{n \in [N]}$ is not $\eps$-equidistributed on $\R/\Z$. Then for all $i = d,d-1,\dots,1$ there are coprime integers $a_i,q_i$,  $q_i \leq \eps^{-C_d}$, such that 
\[ \alpha_i = \frac{a_i}{q_i} + O(\frac{\eps^{-C_d}}{N^i}).\]
\end{lemma}
\proof This is actually a special case of the Quantitative Leibman Dichotomy, Theorem \ref{quant-leib}, although this is a somewhat misleading statement to make since it is also a crucial ingredient in the \emph{proof} of that result. It is proven using Weyl's criterion for equidistribution and Weyl's inequality (see, for example, \cite{vaughan}), and indeed the statement that the lead coefficient $\alpha_d$ is close to rational is essentially equivalent to that inequality. The other coefficients $\alpha_{d-1},\alpha_{d-2}\dots$ may be shown to be almost rational iteratively; the argument is given in detail in \cite[\S 4]{green-tao-nilratner}.\endproof

Secondly we recall a quantitative version of Kronecker's theorem, phrased in language appropriate to \S \ref{step1-sec}. Once again this is a special case of the Quantitative Leibman Dichotomy, and once again it is very well-known.

\begin{lemma}\label{kronecker}
Let $d \geq 1$ be an integer, let $\eps \in (0,1/2)$ be a parameter, and let $\alpha_1,\dots,\alpha_d \in \R/\Z$ be frequencies. Suppose that $((\alpha_1 n,\dots,\alpha_d n)\mdlem{1})_{n \in [N]}$ fails to be $\eps$-equidistributed in the torus $(\R/\Z)^d$. Then the set $\{\alpha_1,\dots,\alpha_d\}$ satisfies an $\eps^{-C_d}$-linear relation up to $\eps^{-C_d}/N$ \textup{(}that is, there are integers $m_1,\dots,m_d$, not all zero, with $|m_i| \leq \eps^{-C_d}$ for all $i$ and $\Vert m_1 \alpha_1 + \dots + m_d \alpha_d\Vert_{\R/\Z} \leq \eps^{-C_d}/N$\textup{)}.
\end{lemma}
\proof This is discussed in detail in \cite[\S 3]{green-tao-nilratner}. Here is a very rough sketch: if the sequence is not $\eps$-equidistributed, there is some Lipschitz function $F : (\R/\Z)^d \rightarrow \C$ with \[ |\E_{n \in [N]} F(\alpha_1 n,\dots,\alpha_d n) - \int_{(\R/\Z)^d} F(\theta) d\theta| \geq \eps \Vert F \Vert_{\Lip}.\] Lipschitz functions are well-approximated in $L^{\infty}$ by their Fourier series; exanding $F$ into such a series, it follows that some exponential sum
\[ \E_{n \in [N]} e(\vec{m} \cdot \vec{\alpha} n)\] has modulus at least $\eps^{C_d}$, where $\vec{m} = (m_1,\dots,m_d)$ and $|m_i| \leq \eps^{-C_d}$. The lemma now follows with an application of the formula for the sum of a geometric series.\endproof

\section{Almost  nilsequences}
\label{approx-nil-app}

The aim of this section is to establish Lemmas \ref{1-step-approx} and \ref{23-step-approx}, which asserted that various objects -- chiefly bracket polynomials -- are $1$-, $2$- and $3$-step almost nilsequences.

\begin{1step-repeat}Suppose that $\alpha,\beta \in [0,1]$ and that $M > 1$ is a complexity parameter. The following are all examples of almost nilsequences of degree $1$ and complexity $O_M(1)$:
\begin{enumerate}
\item  the set of $1$-step Lipschitz nilsequences of complexity at most $M$;
\item  the set of characteristic functions $1_P$, where $P \subseteq [N]$ is a progression of length at least $N/M$;
\item  the set of functions of the form $n \mapsto e(\alpha \{\beta n\})$, with $\alpha \in \R$ and $\beta \in \R/\Z$;
\item  the set of functions of the form $n \mapsto e(\{\alpha n\} \{\beta n\})$, with $\alpha,\beta \in \R/\Z$;
\item  the set of functions of the form $n \mapsto e(\alpha n\lfloor \beta n \rfloor)$, where $\Vert \beta \Vert_{\R/\Z} \leq M/N$.
\end{enumerate}
\end{1step-repeat}
\proof (i) is trivial.

To prove (ii) we first note that $1_P(n)$ can be expressed as the product of $1_I(n)$ and $1_{n \equiv a \mdsub{q}}$, where $I \subseteq N/M$ is an interval and $q \leq M$. The second object is in fact a 1-step nilsequence $F(g(n)\Gamma)$ on $\R/\Z$, the polynomial sequence $g : \Z \rightarrow \R$ being $g(n) = n/q$ and the function $F : \R/\Z \rightarrow [0,1]$ being Lipschitz, equal to $1$ at $a/q$ and supported within $1/10q$ (say) of $a/q$. The first object, $1_I(n)$, is not quite a genuine 1-step nilsequence. However let us observe that any function $\psi : [N] \rightarrow \C$ with Lipschitz constant $O(1/N)$ is a genuine $1$-step nilsequence; indeed we have $\psi(n) = F(g(n)\Gamma)$ on $\R/\Z$, where $g(n) = n/2N$ and $F : \R/\Z \rightarrow \C$ is defined by setting $F(n/2N) := \psi(n)$ for $n \in [N]$ and by Lipschitz extension elsewhere. Now simply note that $1_I$ may be approximated arbitrarily closely, in $L^1[N]$, by functions $\psi$ of this type. Specifically, we may take a sequence of Lipschitz ``tent'' functions $\psi$ which equal $1$ on $I$ and are zero at points distance more than $\eps N$ from $I$.  The claim now follows from Lemma \ref{alg-lem}. 

To establish (iii) we first note that if $\alpha \equiv \alpha' \md{1}$ then $\alpha \{ \beta n\} \equiv \alpha' \{\beta n\} + (\alpha - \alpha')\beta n \md{1}$, and so we may assume that $0 \leq \alpha \leq 1$. Let $\eps > 0$ be arbitrary and define $F : \R/\Z \rightarrow \C$ by $F(x) = e(\alpha\{x\})$ and divide into two cases: either $(\beta n \md{1})_{n \in [N]}$ is $\eps/10$-equidistributed on $\R/\Z$, or it is not. In the former case we take a $100\eps$-Lipschitz function $\tilde F$ which agrees with $F$ outside of the set $\{ x \in \R/\Z : \Vert x \Vert_{\R/\Z} \leq \eps/10\}$ and is bounded by 1 elsewhere. By the assumed equidistribution we obviously have $e(\alpha \{\beta n\}) = F(\beta n) = \tilde F(\beta n)$ for all except at most $\eps N/2$ values of $n$. The result is then immediate.

If, on the other hand, the sequence $(\beta n \md{1})_{n \in [N]}$ fails to be $\eps/10$-equidistributed then by Lemma \ref{weyl-application} with $d = 1$ there is an integer $q \leq \eps^{-C}$ and an $a \in \Z$ such that $\Vert \beta - \frac{a}{q}\Vert_{\R/\Z} \ll \eps^{-C} /N$. This in turn means that we may divide $[N]$ into progressions $P_1 \cup \dots \cup P_m$, $m \ll \eps^{-C}$, on which $n \mapsto e(\alpha \{\beta n\})$ varies by at most $\eps/100$. Since (by part (ii)) functions which are constant on progressions are almost 1-step nilsequences, the result follows (using Lemma \ref{alg-lem} as necessary).

To prove (iv) we use a trick. The function $(x,y) \mapsto e(xy)$ on the square $[0,1]^2$ may be smoothly extended to a periodic function on $[0,2]^2$. By Fourier analysis (cf. \cite[Lemma A.9]{green-tao-u3mobius}) it may then be uniformly approximated to any desired accuracy $\eps$ by a linear combination of frequencies $e((kx + ly)/2)$, $k,l \in \Z$. Thus $e(\{\alpha n\} \{\beta n\})$ may be approximated uniformly by a linear combination of functions of the form $e(k\{\alpha n\}/2)e(l\{\beta n\}/2)$. But such functions are almost 1-step nilsequences by (iii), and the claim follows from Lemma \ref{alg-lem}.

Finally we turn to (v). The condition that $\Vert \beta \Vert_{\R/\Z} \leq M/N$ means that we may divide $[N]$ into subprogressions (in fact subintervals) $P_1 \cup \dots \cup P_m$, $m = O_M(1)$, such that $\lfloor \beta n \rfloor$ is equal to some constant $c_j$ for $n \in P_j$.  The result then follows from (ii) and Lemma \ref{alg-lem}.\endproof

Now we turn to higher degree bracket polynomial phases.

\begin{23step-repeat} Suppose that $\alpha,\beta,\gamma \in [0,1]$. Then the following are all examples of almost nilsequences of degree $s \geq 2$ and complexity $O(1)$:
\begin{enumerate}
\item $n \mapsto e(\lfloor \alpha n \rfloor\beta n)$, of degree $2$;
\item $n \mapsto e(\lfloor \alpha n \rfloor\beta n^2)$, of degree $3$;
\item $n \mapsto e(\lfloor \alpha n \rfloor \lfloor \beta n \rfloor\gamma n)$, of degree $3$.
\end{enumerate}
\end{23step-repeat}
\proof The proofs of all three parts are somewhat similar and proceed along the following lines: each object may be exhibited in a fairly obvious way as a nilsequence $F(g(n)\Gamma)$, where $F$ is, however, only piecewise Lipschitz. If the sequence $(g(n)\Gamma)_{n \in [N]}$ is highly equidistributed then it spends sufficient time away from singularities for one to be able to approximate by $\tilde F(g(n)\Gamma)$, where $\tilde F$ is genuinely Lipschitz. If not then there must be an approximate rational relation between the horizontal frequencies of $g(n)$ (that is, the frequencies occurring in the projection to $G/[G,G]$). This may then be used to approximate the object in question by objects of lower complexity.

To exhibit these arguments as part of a more general theory is not a particularly easy matter and involves a more conceptual understanding of bracket identities such as those in Lemma \ref{brack-identities} and others such as \eqref{3-brack} below. The required theory is implicit in the work of Leibman \cite{leibman} and will be introduced properly in our longer paper to come. 

In this paper we can proceed in an \emph{ad hoc} and slightly calculational way, taking advantage of one or two simplifications specific to the $U^4$ (3-step) case. In a sense, however, these calculations also serve as motivation for the longer paper to come. We begin by recalling the constructions of \S \ref{sec3} leading up to \eqref{5.33}. Specialising to the free 2-step nilpotent group on two generators (essentially the Heisenberg group) we have
\[ F_{[1,2]}(g(n)\Gamma) = e(\lfloor \alpha n \rfloor\beta n) \]
and 
\[ F_{[1,2]}(g'(n)\Gamma) = e(\lfloor \alpha n \rfloor \beta n^2)\]
where $F_{[1,2]} : G/\Gamma \rightarrow \C$ is the basic coordinate function introduced in Definition \ref{def5.3} and $g,g' : \Z \rightarrow G$ are polynomial sequences of degree $2$ and $3$ respectively given in coordinates by $g(n) = (\alpha n,-\beta n,0)$, $g'(n) = (\alpha n, -\beta n^2,0)$. Only the first two coordinates (corresponding to the horizontal torus $G/[G,G]$) are really important.

The discontinuities of $F_{[1,2]}$ are very manageable: the key point, already exploited in \S \ref{step1-sec}, is that for any $\eps > 0$ there is are $\eps^{-C}$-Lipschitz functions $\tilde F : G/\Gamma \rightarrow \C$ and $\Psi : G/\Gamma \rightarrow [0,1]$ such that $\int_{G/\Gamma} \Psi  \leq \eps$ and $|F(x) - \tilde F(x)| \leq \Psi(x)$ pointwise. 

Fix $\eps > 0$. Let us consider statement (i), for which we consider the sequence $(g(n)\Gamma)_{n \in [N]}$. If it is $\eps$-equidistributed then, by the preceding, $e(\{\alpha n\}\beta n) = F(g(n)\Gamma)$ and $\tilde F(g(n)\Gamma)$ are within $2\eps$ in $L^1[N]$. If this is not the case then, by the Quantitative Leibman Dichotomy (Theorem \ref{quant-leib}) there must be some $O(\eps^{-C})$-linear relation, up to $O(\eps^{-C}/N)$, between $\alpha$ and $\beta$. The rest of the argument in this case is essentially identical to that at the very end of \S \ref{step1-sec}; we may find some $\gamma$ such that $\alpha = q_1 \gamma + O(\eps^{-C}/N)$ and $\beta = q_2 \gamma + O(\eps^{-C}/N)$, where $q_1,q_2$ are integers with magnitude at most $\eps^{-C}$. Substituting into $e(\alpha n\lfloor \beta n \rfloor)$ and making repeated use of the bracket identities of Lemma \ref{brack-identities} as well as Lemma \ref{1-step-approx}, one sees that in this case $e(\alpha n\lfloor \beta n \rfloor)$ lies within $\eps$ in $L^1[N]$ of a degree $2$ nilsequence (of step 1) of complexity $O_{\eps}(1)$. Thus in either case we have approximated $e(\alpha n\lfloor \beta n \rfloor)$ within $O(\eps)$ by a degree $2$ polynomial nilsequence of complexity $O_{\eps}(1)$, thereby completing the proof of (i). 

The analysis of (ii) is similar but, obviously, involves consideration of the sequence $(g'(n)\Gamma)_{n \in [N]}$ instead. If the sequence $(g'(n)\Gamma)_{n \in [N]}$ is $\eps$-equidistributed then we are done, as before. If not, the Quantitative Leibman Dichotomy implies that either $\alpha = \frac{a_1}{q_1} + O(\eps^{-C}/N)$ or else $\beta = \frac{a_2}{q_2} + O(\eps^{-C}/N^2)$. In the first case we may then partition $[N]$ into progressions $P_1 \cup \dots \cup P_m$, $m = O_{\eps}(1)$, on which $\lfloor \alpha n \rfloor$ is constant and then apply Lemma \ref{1-step-approx} (ii) to approximate $e(\lfloor \alpha n \rfloor \beta n^2)$ within $O(\eps)$ by a degree $2$ polynomial nilsequence of complexity $O_{\eps}(1)$. In the second case we first apply the bracket identity \eqref{key-bracket} to write 
\begin{equation}\label{another-brack} e(\lfloor \alpha n \rfloor\beta n^2) = e(\alpha \beta n^3) e(-\alpha n[\beta n^2]) e(-\{\alpha n\}\{\beta n^2\}).\end{equation}
The first term here is already a degree $3$ polynomial nilsequence of complexity $O(1)$. In the second term we may partition $[N]$ into progressions $P_1 \cup \dots \cup P_m$, $m = O_{\eps}(1)$, on which $[\beta n^2]$ is constant and then apply Lemma \ref{1-step-approx} (ii) to approximate arbitrarily closely by a degree $1$, nilsequence. The third term, $e(-\{\alpha n\} \{\beta n^2\})$, may be handled using the same trick as in the proof of Lemma \ref{1-step-approx} (iv). This reduces matters to handling $e(\theta \{ \theta' n\})$ (already known to be a degree $1$ almost nilsequence by Lemma \ref{1-step-approx} (iii)) and $e(\theta \{ \theta' n^2\})$. By an argument almost identical to that used in the proof of Lemma \ref{1-step-approx} (iii), only using Lemma \ref{weyl-application} with $d = 2$ instead, this second object may be shown to be an degree $2$ almost nilsequence.  Using Lemma \ref{alg-lem} to put everything together, we obtain the claim.

We turn now to the proof of (iii), which is important in the sense that it is the only place in our paper where a 3-step nilmanifold is actually constructed!

Specifically, we let ${\mathfrak g}$ be the free $3$-step Lie algebra generated by three generators $e_1,e_2,e_3$, or equivalently
\begin{align*}
G:= & \{e_1^{t_1}e_2^{t_2}e_3^{t_3} \\ & e_{21}^{t_{21}}e_{211}^{t_{211}}e_{31}^{t_{31}}e_{311}^{t_{311}}
e_{32}^{t_{32}}e_{322}^{t_{322}}e_{212}^{t_{212}}e_{312}^{t_{312}}e_{213}^{t_{213}}e_{313}^{t_{313}}e_{323}^{t_{323}}: t_i,t_{ij},t_{ijk} \in \R, 1\leq i,j,k,\leq 3\}.
\end{align*}
subject to the relations $e_i^{-1}e_j^{-1} e_ie_j = [e_i, e_j] = e_{[i,j]}$ for $1\leq j < i \leq 3$, 
$[[e_i,e_j],e_k]=e_{ijk}$, and the Jacobi relation $[[e_i,e_j],e_k][[e_j,e_k],e_i][[e_k,e_i],e_j]=1$. 
Inside $G$ we take the standard lattice 
\begin{align*}
\Gamma := &
\{e_1^{n_1}e_2^{n_2}e_3^{n_3} \\ &e_{21}^{n_{21}}e_{211}^{n_{211}}e_{31}^{n_{31}}
e_{311}^{n_{311}}
e_{32}^{n_{32}}e_{322}^{n_{322}}e_{212}^{n_{212}}e_{312}^{n_{312}}e_{213}^{n_{213}}e_{313}^{n_{313}}e_{323}^{n_{323}}: n_i,n_{ij},n_{ijk} \in \Z, 1\leq i,j,k,\leq 3\}.
\end{align*}
Then $G/\Gamma$ is the free $3$-step nilmanifold on $3$ generators.   We take $G_\bullet$ to be the lower central series on $G$

We abbreviate $e_1^{t_1} \ldots e_{323}^{t_{323}}$ as $(t_1,\ldots,t_{323})$.  A computation yields the multiplication law
\[
(t_1,\ldots,t_{323})\star (u_1,\ldots,u_{323})=(s_1,\ldots,s_{323})
\]
where $s_i=t_i+u_i$ for $i=1,2,3$, $s_{ij}=t_{ij}+u_{ij}+t_iu_j$ for $1 \leq j <i \leq 3$, and 
$s_{312}=t_{312}+u_{312}+t_{32}u_1+t_{31}u_2+t_3u_1u_2$; we will ignore the other coordinates, as
they will not be needed in this calculation. 

Using this law, we see that for any real numbers $t_1,\ldots,t_{323}$, one has
$$
(t_1,\ldots,t_{323})\Gamma=(s_1,\ldots,s_{323})\Gamma$$
where
\begin{align}
\nonumber s_i &:= \{ t_i \} \hbox{ for } i=1,2,3; \\
\nonumber s_{ij} &:= \{ t_{ij}- t_{i}[t_j]\} \hbox{ for } 1 \leq j<i \leq 3; \\
s_{312} &:= \{t_{312}-t_{32} [t_1]-t_{31}[t_2] +t_3 [t_1][t_2]\},\label{coords-comp}
\end{align}
and with the other coordinates $s_{ijk} \in [0,1]$ being explicitly computable, but not relevant for this discussion.  Thus if we let
$$ g(n) := e_1^{\alpha n} e_2^{\beta n} e_3^{\gamma n}$$
and let $F: G/\Gamma \to \C$ be the $3$-step basic coordinate function function
$$ F((s_1,\ldots,s_{323})\Gamma) := e( s_{312} )$$
for $s_1,\ldots,s_{323} \in [0,1]$, then one sees that $e(\lfloor \alpha n \rfloor\lfloor \beta n \rfloor\gamma n)$ is equal to $F(g(n) \Gamma )$ times objects already known to be almost nilsequences by earlier parts.

This concludes the argument unless $(g(n)\Gamma)_{n \in [N]}$ spends too much time near the singularities of $F$, which are at the points $s_j = 0$ and $s_j = 1$, $j = 1,2,3$. There will be no problem unless\footnote{This observation, which is stronger than saying that the abelianization $((\alpha n,\beta n,\gamma n) \md{1})_{n \in [N]}$ $\subseteq (\R/\Z)^3$ is not equidistributed, is somewhat specific to the $3$-step situation we are working with and represents something of a simplification over the argument required in general.} one of the sequences $(\alpha n \md{1})_{n \in [N]}$, $(\beta n \md{1})_{n \in [N]}$, $(\gamma n \md{1})_{n \in [N]}$ fails to be $\eps$-equidistributed. If $(\alpha n \md{1})_{n \in [N]}$ is not $\eps$-equidist- ributed then, by the now-familiar application of Lemma \ref{weyl-application} with $d = 1$, we may partition $[N]$ as a union $P_1 \cup \dots \cup P_m$ of at most $\eps^{-C}$ progressions such that $\lfloor \alpha n \rfloor$ is constant on $P_i$. We may then conclude using part (i) and Lemma \ref{1-step-approx} (ii).
An identical argument works if $(\beta n \md{1})_{n \in [N]}$ fails to be $\eps$-equidistributed.

The final case is when $(\gamma n \md{1})_{n \in [N]}$ fails to be $\eps$-equidistributed. In this case we note that 
\begin{equation}\label{3-brack} \{\alpha n\}\{\beta n\} \{\gamma n\} = (\alpha n - \lfloor \alpha n \rfloor)(\beta n - \lfloor \beta n \rfloor)(\gamma n - \lfloor \gamma n\rfloor)\end{equation} so that 
\begin{align*} e(\lfloor \alpha n \rfloor\lfloor \beta n \rfloor\gamma n)  & = e(\{\alpha n\} \{\beta n\}\{\gamma n\})e(-\lfloor \alpha n \rfloor\beta n\lfloor\gamma n\rfloor) \times \\ & \times e(-\alpha n\lfloor \beta n \rfloor\lfloor \gamma n\rfloor)e(\alpha\beta n^2\lfloor \gamma n\rfloor) e(\alpha \gamma n^2\lfloor \beta n \rfloor) e(\beta \gamma n^2 \lfloor \alpha n \rfloor).\end{align*}
Each of the terms on the right except the first can be handled using part (ii) or by those instances of part (iii) already established. To deal with the first term $e(\{\alpha n\}\{\beta n\} \{\gamma n\})$ one may proceed exactly as in Lemma \ref{1-step-approx} (iv) to show that this is in fact an degree $1$ almost nilsequence.  Applying Lemma \ref{alg-lem} to collect terms, we obtain the claim.\endproof

The main business of the paper is now concluded. The remaining two appendices were promised in the introduction but are not necssary for the proof of Theorem \ref{mainthm}.

\section{The strong inverse conjecture}\label{strong-appendix}

We have shown, in Theorem \ref{mainthm}, that a $1$-bounded function $f : [N] \rightarrow \C$ with $\Vert f \Vert_{U^4} \geq \delta$ correlates with a degree $3$ polynomial nilsequence $F(g(n)\Gamma)$. As we remarked after the statement of Conjecture \ref{gis-conj}, this does not quite establish the result used in (for example) \cite{green-tao-linearprimes}, where correlation with a nilsequence $F(g^nx\Gamma)$ was used. In this section we shall refer to \emph{linear nilsequences} to distinguish objects of this last type from more general polynomial nilsequences.

In this section we indicate, very briefly, how our arguments may be modified to obtain this apparently stronger statement. In the longer paper to come we will provide a quite general proof that Conjecture \ref{gis-conj} implies this strong variant. Let us recall once more, however, our view that this is the ``wrong'' perspective and that \cite{green-tao-linearprimes} works, with rather minimal changes, in the context of polynomial nilsequences.

We need only show that large $U^4$-norm entails correlation with \emph{almost} linear nilsequences, defined in exact analogy with Definition \ref{almost-nil-def}.  We already have correlation with almost polynomial sequences, so it will suffice to show that the almost polynomial sequences used in the paper are also almost linear sequences of the same degree. 

Clearly, any degree $1$ almost nilsequence is already an almost linear $1$-step nilsequence, and an inspection of the previous appendix shows that $e(\lfloor \alpha n \rfloor \beta n)$ is an almost linear $2$-step nilsequence, modulo a quadratic phase $e( \gamma n^2 )$, and similarly $e(\lfloor \alpha n \rfloor \lfloor \beta n\rfloor \gamma n)$ is an almost linear $3$-step nilsequence modulo phases such as $e( \gamma n^3 )$ and $e( \lfloor \gamma n \rfloor \delta n^2 )$.  As Lemma \ref{alg-lem} is clearly also valid for almost linear nilsequences, one only needs to verify three remaining claims, for any real numbers $\alpha,\beta$:
\begin{itemize}
\item $e(\alpha n^2)$ is an almost linear $2$-step nilsequence;
\item $e(\alpha n^3)$ is an almost linear $3$-step nilsequence; and
\item $e(\lfloor \alpha n\rfloor \beta n^2)$ is an almost linear $3$-step nilsequence.
\end{itemize}

We look first at $e(\alpha n^2)$ and consider once again the 2-step nilpotent group on $2$ generators (Heisenberg group); looking all the way back to \eqref{5.33} and taking $g = (2\alpha,1,0)$ one may compute that $F_{[1,2]}(g^n\Gamma) = e(\alpha n^2 + \theta n)$ for some $\theta \in \R/\Z$. Now $F_{[1,2]}(t_1,t_2,t_{12})$ is discontinuous when $t_{12} = 0$ or $1$. If we wish to approximate $e(\alpha n^2 + \theta n)$ within $\eps$ (in $L^1[N]$) by a \emph{Lipschitz} linear nilsequence, we must show (for example) that there are no more than $10\eps N$ values of $n \in [N]$ for which $\alpha n^2 + \theta n \md{1}$ is within $\eps$ of $0$. But if this is not the case then, by Lemma \ref{weyl-application}, we have $\alpha = a/q + O(\eps^{-C}/N^2)$, at which point we can split $[N]$ into $\eps^{-C}$ progressions on which $e(\alpha n^2 + \theta n)$ is within $O(\eps)$ of a linear phase. One may then proceed using Lemma \ref{1-step-approx}.

Now we turn to the $3$-step objects $e(\alpha n^3)$ and $e(\lfloor \alpha n \rfloor\beta n^2)$, which require some slightly more careful calculations on the free $3$-step nilmanifold are required. With the notation for the free $3$-step nilpotent Lie group as in the preceding section, let $g=e_1^{\alpha}e_2^{\beta}e_3^{\gamma}$.
Then one can check that 
\[ g^n=e_1^{n\alpha} e_2^{n\beta} e_3^{n\gamma}  e_{21}^{\binom{n}{2}\alpha\beta} 
e_{31}^{\binom{n}{2}\alpha\gamma} e_{32}^{\binom{n}{2}\beta\gamma} \cdots e_{312}^{\alpha \beta \gamma (2\binom{n}{3}+\binom{n}{2})} \cdots\]
and hence one may compute (cf. \eqref{coords-comp})
\begin{equation}\label{general-bracket}
F_{312}(g^n \Gamma)= \textstyle e\big(\alpha \beta \gamma \left(2\binom{n}{3}+\binom{n}{2}\right) - \binom{n}{2} \beta\gamma\lfloor \alpha n \rfloor- \binom{n}{2} \alpha \gamma \lfloor \beta n \rfloor + n\gamma[\alpha n ]\lfloor \beta n \rfloor \big)\end{equation}
Taking $\beta = \gamma = 1$ and replacing $\alpha$ by $6\alpha$ gives $F_{312}(g^n\Gamma) = e(\alpha n^3+q(n))$ for some quadratic $q$. The discontinuities of $F_{312}$ may be handled as for $F_{[1,2]}$ above, and so we see that $e(\alpha n^3 + q(n))$ is an almost $3$-step linear nilsequence for some quadratic $q$. Since we can already obtain pure quadratic and linear phases as almost linear nilsequences of step less than $3$, it follows that $e(\alpha n^3)$ itself is an almost $3$-step linear nilsequence.

Next, we take $\beta = 1$ and replace $\gamma$ by $-2\gamma$. Taking into account objects already known to be almost linear nilsequences, we have now obtained $e(\gamma n^2 [n\alpha])$ as a $3$-step almost linear nilsequence. Applying \eqref{another-brack}, we see that to obtain the desired object $e([n\alpha]\gamma n^2)$ it suffices to examine $e(\{\alpha n\}\{\gamma n^2\})$. By the trick used in the proof of Lemma \ref{1-step-approx} (iv), it suffices in turn to handle $e(\theta \{\theta' n\})$ and $e(\theta \{\theta' n^2\})$. The first of these is an almost 1-step (linear) nilsequence by Lemma \ref{1-step-approx} (iii). To handle the second, proceed in the same way as in the proof of Lemma \ref{1-step-approx} (iii) but in the obvious places substitute the fact (established above of course) that pure quadratic phases are $2$-step linear nilsequences, together with the case $d = 2$ of Lemma \ref{weyl-application}.\endproof

\section{Necessity of the inverse conjectures}\label{nec-appendix}

In this appendix we sketch a rather short proof of Proposition \ref{inv-nec}, which asserted that functions which correlate with a degree $s$ polynomial nilsequence must have large $U^{s+1}$-norm. Since \emph{linear} nilsequences are merely special cases of polynomial ones, this kind of argument could substitute in, for example, \cite[Sec. 10]{green-tao-linearprimes}, where a rather more complicated approach was taken. 

\begin{inv-nec-again}
Suppose that $f :[N] \rightarrow \C$ is a $1$-bounded function, that $(F(g(n)\Gamma))_{n \in \Z}$ is a polynomial nilsequence of degree $s$ and complexity $O_{\delta}(1)$, and that \[ |\E_{n \in [N]} f(n) \overline{F(g(n)\Gamma)}| \geq \delta.\] Then $\Vert f \Vert_{U^{s+1}} \gg_{\delta} 1$.
\end{inv-nec-again}
\emph{Sketch proof.} The argument is only a sketch in that we do not address such issues as the complexity of the nilsequences involved. We leave this as a (not particularly interesting) exercise to the reader, most of the details of which may be found in \cite{green-tao-nilratner} where these complexity issues are discussed in detail. We proceed by induction on $s$, the claim being obvious when $s = 0$.  
Let  $f:[N] \to \C$ be a $1$-bounded function, and let $g: \Z \to G$ be a polynomial sequence 
of degree $s$ adapted to the filtration $G_{\bullet}$.  Let $F(g(n)\Gamma)$ be a polynomial nilsequence
of complexity $O_{\delta}(1)$.  Assume that 
\begin{equation}\label{to-square}|\E_{n \in [N]} f(n)F(g(n)\Gamma)| \gg_{\delta} 1.\end{equation} By decompositing $F$ into vertical characters as in \cite[Lemma 3.7]{green-tao-nilratner}, we may assume that $F$ has a vertical frequency: that is, there is some nontrivial character $\xi : G_{(s)}/G_{(s)} \cap \Gamma$ such that  
\[ F(g_s x) = e(\xi(g_s))F(x)\] for all $g_s \in G_{(s)}$ and $x \in G/\Gamma$.

By taking the modulus squared of \eqref{to-square} and making the substitution $n' = n + h$ we see that 
\[
 \E_{n \in [N]} \Delta_h f (n) F(g(n+h)\Gamma) \overline{F(g(n)\Gamma)} \gg_{\delta} 1
\]
for $\gg_{\delta} N$ values of $h \in [N]$.

However for each fixed $h$ the ``derivative'' $n \mapsto F(g(n+h)\Gamma) \overline{F(g(n)\Gamma)}$ of the degree $s$-step nilsequence $F(g(n)\Gamma)$ is a Lipschitz polynomial nilsequence of degree $(s-1)$, the underlying nilmanifold being
\[ \overline{(G^{\square})} = (G \times_{G_{(2)}} G)/G^{\triangle}_{(s)},\] where 
$G \times_{G_{(2)}} G=\{(g,h):g,h \in G, gh^{-1} \in G_{(2)}\}$, and $G^{\triangle}_s=\{(g_s,g_s): g_s \in G_{(s)}\}$. For details of this theory see Section 7 of \cite{green-tao-nilratner}.

We now invoke our induction hypothesis to conclude that  
\[ \Vert \Delta_h f \Vert_{U_{s}} \gg_{\delta} 1\] for $\gg_{\delta} N$ values of $h$.

Noting that 
\[ \Vert f \Vert^{2^{s+1}}_{U_{s+1}} = \E_{h \in \Z/N'\Z} \|\Delta_h f\|^{2^{s}}_{U_{s}},\] we are done.  \endproof
 
It is perhaps worth reiterating the main point of the above argument, since it explains the importance of nilsequences in the whole theory: the derivative of a degree $s$ polynomial nilsequence with a vertical character is a degree $(s-1)$ polynomial nilsequence.

\providecommand{\bysame}{\leavevmode\hbox to3em{\hrulefill}\thinspace}

\end{document}